\documentclass[11pt,leqno]{amsart}

\usepackage{parskip}
\allowdisplaybreaks
\newcommand{\res}{\mathbin{\hspace{0.1em}\vrule height 1.3ex depth 0pt width 0.13ex\vrule height 0.13ex depth 0pt width 1.0ex}} 

\usepackage[left = 2.54cm, right = 2.54cm, top = 2.54cm, bottom = 2.54cm]{geometry}
\usepackage{amsmath}
\usepackage{amsfonts}
\usepackage{amssymb}
\usepackage{amsthm}
\usepackage{mathtools}
\usepackage{graphicx,caption, color, hyperref}
\usepackage{enumitem}
\usepackage{multicol}
\usepackage{bbm}
\usepackage{stmaryrd}
\usepackage{tikz}
\usepackage{enumitem}
\usepackage{bm}
\usepackage{pict2e}
\usepackage{subcaption}
\usepackage[mathscr]{euscript}

\usepackage{indentfirst}
\newlist{multienum}{enumerate}{1}
\setlist[multienum]{
    label=\alph*),
    before=\begin{multicols}{2},
    after=\end{multicols}
}
\usepackage{hyperref}
\usepackage[nameinlink]{cleveref}
\hypersetup{
    colorlinks=true,
    linkcolor=black,
    filecolor=black,      
    urlcolor=black,
    citecolor={black}
    }

\usepackage{enumitem}

\numberwithin{equation}{section}
\newcounter{theorem}
\newtheorem{thm}[theorem]{Theorem}

\newtheorem{cor}[theorem]{Corollary}
\newtheorem{prop}[theorem]{Proposition}
\newtheorem{lem}[theorem]{Lemma}
\newtheorem{theorem}[theorem]{Theorem}

\numberwithin{theorem}{section}

\theoremstyle{definition}
\newtheorem{defn}[theorem]{Definition}

\newtheorem{claim}[theorem]{Claim}

\newtheoremstyle{newtheoremstyledefn}
{3pt}
{3pt}
{}
{\parindent}
{\bfseries}
{.}
{0.5em}
{} 

\theoremstyle{newtheoremstyledefn}
\newtheorem{remark}[theorem]{Remark}

\newcommand{\A}{\mathcal{A}}
\newcommand{\B}{\mathcal{B}}

\newcommand{\C}{\mathbb{C}}

\newcommand{\CD}{\mathcal{D}}

\newcommand{\D}{\mathbb{D}}

\newcommand{\G}{\mathcal{G}}
\renewcommand\H{\mathcal{H}}

\newcommand{\CK}{\mathcal{K}}

\newcommand{\R}{\mathbb{R}}
\newcommand{\CR}{\mathcal{R}}

\newcommand{\Z}{\mathbb{Z}}
\newcommand{\CZ}{\mathcal{Z}}

\newcommand{\reg}{\textnormal{reg}}

\newcommand{\al}{\alpha}
\renewcommand{\dim}{\textnormal{dim}}

\newcommand{\smallnorm}[1]{\lVert#1\rVert}
\newcommand{\sing}{\textnormal{sing}}

\newcommand{\dist}{\textnormal{dist}}

\newcommand{\spt}{\textnormal{spt}}

\newcommand{\graph}{\textnormal{graph}}

\newcommand{\Lip}{\textnormal{Lip}}

\newcommand{\mb}[1]{\llbracket{#1}\rrbracket}

\newcommand{\vphi}{\varphi}
\renewcommand\phi{\vphi}
\newcommand{\eps}{\varepsilon}
\renewcommand{\bf}[1]{\mathbf{#1}}
\renewcommand{\d}{\mathrm{d}}

\renewcommand{\Re}{\textnormal{Re}}
\renewcommand{\Im}{\textnormal{Im}}
\makeatletter

\address{\textnormal{Department of Pure Mathematics and Mathematical Statistics, University of Cambridge}}
\email{sas258@cam.ac.uk}

\usepackage{fancyhdr}
\title{Regularity of $Q$-valued $C^1$ Harmonic Functions with Applications to Stable Minimal Hypersurfaces}
\author{Sidney Stanbury}
\date{}
\begin{document}
\begin{abstract}
We show that the derivative of a $Q$-valued $C^1$ function on an $n$-dimensional domain that is locally harmonic on its regular set is locally Dirichlet minimising, and that its branch set coincides with that of its derivative. As a consequence, we get that a $Q$-valued $C^1$ harmonic function is $C^{1,\al}$ and that its branch set is countably $(n-2)$-rectifiable when $n\geq 3$ and discrete when $n=2$. 

Our methods are directly applicable to the regularity theory of stable minimal hypersurfaces. We show that a multi-valued Lipschitz function whose graph is a two-sided, properly immersed, stable minimal $n$-dimensional hypersurface with locally $\H^{n-2}$-finite singular set, and whose Lipschitz constant is sufficiently small is in fact multi-valued $C^{1,\al}$. As a consequence, we obtain optimal regularity conclusions for a branched sheeting theorem of Bellettini for singular immersions. 

When $n=2$, we obtain stronger conclusions and show that locally about a branch point, a $Q$-valued $C^1$ harmonic function is the real part of a Dirichlet minimiser which is itself obtained from a finite number of holomorphic functions via classical branching, thus giving a full series expansion for the function near the branch point.  Furthermore, we prove that the class of branched, two-sided, properly immersed, stable minimal surfaces with isolated singularities is compact under varifold convergence.
\end{abstract}
\maketitle
\tableofcontents
\section{Introduction}
\subsection{$Q$-valued $C^1$ Harmonic Functions}\label{sec: C^1 harmonic} Let $n,Q\geq 2$ and $k\geq 1$. Suppose $\Omega\subseteq\R^n$ is a domain and $w\in C^1(\Omega;\A_Q(\R^k))$. If $U\subseteq\Omega$ is open, we say that $w$ is \emph{locally harmonic} in $U$ if every point in $U$ is contained in an open ball in which $w$ decomposes as the sum of $Q$ single-valued harmonic functions. We say that $w$ is a $Q$-valued $C^1$ \emph{harmonic function} if $w$ is locally harmonic in its regular set $\CR_w$; namely, the set of points in $\Omega$ contained in an open ball in which $w$ decomposes as the sum of $Q$ single-valued $C^1$ functions.

The condition for a $Q$-valued $C^1$ function $w:B^n_1(0)\to\A_Q(\R^k)$ to be a $Q$-valued $C^1$ harmonic function is purely a condition on the regular set of the function; it only requires $w$ to satisfy an equation on its regular set $\CR_w$, and thus \emph{does not} require $w$ to satisfy any variational property across its branch set $\B_w:=B^n_1(0)\setminus\CR_w$. We stress that this is in contrast to many other notions of multi-valued harmonic functions in the literature, for which this is not the case (cf.~Remark \ref{remark: Literaturature specification}). We note that the regular set $\CR_w$ is an open and dense subset of $B^n_1(0)$, and thus this condition is never vacuously true (cf.~Remark \ref{remark: separation principle}). However, an open and dense subset of $B^n_1(0)$ can have arbitrarily small measure, so it is not a priori clear whether this condition is sufficiently strong to give information about the branch set. For instance, $w$ is Dirichlet stationary in $\CR_w$, but it is not clear whether the corresponding variational identities can be extended across $\B_w$. 

In the two-valued $C^{1,\al}$ harmonic setting, the combined works of \cite{SW16,KW13,Kru14ConstFreq} establish optimal dimension bounds on the branch set, optimal regularity conclusions, and a fine analysis of the branch set, such as an asymptotic expansion for the function about $\H^{n-2}$-a.e. branch point, and countable $(n-2)$-rectifiability of the branch set. As far as the author is aware, prior to this work, no such results were known for $Q$-valued $C^1$ harmonic functions, even in the case when $Q=2$.

The following is our main theorem for $Q$-valued $C^1$ harmonic functions:
\begin{thm}\label{thm: derivative is Dir min}
If $w:B^n_1(0)\to\A_Q(\R^k)$ is a $Q$-valued $C^1$ harmonic function, then $Dw$ is locally Dirichlet minimising and $\B_w=\B_{Dw}$.
\end{thm}
As a first corollary to Theorem \ref{thm: derivative is Dir min}, we obtain the following regularity improvement result for $Q$-valued $C^1$ harmonic functions:
\begin{thm}\label{thm: main-reg}
If $w:B^n_1(0)\to\A_Q(\R^k)$ is a $Q$-valued $C^1$ harmonic function, then $w$ is $C^{1,\al}$ for some $\al=\al(n,Q)\in(0,1]$ with $\al(2,Q)=1/Q$ and $\al(n,2)=1/2$. 
\end{thm}

As a second corollary to Theorem \ref{thm: derivative is Dir min}, we obtain the following description of the branch set of a $Q$-valued $C^1$ harmonic function:
\begin{thm}\label{thm: rectifiability of branch set}
Suppose $w:B^n_1(0)\to\A_Q(\R^k)$ is a $Q$-valued $C^1$ harmonic function. Then $\B_w$ is countably $(n-2)$-rectifiable when $n\geq 3$, and discrete when $n=2$.
\end{thm}

As a third corollary to Theorem \ref{thm: derivative is Dir min}, we obtain the following local finiteness of $\H^{n-2}$-measure and Minkowski-type estimates for the branch set in the two-valued setting:
\begin{thm}\label{thm: Minkowski-type estimates in two-valued setting}
If $w:B^n_1(0)\to\A_2(\R^k)$ is a two-valued $C^1$ harmonic function, then we have $\H^{n-2}(\B_w\cap K)<\infty$ for every compact set $K\subset B^n_1(0)$. Furthermore, we have the stronger Minkowski-type estimate 
$$
\H^n(N_r(\B_w)\cap K)\leq C(K,w)r^2,\quad\text{ for all }r\in(0,1),
$$
where $N_r(\B_w)$ denotes the $r$-neighbourhood of $\B_w$ in $B^n_1(0)$.
\end{thm}

\begin{remark}\label{remark: main theorem remarks}
We digress to make a few remarks.
\begin{itemize}[leftmargin=2em]
\item[(a)] A key ingredient in the proof of Theorem \ref{thm: derivative is Dir min} is an analogous result we establish for $W^{2,2}_{\textnormal{loc}}$ functions that are stationary for outer variations (cf.~Proposition \ref{prop: W^2,2 Dir stationary is C^1,alpha}). This latter result gives an alternative proof of the homogeneity gaps in \cite[Section 3.4]{LucheseSchauder} which bypasses the need for aperture estimates and matrix-valued measures. It then follows by standard arguments (cf.~\cite[Section 4]{LucheseSchauder}) that these homogeneity gaps lead to a $Q$-valued Schauder theory for $C^{1,\al}$ weak solutions of a class of quasilinear elliptic systems. In particular, it follows from planar frequency arguments (cf.~\cite[Section 3]{KMW}) that if $u\in C^{1,\al}(B^2_1(0);\A_Q(\R^k))$ is such that the integral varifold $\bf{v}(u)$ associated to the graph of $u$ is stationary in $\R^k\times B^2_1(0)$, then $u$ is $C^{1,1/Q}$. Similarly, we recover the optimal result of \cite{SW16} which establishes $C^{1,1/2}$ regularity for two-valued $C^{1,\al}$ stationary graphs.

\item[(b)] Combining Theorem \ref{thm: derivative is Dir min} with \cite[Theorems A \& B]{KW17} gives an asymptotic for $Dw$ near $\H^{n-2}$-a.e. branch point, and a decomposition for $\B_w$ analogous to those obtained for locally Dirichlet minimising functions therein.

\item[(c)] If $w$ is a $Q$-valued $C^1$ function, then the set $\B_{Dw}$ is generally not contained in $\B_w$. Furthermore, for $Q\geq 3$, there exists a $Q$-valued $C^1$ function $w_Q$ such that $\B_{w_Q}$ is not contained in $\B_{Dw_Q}$. In contrast, the branch set of a two-valued $C^1$ function is always contained in that of its derivative (cf.~Appendix \ref{app: containments of branch sets}).

\item[(d)] \cite{SW16} obtained the sharp regularity conclusion that a two-valued $C^{1,\al}$ harmonic function is $C^{1,1/2}$. Theorem \ref{thm: main-reg} extends this optimal regularity improvement under the strictly weaker $C^1$ assumption. In contrast to \cite{SW16}, Theorem \ref{thm: main-reg} is obtained \emph{without} appealing to a multi-valued Schauder theory. Furthermore, the regularity improvement given by Theorem \ref{thm: main-reg} is optimal when $n=2$. As far as the author is aware, it remains open what the optimal value for $\al(n,Q)$ is when both $n\geq 3$ and $Q\geq 3$.

\item[(e)] Similarly, Theorem \ref{thm: main-reg} implies that the results of \cite{Kru14ConstFreq} for the branch set of a two-valued $C^{1,\al}$ harmonic function hold under the strictly weaker $C^1$ assumption. In particular, the branch set of a two-valued $C^1$ harmonic function contains a relatively open and dense subset on which the branch set is locally a codimension-$2$ real analytic submanifold. 

\item[(f)] Theorem \ref{thm: rectifiability of branch set} extends \cite[Theorem A]{KW13} to the arbitrary $Q$-valued and strictly weaker $C^1$ setting. In fact, Theorem \ref{thm: Minkowski-type estimates in two-valued setting} improves the countable $(n-2)$-rectifiability conclusion therein by establishing local finiteness of $\H^{n-2}$-measure and Minkowski-type estimates for the branch set. Furthermore, Theorem \ref{thm: rectifiability of branch set} and \cite[Appendix]{SW16} imply that either $\B_w$ is empty, or $\dim_\H\B_w= n-2$ and $\H^{n-2}(\B_w)>0$. Moreover, Theorem \ref{thm: rectifiability of branch set} implies that the average of a $Q$-valued $C^1$ harmonic function is a single-valued (analytic) harmonic function.

\item[(g)] Theorem \ref{thm: rectifiability of branch set} and \cite[Lemma 3.2]{MinterCampanato} establish the `Classification Hypothesis' of \cite[p.~30]{MinterCampanato}. As a consequence, the statements of \cite[Theorems 3.1 \& 3.2]{MinterCampanato} hold unconditionally.

\item[(h)] For $Q\geq 3$, Theorem \ref{thm: derivative is Dir min} and \cite[Theorem 1.2]{DLMSVDirRectifiability} give local finiteness of $\H^{n-2}$-measure and Minkowski-type estimates for the set $\{Dw=Q\mb{Dw_a}\}$, provided that $Dw\not\equiv Q\mb{Dw_a}$. Furthermore, by Theorem \ref{thm: derivative is Dir min}, if $Dw\equiv Q\mb{Dw_a}$, then $\B_w$ is empty.
\end{itemize}
\end{remark}

The proof of Theorem \ref{thm: derivative is Dir min} follows from an inductive argument on $Q$. Interestingly, a necessary ingredient in the proof of the induction step is that the branch set of a $(<Q)$-valued $C^1$ harmonic function has vanishing $2$-capacity, rather than being merely of Hausdorff dimension at most $n-2$. In particular, it is unclear whether it would be possible to obtain the dimension bound given by Theorem \ref{thm: rectifiability of branch set} without establishing (or appealing to) a vanishing $2$-capacity conclusion.

It follows from \cite[Lemma 6.9]{DLMSVDirRectifiability} that any two $Q$-valued locally Dirichlet minimising functions on $B^n_1(0)$ which agree in a non-empty open set must agree on all of $B^n_1(0)$. As a consequence of our arguments, we obtain the following weak unique continuation for $Q$-valued $C^1$ harmonic functions:
\begin{thm}\label{thm: weak unique continuation}
Suppose $w,v:B^n_1(0)\to\A_Q(\R^k)$ are $Q$-valued $C^1$ harmonic functions. If $U\subseteq B^n_1(0)$ is non-empty and open such that $w\equiv v$ in $U$, then $w\equiv v$ in $B^n_1(0)$.
\end{thm}
As a corollary to Theorem \ref{thm: weak unique continuation}, we obtain the following sharp quantitative unique continuation for $Q$-valued $C^1$ harmonic functions:\footnote{See Section \ref{sec: notation and terminology} for the definition of the $1$-jet of a $Q$-valued $C^1$ function.}
\begin{thm}\label{thm: quantitative unique continuation}
Suppose $w,v:B^n_1(0)\to\A_Q(\R^k)$ are $Q$-valued $C^1$ harmonic functions.
\begin{enumerate}
\item If $\dim_\H\{w=v\}>n-1$, then $w\equiv v$ in $B^n_1(0)$.
\item If $\dim_\H\{j_1w=j_1v\}>n-2$, then $w\equiv v$ in $B^n_1(0)$.
\end{enumerate}
\end{thm}

In order to state our next result, we will need some further notation. Denote $\D:=\{z\in \C:|z|< 1\}$. For $\zeta=\sum_{i=1}^Q\mb{\zeta^i}\in\A_Q(\C^k)$ we denote the \emph{real part} of $\zeta$ by
$$
\Re(\zeta):=\sum_{i=1}^Q\mb{\Re(\zeta^i)}\in\A_Q(\R^k).
$$
Suppose $f:\D\to\A_Q(\C^k)$ is a complex $Q$-valued function. We define the \emph{real part} of $f$ to be the function $\Re (f):\D\to\A_Q(\R^k)$ defined by $\Re (f)(z):=\Re(f(z))$. 

It is well-known that any (real) single-valued harmonic function on a simply connected planar domain is the real part of a holomorphic function. Theorem \ref{thm: rectifiability of branch set} gives the following (local) generalisation of this fact to the $Q$-valued setting:

\begin{thm}\label{thm: real part of Dir min} 
Suppose $w:\D\to\A_Q(\R^k)$ is a $Q$-valued $C^1$ harmonic function. Then, for every $x\in \B_w$ there exists $r\in(0,1-|x|),N\in\{1,\dots,Q\}, Q_1,\dots,Q_N\in\{1,\dots,Q\}$ and holomorphic functions $F_1,\dots,F_N:\D\to\C^k$ such that the function $g:\D\to\A_Q(\C^k)$ defined by
\begin{equation*}
g(z):=\sum_{j=1}^N\sum_{\zeta^{Q_j}=z}\mb{F_j(\zeta)} 
\end{equation*}
is Dirichlet minimising, and $w(x+r\cdot)\equiv \Re (g)$ in $\D$.
\end{thm}
\begin{remark}
Let $\mathscr{V}_Q$ denote the collection of complex varieties of the form
$$
\{(z,\eta)\in \D\times\C^k:\exists\,\zeta\in \D,\zeta^q=z,\eta=f(\zeta)\},
$$
where $q\in\{1,\dots,Q\}$ and $f:\D\to\C^k$ is holomorphic. Then by Theorem \ref{thm: real part of Dir min}, any two-dimensional $Q$-valued $C^1$ harmonic function is locally described about any of its branch points by a finite collection of varieties in $\mathscr{V}_Q$ via classical branching.
\end{remark}
As a corollary to Theorem \ref{thm: real part of Dir min}, we obtain the following series expansion for a two-dimensional $Q$-valued $C^1$ harmonic function near a branch point:
\begin{cor}\label{cor: Full Series Expansion in Disc}
Suppose $w:\D\to\A_Q(\R^k)$ is a $Q$-valued $C^1$ harmonic function. Then, for every $x\in \B_w$ there exists $r\in(0,1-|x|),N\in\{1,\dots,Q\}, Q_1,\dots,Q_N\in\{1,\dots,Q\}$, and, for each $j\in\{1,\dots,N\}$ there exists a single-valued harmonic function $h_j:\D\to\R^k$ and a sequence $(a_{j,m})\subset\C^k$, such that
$$
w(x+rz)=\sum_{j=1}^{N}\sum_{\zeta^{Q_j}=z}\biggl\llbracket h_j(z)+\Re\biggl(\sum_{m\ge Q_j+1:Q_j\nmid m}a_{j,m}\zeta^m\biggr)\biggr\rrbracket,\quad\text{ for all }z\in\D.
$$
In particular, if $w\in C^{1,\beta}$ for some $\beta>1-1/Q$, then $w\in C^{2,1/Q}$.
\end{cor}

It follows from our proof of Theorem \ref{thm: real part of Dir min} (cf.~Claim \ref{claim: real part of Dir min}) and the regularity theory for $Q$-valued Dirichlet minimisers (cf.~\cite[Theorem 2.13]{Almgren00}), that if $w\in L^\infty_{\textnormal{loc}}(\D;\A_Q(\R^k))$ is locally harmonic away from a discrete set, then $w$ is $C^{0,1/Q}$. Furthermore, by \cite[Proposition 3.1]{Spadaro2010ComplexVarieties} (cf.~\cite[Theorem 6.1]{DLS14GradientL^pEstimates}), $w\in W^{1,p}_{\textnormal{loc}}(\D;\A_Q(\R^k))$ for all $p\in[1,2Q/(Q-1))$. In particular, if $w$ is additionally assumed to be in $W^{1,\infty}_{\textnormal{loc}}(\D;\A_Q(\R^k))$, then applying the above with $Dw$ gives that $w$ is a $Q$-valued $C^{1,1/Q}$ harmonic function and $w\in W^{2,p}_{\textnormal{loc}}(\D;\A_Q(\R^k))$, for all $p\in[1,2Q/(Q-1))$. As a special case, we obtain the following sharp integrability estimate for 2-dimensional $Q$-valued $C^1$ harmonic functions:
\begin{thm}
If $w:\D\to\A_Q(\R^k)$ is a $Q$-valued $C^1$ harmonic function, then we have that $w\in W_{\textnormal{loc}}^{2,p}(\D;\A_Q(\R^k))$ for all $p\in[1,2Q/(Q-1))$.
\end{thm}
We have the following partial generalisation of Corollary \ref{cor: Full Series Expansion in Disc} to arbitrary dimensions:
\begin{thm}\label{thm: Uniqueness of Tangent Maps}
Suppose $w:B^n_1(0)\to\A_Q(\R^k)$ is a $Q$-valued $C^1$ harmonic function. For $\H^{n-2}$-a.e. $z\in\B_w$, there exists $\rho_z\in(0,1-|z|)$, an integer $M_z\geq 1$, positive integers $p^{(1)}_z,\dots,p^{(M_z)}_z$ with $\sum_{\ell=1}^{M_z}p_z^{(\ell)}=Q$, positive integers $N_z^{(1)},\dots, N_z^{(M_z)}$, relatively prime pairs of positive integers $(k_z^{(1)}, q_z^{(1)}),\dots,(k_z^{(M_z)}, q_z^{(M_z)})$ with $k^{(\ell)}_z\geq q^{(\ell)}_z+1$ and $N^{(\ell)}_zq^{(\ell)}_z\leq p^{(\ell)}_z$ for each $\ell\in\{1,\dots,M_z\}$, and orthogonal rotations $Q^{(1)}_z,\dots,Q^{(M_z)}_z$ of $\R^n$ such that for each $x\in B^n_{\rho_z}(0)$, we have
$$
w(z+x)=\sum_{\ell=1}^{M_z}\sum_{j=1}^{p^{(\ell)}_z}\mb{h^{(\ell)}_z(z+x)+\phi_j^{(z,\ell)}(Q^{(\ell)}_zx)+\eps_j^{(z,\ell)}(x)},
$$
where for each $\ell\in\{1,\dots,M_z\}$:
\begin{itemize}
\item[(i)] $h^{(\ell)}_z:B^n_{\rho_z}(z)\to\R^k$ is a single-valued harmonic function;
\item[(ii)] either $p^{(\ell)}_z=1$ and $\phi_1^{(z,\ell)}\equiv0$, or $p^{(\ell)}_z\geq 2$ and $\phi^{(z,\ell)}:\R^n\to\A_{p^{(\ell)}_z}(\R^k)$ is a non-zero, homogeneous, average-free $p^{(\ell)}_z$-valued $C^1$ harmonic function given by
$$
\phi^{(z,\ell)}(x)=\sum_{j=1}^{p^{(\ell)}_z}\mb{\phi_j^{(z,\ell)}}=(p^{(\ell)}_z-N^{(\ell)}_zq^{(\ell)}_z)\mb{0}+\sum_{m=1}^{N^{(\ell)}_z}\textnormal{Re}(c_m^{(z,\ell)}(x_1+ix_2)^{k^{(\ell)}_z/q^{(\ell)}_z})
$$
for some $c_m^{(z,\ell)}\in\C^k\setminus\{0\}$;
\item[(iii)] for each $x\in B^n_{\rho_z}(0)$, $\eps^{(z,\ell)}_j(x)\in\R^k$ such that $\sum_{j=1}^{p^{(\ell)}_z}\eps_j^{(z,\ell)}(x)=0$ and
$$
\sigma^{-n}\int_{B^n_\sigma(0)}\sum_{j=1}^{p^{(\ell)}_z}|\eps_j^{(z,\ell)}|^2= o(\sigma^{2k^{(\ell)}_z/q^{(\ell)}_z}),\quad\text{ for all }\sigma\in(0,\rho_z).
$$
\end{itemize}
\end{thm}
\begin{remark}
The proof of Theorem \ref{thm: Uniqueness of Tangent Maps} shows that $w$ has a unique tangent map at every point of $\CK_w$ at which the (Dirichlet minimising) function $Dw$ has a unique tangent map (cf.~Section \ref{sec: Proof of Uniqueness of Tangent Maps}). In particular, by \cite[Theorem 5.3]{DLS11Qvalued} (c.f.~\cite[Theorem A]{KW17}), when $n=2$, $w$ has a unique tangent map at every point of $\CK_w$. Of course, this uniqueness also follows from Theorem \ref{thm: real part of Dir min}.
\end{remark}
\subsection{Applications to Stable Minimal Hypersurfaces}
The philosophy emerging from Section \ref{sec: C^1 harmonic} is that the derivative of a multi-valued function satisfying a given variational property should minimise a corresponding functional. In particular, the methods of the present paper have direct applications to the regularity theory of stable minimal hypersurfaces. The following is the key result connecting our methods to the latter theory:
\begin{thm}\label{thm: Lip+W^2,2 Stationary + Lip small implies C^1,alpha}
Fix $n\geq2$ and $k,Q\geq 1$. There exists $\eps=\eps(n,Q)\in(0,1)$ such that the following holds. Suppose $u\in W^{2,2}_{\textnormal{loc}}(B^n_1(0);\A_Q(\R^k))$ is Lipschitz and that the integral varifold $\bf{v}(u)$ is stationary in $\R^k\times B^n_1(0)$. If $\Lip(u)\leq\eps$, then $u\in C^{1,\al}(B^n_{1/2}(0);\A_Q(\R^k))$ for some $\al=\al(n,Q)\in(0,1]$ with $\al(2,Q)=1/Q$ and $\al(n,2)=1/2$.
\end{thm}
As a corollary to Theorem \ref{thm: Lip+W^2,2 Stationary + Lip small implies C^1,alpha}, we obtain the following optimal regularity improvement for the conclusions of \cite[Theorem 5]{Bel25}\footnote{Below, $\bf{v}(u)$ denotes the integral varifold associated to the graph of $u$.}:
\begin{thm}\label{thm: Bel25 Theorem 5 reg improvement}
Fix $n\geq 2$ and $Q\geq 1$. There exists $\eps=\eps(n,Q)\in(0,1)$ such that the following holds. Suppose that $M$ is a properly immersed, two-sided, stable minimal hypersurface in $C_{2R}\setminus\Sigma$, where $\Sigma$ is closed in $C_{2R}:=(-2R,2R)\times B^n_{2R}(0)$ and $\textnormal{cap}_2(\Sigma\cap\overline{M})=0$. Suppose further that $M\cap C_{R/2}$ is non-empty, $(\omega_n(2R)^n)^{-1}\H^n(M\cap C_{2R})\leq Q+1/2, \sup_{M\cap C_{2R}}|x_{n+1}|<R/4$ and
$$
E_M(R)^2:=R^{-n}\int_{M\cap C_R}(1-(\nu_M\cdot e_{n+1})^2)\leq\eps.
$$
Then, there exists $q\in\{1,\dots,Q\}$ and a $q$-valued $C^{1,\al}$ function $u:B^n_{R/4}(0)\to\A_q(\R)$ such that
$$
\overline{M}\cap C_{R/4}=\bf{v}(u),
$$
for some $\al=\al(n,Q)\in(0,1]$ with $\al(2,Q)=1/Q$ and $\al(n,2)=1/2$.
\end{thm}

In fact, by the proof of Theorem \ref{thm: Bel25 Theorem 5 reg improvement}, we can improve the regularity of the Lipschitz multi-valued functions obtained in \cite[Theorem 1.9]{WZ26StableMinHypNonBranchedSingSet} to multi-valued $C^{1,\al}$. This regularity improvement combined with \cite{LucheseSchauder,KMW} and \cite[Theorem 1.3]{WZ26StableMinHypNonBranchedSingSet} gives the following dimension bound for the varifold limit of a sequence of stable minimal hypersurfaces which are properly-immersed away from a set of vanishing $2$-capacity\footnote{Below, $\sing_\star(V):=\spt\smallnorm{V}\setminus\reg_\star(V)$, where $\reg_\star(V)$ denotes the smoothly \textit{immersed} regular set of an $n$-dimensional stationary integral varifold $V$ in $B^{n+1}_2(0)$ (cf.~\cite[Definition 1.2]{WZ26StableMinHypNonBranchedSingSet}).}:
\begin{cor}\label{cor: SMH dim bound of limit}
Fix $\Lambda>0$. Suppose that $(V_j)_j$ is a sequence of $n$-dimensional stationary integral varifolds in $B^{n+1}_2(0)$ which satisfy $0\in\spt\smallnorm{V_j}$ and $\smallnorm{V_j}(B^{n+1}_2(0))\leq\Lambda$. Suppose that for every $j\geq 1$, $M_j:=\reg_\star(V_j)$ is a properly immersed, two-sided, stable minimal hypersurface in $B^{n+1}_2(0)\setminus\Sigma_j$, where $\Sigma_j:=\overline{M_j}\setminus M_j$ is such that $\textnormal{cap}_2(\Sigma_j;\spt\smallnorm{V_j})=0$. Then, after passing to a subsequence, there exists a stationary integral varifold $V$ in $B^{n+1}_2(0)$ such that $V_j\rightharpoonup V$ and $\dim_\H\sing_\star(V)\leq n-2$.
\end{cor}
When $n=2$, our methods lead to the following compactness theorem for the class of branched, two-sided, properly immersed, stable minimal surfaces with isolated singularities (cf.~\cite[p.~6]{Bel25}):
\begin{thm}\label{thm: Compactness Branched SMH n=2}
Fix $\Lambda>0$. Suppose that $(M_j)_j$ is a sequence of properly immersed, two-sided, stable minimal surfaces in $B^3_2(0)\setminus\Sigma_j$, where $\Sigma_j\subset B^3_2(0)$ is closed with $\textnormal{cap}_2(\Sigma_j;\overline{M_j})=0$. Suppose further that $0\in\overline{M_j}$ and $\H^2(M_j\cap B^3_2(0))\leq \Lambda$ for all $j\geq 1$. Then, after passing to a subsequence, there exists a discrete set $\Sigma\subset B^3_2(0)$ and a properly immersed, two-sided, stable minimal hypersurface $M$ in $B^3_2(0)\setminus\Sigma$ such that $ M_j$ converges to $M$ in both the varifold topology of $B^3_2(0)$ and the multi-valued $C^1$ topology of $B^3_2(0)\setminus\Sigma$. Furthermore, we have
$$
\H^0(\sing_\star(M)\cap\Omega)\leq\liminf_{j\to\infty}\H^0(\sing_\star(M_j)\cap\Omega),\quad\text{ for every open }\Omega\Subset B^3_2(0).
$$
\end{thm}
\begin{remark} 
We digress to make a few remarks.
\begin{itemize}[leftmargin=2em]
\item[(a)] It remains open whether the non-immersed singular set of $V$ in Corollary \ref{cor: SMH dim bound of limit} has vanishing $2$-capacity when $n\geq 3$. We plan to address this question in future work by building on the methods of the present paper and adapting the geometric measure theory techniques developed in \cite{Simon93Cylindrical,KW13,KW17,KW21}. 
\item[(b)] We note that a discrete subset of $B^3_2(0)$ has zero $2$-capacity. In particular, Theorem \ref{thm: Compactness Branched SMH n=2} implies that the class of branched, two-sided, properly immersed, stable minimal surfaces in $B^3_2(0)$ with locally finite singular set is compact in the varifold topology. 
\item[(c)] Theorem \ref{thm: Compactness Branched SMH n=2} implies that if $M$ is a properly immersed, two-sided, stable minimal surface in $B^3_2(0)\setminus\Sigma$, where $\Sigma\subset B^3_2(0)$ is closed with $\textnormal{cap}_2(\Sigma)=0$, then, provided $\Sigma$ does not contain any removable singularities of $\overline{M}$, $\Sigma$ must be discrete. This conclusion is optimal and isolated singularities are unavoidable in general (cf.~\cite{Henneberg1875}, \cite[p.77-78]{Lawson77lectures}, \cite[Section 4]{MoyaPerez23MinSurfaces} \& \cite[Corollary 1]{SW07}). In fact, one can use the Weierstrass–Enneper representation to construct simple branched examples (cf.~\cite[Section 8]{Osserman86Survey} with $f(\zeta)=2\zeta$ and $g(\zeta)=\zeta$): the image of the map
$$
\D\ni\zeta\mapsto\Bigl(\textnormal{Re}\Bigl(\zeta^2-\frac{\zeta^4}{2}\Bigr),\textnormal{Im}\Bigl(\zeta^2+\frac{\zeta^4}{2}\Bigr),\frac{4}{3}\textnormal{Re}(\zeta^3)\Bigr)
$$
defines a branched, two-sided, properly immersed, stable minimal surface in $B^3_{1/8}(0)\setminus\{0\}$.
\end{itemize}
\end{remark}

Furthermore, combining Theorem \ref{thm: Lip+W^2,2 Stationary + Lip small implies C^1,alpha} with \cite[Corollary 5.2]{LucheseSchauder}, we obtain the following small-slope Bernstein theorem for stationary $W_{\textnormal{loc}}^{2,2}$ Lipschitz graphs:
\begin{cor}\label{cor: W^2,2+Lip Small-Slope Bernstein}
Fix $n\geq 2$ and $k,Q\geq 1$. There exists $\eps=\eps(n,k,Q)\in(0,1)$ such that the following holds. Suppose $u\in W^{2,2}_{\textnormal{loc}}(\R^n;\A_Q(\R^k))$ is Lipschitz and that the integral varifold $\bf{v}(u)$ is stationary in $\R^{n+k}$. If $\Lip(u)\leq\eps$, then $u$ is affine.
\end{cor}

Finally, we show that the smallness assumption on the Lipschitz constant in Theorem \ref{thm: Lip+W^2,2 Stationary + Lip small implies C^1,alpha} can be removed for the price of $C^1$ regularity:
\begin{thm}\label{thm: C^1+W^2,2 implies C^1,alpha}
Fix $n\geq2$ and $k,Q\geq 1$. Suppose $u\in C^1(B^n_1(0);\A_Q(\R^k))\cap W^{2,2}_{\textnormal{loc}}(B^n_1(0);\A_Q(\R^k))$ and that the integral varifold $\bf{v}(u)$ is stationary in $\R^k\times B^n_1(0)$. Then $u\in C^{1,\al}(B^n_{1/2}(0);\A_Q(\R^k))$ for some $\al=\al(n,Q)\in(0,1]$ with $\al(2,Q)=1/Q$ and $\al(n,2)=1/2$.
\end{thm}
\begin{remark}
We make a few further remarks.
\begin{itemize}[leftmargin=2em]
\item[(a)] By \cite{SW07}, the conclusions of Theorems \ref{thm: Lip+W^2,2 Stationary + Lip small implies C^1,alpha}, \ref{thm: Bel25 Theorem 5 reg improvement} and \ref{thm: C^1+W^2,2 implies C^1,alpha} are optimal when $Q=2$.
\item[(b)] We can assume without loss of generality that $\Sigma$ in Theorem \ref{thm: Bel25 Theorem 5 reg improvement} does not contain any removable singularities of $\overline{M}$. Hence, \cite{SW16,KW21} implies that when $Q=2$, $\Sigma$ is discrete when $n=2$, and countably $(n-2)$-rectifiable when $n\geq 3$.
\item[(c)] The $W^{2,2}_{\textnormal{loc}}$ assumption in Theorem \ref{thm: Lip+W^2,2 Stationary + Lip small implies C^1,alpha} and Corollary \ref{cor: W^2,2+Lip Small-Slope Bernstein} is necessary in general; there exists a sequence of non-affine Lipschitz stationary graphs $u_j:\R^2\to\A_2(\R^2)$ such that $\Lip(u_j)\downarrow0$ and each $u_j$ is not differentiable at the origin (cf.~\cite[Example 1.2]{BK17}).
\item[(d)] The condition in Theorem \ref{thm: C^1+W^2,2 implies C^1,alpha} that $\bf{v}(u)$ is stationary can be replaced with the weaker condition that every point of $\CR_u$ is contained in an open ball in which $u$ decomposes as the sum of $Q$ single-valued $C^1$ functions solving the minimal surface system.
\item[(e)] By \cite{KMW,LucheseSchauder}, if $u:B^n_1(0)\to\A_Q(\R^k)$ is in Theorem \ref{thm: Lip+W^2,2 Stationary + Lip small implies C^1,alpha} or Theorem \ref{thm: C^1+W^2,2 implies C^1,alpha}, then $\dim_\H\B_u\leq n-2$ (cf.~Remark \ref{remark: main theorem remarks}\,(a)). Furthermore, when $Q=2$, it follows from \cite{SW16,KW21} that $\B_u$ is discrete when $n=2$, and countably $(n-2)$-rectifiable when $n\geq 3$.
\end{itemize}
\end{remark}

\textbf{AI Statement:} These results were obtained by combining ideas of the author with a publicly available LLM over the course of several months. The proof of Theorem \ref{thm: graph of derivative is calibrated in 2d} was strongly supported by AI, and its initial form had very limited human input. Another significant input given by AI was the observation that the functional $\mathcal{A}_U$ in the proof of Proposition \ref{prop: W^2,2 Dir stationary is C^1,alpha} is a null Lagrangian. Other uses of AI, although beneficial, were much more localised, and were of a technical and collaborative nature to assist in the process of developing the author's ideas. Among these uses, parts of some of the proofs in Appendices \ref{app: containments of branch sets}, \ref{app: Multi-Valued Functions} and \ref{app: Obstruction to Calibrations in Higher Dimensions} were the only arguments that benefited from significant AI assistance. AI assistance was also used to check for typos and errors. The paper does not contain AI written text and the author takes full responsibility for its correctness.

\textbf{Acknowledgements:} This research was supported by the Engineering and Physical Sciences Research Council (EPSRC) - EP/W524633/1. The author is grateful to Paul Minter and Neshan Wickramasekera for many helpful discussions. After completing the present work, the author was made aware that Theorems \ref{thm: Lip+W^2,2 Stationary + Lip small implies C^1,alpha}, \ref{thm: Bel25 Theorem 5 reg improvement} and Corollary \ref{cor: SMH dim bound of limit} have been independently obtained by Mattia Luchese using alternative methods (cf.~\cite[Corollaries 6, 9 \& Theorem 12]{Luch26SmallLap}).

\subsection{Preliminaries}\label{sec: notation and terminology}
In this section, we discuss notation and terminology which will be used throughout the paper. We refer the reader to \cite{Almgren00,DLS11Qvalued} for background on multi-valued functions.

Fix $n\geq 2$ and $k,Q\geq 1$. We begin with the notion of a differentiable $Q$-valued function:
\begin{defn}
Suppose $w:B^n_1(0)\to\A_Q(\R^k)$ and $x\in B^n_1(0)$. We say that $w$ is \emph{differentiable at} $x$ if there exist $Q$ linear maps $L_x^1,\dots,L_x^Q:\R^n\to\R^{k}$ such that
\begin{equation}\label{eq: defn of differentiability}
\lim_{y\to x}|x-y|^{-1}\G\biggl(w(y),\sum_{i=1}^Q\mb{w^i(x)+L_x^i\cdot(y-x)}\biggr)=0,
\end{equation}
We say that $w$ is \emph{differentiable} if $w$ is differentiable at every point of $B^n_1(0)$.
\end{defn}
\textbf{Remark:} The above defined notion of differentiability is strictly weaker than the notion defined in \cite[Definition 1.9]{DLS11Qvalued}, which is often referred to as `strongly affine approximable'.

We have the following lemma:
\begin{lem}\label{lem: well-defined derivative}
Suppose $w:B^n_1(0)\to\A_Q(\R^k)$, $x\in B^n_1(0)$ and that for each $j\in\{1,2\}$, $L^j_1,\dots,L^j_Q:\R^n\to\R^k$ are linear maps such that \eqref{eq: defn of differentiability} holds with each $L^j_i$ in place of $L^i_x$. Then there exists a permutation $\sigma\in S_Q$ such that $(w^i(x),L^1_i)=(w^{\sigma(i)}(x),L^2_{\sigma(i)})$ for all $i\in\{1,\dots,Q\}$.
\end{lem}

\textbf{Remark:} Lemma \ref{lem: well-defined derivative} follows, for instance, from \cite[Lemma 3.10]{Aiex26C^1Multi-Valued}. For completeness, we give an elementary proof of the lemma in Appendix \ref{app: Multi-Valued Functions}.

We introduce some further notation which arises from Lemma \ref{lem: well-defined derivative}:
\begin{defn}
Suppose $w:B^n_1(0)\to\A_Q(\R^k)$ is differentiable at $x\in B^n_1(0)$. By Lemma \ref{lem: well-defined derivative}, the $Q$-tuple $Dw(x):=\sum_{i=1}^Q\mb{L_x^i}\in\A_Q(\R^{n\times k})$ is well-defined. If $w$ is differentiable, we call the map $Dw:B^n_1(0)\to\A_Q(\R^{n\times k})$ the \emph{derivative} of $w$. We define the \emph{$1$-jet} of $w$ to be the function $j_1w:B^n_1(0)\to\A_Q(\R^k\times\R^{n\times k})$ defined by 
\begin{equation}\label{eq: defn of 1-jet}
j_1w(x)=\sum_{i=1}^Q\mb{(w^i(x),L^i_x)}.
\end{equation}
For $x\in B^n_1(0)$ and $h\in \R^n$, we write
\begin{equation}\label{eq: defn of action of 1-jet}
j_1w(x)[h]:=\sum_{i=1}^Q\mb{w^i(x)+L^i_x\cdot h}\in\A_Q(\R^k).
\end{equation}
By Lemma \ref{lem: well-defined derivative}, the expressions in \eqref{eq: defn of 1-jet} and \eqref{eq: defn of action of 1-jet} are well-defined.
\end{defn}

We now define $Q$-valued $C^1$ functions:
\begin{defn}\label{defn: C^1,alpha functions}
Suppose $w:B^n_1(0)\to\A_Q(\R^k)$ is differentiable. We say $w\in C^1(B^n_1(0);\A_Q(\R^k))$ if the function $j_1w$ is continuous. For $B^n_\sigma(y)\Subset B^n_1(0)$ we define
\begin{equation*}
|w|_{0;B^n_\sigma(y)}:=\sup_{x\in B^n_\sigma(y)}|w(x)|\quad\text{ and } \quad |w|_{1;B^n_\sigma(y)}:=|w|_{0;B^n_\sigma(y)}+|Dw|_{0;B^n_\sigma(y)}.
\end{equation*}
Let $\al\in(0,1]$. If $w\in C^1(B^n_1(0);\A_Q(\R^k))$, we say that $w\in C^{1,\al}(B^n_1(0);\A_Q(\R^k))$ if the function $Dw$ is $C^{0,\al}$. For $B^n_\sigma(y)\Subset B^n_1(0)$, we define
\begin{equation*}
[Dw]_{\al;B^n_\sigma(y)}:=\sup_{\underset{x\neq z}{x,z\in B^n_\sigma(y)}}\frac{\G(Dw(x),Dw(z))}{|x-z|^\al}\quad\text{ and }\quad |w|_{1,\al;B^n_\sigma(y)}:=|w|_{1;B^n_\sigma(y)}+[Dw]_{\al;B^n_\sigma(y)}.
\end{equation*}
We say that $w\in C^{2,\al}(B^n_1(0);\A_Q(\R^k))$ if $j_1w\in C^{1,\al}(B^n_1(0);\A_Q(\R^k\times\R^{n\times k}))$.
\end{defn}

Suppose $w\in C^1(B^n_1(0);\A_Q(\R^k))$. The single-valued $C^1$ function $w_a:B^n_1(0)\to\R^k$ defined by $w_a(x):=\frac{1}{Q}\sum_{\ell=1}^Qw^\ell(x)$ is called the \emph{average} of $w$. The function $w_f:B^n_1(0)\to\A_Q(\R^k)$ defined by $w_f(x):=\sum_{\ell=1}^Q\mb{w^\ell(x)-w_a(x)}$ is called the \emph{average-free part} of $w$. Note that $w_f$ is a $Q$-valued $C^1$ function. We also define the \emph{density $Q$ critical set} of $w$ to be the set
$$
\CK_w:=\{z\in B^n_1(0):j_1w(z)=Q\mb{j_1w_a(z)}\}.
$$ 
\begin{remark}\label{remark: Q-valued C^1,alpha functions} We digress to make a few remarks.
\begin{itemize}
\item[(a)] If $w_1,\dots,w_Q:B^n_1(0)\to\R^k$ are single-valued $C^{1,\al}$ functions, then the $Q$-valued function $w:=\sum_{i=1}^Q\mb{w_i}$ is a $Q$-valued $C^{1,\al}$ function and $j_1w\equiv\sum_{i=1}^Q\mb{(w_i,Dw_i)}$.
\item[(b)] If $w\in C^{1,\al}(B^n_1(0);\A_Q(\R^k))$, then $j_1w$ is $C^{0,\al}$ and we have $[j_1w]_{\al;B^n_\sigma(y)}\leq C|Dw|_{0,\al;B^n_\sigma(y)}$ for all $B^n_\sigma(y)\Subset B^n_1(0)$ where $C=C(n,k,Q)\in(0,\infty)$.
\item[(c)] By \cite[Theorem 5.5]{Aiex26C^1Multi-Valued}, if $w:B^n_1(0)\to\A_Q(\R^k)$ is differentiable and $Dw$ is continuous, then $w\in C^1$.
\end{itemize}
\end{remark}

Finally, we record the following separation principle (cf.~\cite[Proposition 2.6]{LucheseSchauder}):
\begin{prop}\label{prop: separation principle}
Suppose $Q\geq 2$ and $w\in C^1(B^n_1(0);\A_Q(\R^k))$. Then for every $x\in B^n_1(0)\setminus\CK_w$ there exists $\rho\in(0,1-|x|)$ and $Q_1,Q_2\in\{1,\dots, Q-1\}$ such that $Q_1+Q_2=Q$, and such that for each $i\in\{1,2\}$ there exists $w_i\in C^1(B^n_\rho(x);\A_{Q_i}(\R^k))$ such that $w\equiv w_1+w_2$ in $B^n_\rho(x)$.
\end{prop}
\begin{remark}\label{remark: separation principle}
By Proposition \ref{prop: separation principle}, it follows by induction on $Q$ that if $w:B^n_1(0)\to\A_Q(\R^k)$ is $C^1$, then $\CR_w$ is an open and dense subset of $B^n_1(0)$.
\end{remark}
Finally, we define $Q$-valued $W^{2,2}$ functions.
\begin{defn}\label{defn: Q-valued W^2,2}
Suppose $w\in W^{1,2}_{\textnormal{loc}}(B^n_1(0);\A_Q(\R^k))$. We say that $w\in W^{2,2}_{\textnormal{loc}}(B^n_1(0);\A_Q(\R^k))$ if $j_1w\in W^{1,2}_{\textnormal{loc}}(B^n_1(0);\A_Q(\R^k\times\R^{n\times k}))$.
\end{defn}
\textbf{Remark:} In Definition \ref{defn: Q-valued W^2,2}, the function $j_1w$ is the \textit{approximate $1$-jet} of $w$ defined as in \eqref{eq: defn of 1-jet} using the approximate derivative of $w$ (cf.~\cite[Corollary 2.7]{DLS11Qvalued}).

We will need the following property of $Q$-valued $W^{2,2}_{\textnormal{loc}}$ functions:
\begin{prop}\label{prop: W^2,2 + Dw cont implies C^1}
If $w\in W^{2,2}_\textnormal{loc}(B^n_1(0);\A_Q(\R^k))$ and $Dw\in C^0$, then $w\in C^1$.
\end{prop}
\begin{remark}
The example $v(x,y):=\mb{1-|x|}+\mb{|x|}$ shows that the above $W^{2,2}_\textnormal{loc}$ condition cannot be weakened to the $W^{1,\infty}_{\textnormal{loc}}$ condition. We leave the proof of Proposition \ref{prop: W^2,2 + Dw cont implies C^1} to Appendix \ref{app: Multi-Valued Functions}.
\end{remark}

Let $\CD^n$ denote the space of compactly supported $n$-forms in $\R^k\times B^n_1(0)$. We will need the following:
\begin{prop}\label{prop: current of Q-valued W^1,2 function}
Set $r:=\min\{n,k\}$. Let $f\in W^{1,r}_{\textnormal{loc}}(B^n_1(0);\A_Q(\R^k))$ and define $\bf{T}_f:\CD^n\to\R$ by 
$$
\bf{T}_f(\omega):=\int_{B^n_1(0)}\sum_{\ell=1}^Q\langle\omega(f^\ell(x),x),(Df^\ell(x),I_n)_{\#}(e_1\wedge\cdots\wedge e_n)\rangle\,\d x.
$$
Then $\bf{T}_f$ is a well-defined $n$-dimensional integer rectifiable current in $\R^k\times B^n_1(0)$. Furthermore, if $f$ is continuous, then we have $\spt\smallnorm{\bf{T}_f}=\graph(f)$.
\end{prop}
\textbf{Remark:} We leave the proof of the Proposition \ref{prop: current of Q-valued W^1,2 function} to Appendix \ref{app: Multi-Valued Functions}.
\section{Dirichlet Stationary Functions}\label{sec: W^2,2 outer stationary}
Fix $n\geq 2$ and $k,Q\geq 1$. Suppose $w\in W^{1,2}_{\textnormal{loc}}(B^n_1(0);\A_Q(\R^k))$. We say that $w$ is stationary for \emph{outer variations} if for every $\psi\in C^{1}_c(B^n_1(0)\times\R^k;\R^k)$ we have
\begin{equation*}
\int\sum_{\ell=1}^QDw^\ell(x):D(\psi(x,w^\ell(x)))=0.
\end{equation*}
We say that $w$ is stationary for \emph{inner variations} if for every $\zeta\in C^1_c(B^n_1(0);\R^n)$ we have
\begin{equation*}
\int\sum_{\ell=1}^Q\sum_{\kappa=1}^k(|Dw^{\ell,\kappa}|^2\delta_{ij}-2D_iw^{\ell,\kappa}D_jw^{\ell,\kappa})D_i\zeta_j=0.
\end{equation*}
We say that $w\in W^{1,2}_{\textnormal{loc}}(B^n_1(0);\A_Q(\R^k))$ is \emph{Dirichlet stationary} if it is stationary for both outer and inner variations.
\begin{remark}\label{remark: Literaturature specification}
Some authors refer to Dirichlet stationary functions as multi-valued harmonic functions. We stress that this notion is different to our notion of a $Q$-valued $C^1$ harmonic function.
\end{remark}

Our main result for functions that are stationary for outer variations is the following:
\begin{thm}\label{thm: reg of C^1 outer stat}
Suppose $w\in C^1(B^n_1(0);\A_Q(\R^k))$ is stationary for outer variations. Then $w$ is a $Q$-valued $C^1$ harmonic function.
\end{thm}
\begin{remark}
In particular, $w$ satisfies the conclusions of Theorems \ref{thm: derivative is Dir min}, \ref{thm: main-reg} and \ref{thm: rectifiability of branch set}. We leave the proof of Theorem \ref{thm: reg of C^1 outer stat} to Section \ref{sec: proof of reg of C^1 outer stat}.
\end{remark}

A key ingredient in the proof of Theorem \ref{thm: reg of C^1 outer stat} is the following regularity improvement for $W^{2,2}_{\textnormal{loc}}$ functions that are stationary for outer variations\footnote{The $C^{1,\al}$ regularity conclusion in Proposition \ref{prop: W^2,2 Dir stationary is C^1,alpha} has been established independently by Mattia Luchese using alternative methods (cf.~\cite[Corollary 4]{Luch26SmallLap}).}:
\begin{prop}\label{prop: W^2,2 Dir stationary is C^1,alpha}
Suppose $w\in W^{2,2}_{\textnormal{loc}}(B^n_1(0);\A_Q(\R^k))$ is stationary for outer variations. Then $Dw$ is locally Dirichlet minimising. In particular, there exists $\al=\al(n,Q)\in(0,1]$ with $\al(2,Q)=1/Q$ and $\al(n,2)=1/2$ such that $w\in C^{1,\al}(B^n_1(0);\A_Q(\R^k))$.
\end{prop}
\begin{remark}\label{remark: local Dir minimiser} 
The argument below shows that any function $G\in W^{1,2}_{\textnormal{loc}}(B^n_1(0);\A_Q(\R^{n\times k}))$ satisfying \eqref{eq: Hodge structure} is locally Dirichlet minimising.
\end{remark}
\begin{proof}
Set $G:=Dw$ and write $G=\sum_{\ell=1}^Q\mb{G^\ell}$ where $G^\ell=(G^\ell_1,\dots,G^\ell_n)$ for some $G^\ell_i\in\R^k$. It follows from \cite[Lemma~3.10]{LucheseSchauder} that
\begin{equation}\label{eq: Hodge structure}
D_iG_j^\ell=D_jG_i^\ell\quad\text{ and }\quad\sum_{i=1}^nD_iG_i^\ell=0,\quad\H^n\text{-a.e. in }B^n_1(0)\text{ for all }\ell,i,j.
\end{equation}
For an open set $U\Subset B^n_1(0)$ and a map $F\in W^{1,2}(U;\A_Q(\R^{n\times k}))$, we define the following energies
\begin{align*}
\CD_U(F)&:=\int_U |DF|^2,\qquad \A_U(F):=2\int_U\sum_{\ell=1}^Q\sum_{1\leq i<j\leq n}(D_iF_i^\ell\cdot D_jF_j^\ell-D_jF_i^\ell\cdot D_iF_j^\ell),\\
\H_U(F)&:=\int_U\sum_{\ell=1}^Q\left(\left|\sum_{i=1}^nD_iF_i^\ell\right|^2+\sum_{1\leq i<j\leq n}|D_iF_j^\ell-D_jF_i^\ell|^2\right).
\end{align*}
Note that for any square matrix $M$, writing $|M|:=\textnormal{tr}(MM^T)^{1/2}$, we have 
$$
|M|^2-\textnormal{tr}(M^2)=\frac{1}{2}|M-M^T|^2,
$$ 
and thus, we have that
\begin{equation}\label{eq: energies identity}
\CD_U(F)+\A_U(F)=\H_U(F).
\end{equation}
In particular, it follows from \eqref{eq: Hodge structure} and \eqref{eq: energies identity} that
\begin{equation}\label{eq: Hodge energy of G}
\H_U(G)=0\quad\text{ and }\quad\CD_U(G)=-\A_U(G),\text{ whenever }U\Subset B^n_1(0).
\end{equation}
Suppose $U\Subset B^n_1(0)$ is open and $H\in W^{1,2}(U;\A_Q(\R^{n\times k}))$ is such that $H\equiv Q\mb{0}$ outside a compact subset of $U$. Then \cite[Lemma~3.9]{LucheseSchauder} implies that
\begin{equation}\label{eq: consequence of paired Stokes}
\A_U(H)=0.
\end{equation}
Fix a ball $B\Subset B^n_1(0)$. Suppose that $F\in W^{1,2}(B;\A_Q(\R^{n\times k}))$ is such that $F|_{\partial B}=G|_{\partial B}$. Let $U$ be open with $B\Subset U\Subset B^n_1(0)$, and set
$$
\widetilde{F}:=
\begin{cases}
F,\text{ in }B,\\
G,\text{ in }U\setminus B.
\end{cases}
$$
Then $\widetilde{F}\in W^{1,2}(U;\A_Q(\R^{n\times k}))$. Let $\chi\in C_c^\infty(U)$ be such that $\chi=1$ in a neighbourhood of $\overline B$. By \eqref{eq: consequence of paired Stokes}, we have that $\A_U(\chi\widetilde{F})=\A_U(\chi G)=0$, which gives
\begin{equation}\label{eq:s commutor determined by trace}
\A_B(F)=\A_B(G).
\end{equation}
Combining \eqref{eq: energies identity}, \eqref{eq: Hodge energy of G}, and \eqref{eq:s commutor determined by trace}, we obtain
$$
\int_B|DF|^2-\int_B|DG|^2=\H_B(F)\geq 0.
$$
Since $F$ was arbitrary, $G=Dw$ is locally Dirichlet minimising. Hence, \cite[Theorem 2.13]{Almgren00} implies that $G\in C^{0,\al}(B^n_1(0);\A_Q(\R^{n\times k}))$, for some $\al=\al(n,Q)\in(0,1]$ with $\al(2,Q)=1/Q$ and $\al(n,2)=1/2$. It now follows from Proposition \ref{prop: W^2,2 + Dw cont implies C^1} that $w\in C^{1,\al}$, which completes the proof.
\end{proof}

\section{Proofs of the Main Theorems}

\subsection{Proofs of Theorems \ref{thm: derivative is Dir min}, \ref{thm: main-reg}, \ref{thm: rectifiability of branch set} and \ref{thm: Minkowski-type estimates in two-valued setting}}\label{sec: Proofs of reg thms for Q-valued C^1 harmonic functions}
Fix $n\geq 2$ and $k\geq 1$. We will work with the following inductive hypothesis for $Q\geq 1$:

\textbf{Hypothesis $(\star)_{Q}$:} Either $Q=1$, or $Q\geq 2$ and the following hold. If $q\in\{1,\dots, Q-1\}$ and $v\in C^1(B^n_1(0);\A_{q}(\R^k))$ is a $q$-valued $C^1$ harmonic function, then we have:
\begin{itemize}
\item[(i)] $\textnormal{cap}_2(\B_{v})=0$, i.e. $\B_{v}$ has vanishing $2$-capacity;
\item[(ii)] $v \in W^{2,2}_{\text{loc}}(B^n_1(0);\A_{q}(\R^k))$.
\end{itemize}

\textbf{Remark:} Hypothesis $(\star)_Q$ additionally depends on $n$ and $k$, and we have suppressed this dependence in the notation. 

We will establish the following:
\begin{prop}\label{prop: ind hyp always true}
Hypothesis $(\star)_Q$ holds for all $Q\geq 1$.
\end{prop}
We have a first lemma:

\begin{lem}\label{lem: separation harmonic assuming induction hyp}
Let $Q\geq 2$ and suppose Hypothesis $(\star)_{Q}$ holds. Let $w:B^n_1(0)\to\A_{Q}(\R^k)$ be a $Q$-valued $C^1$ harmonic function. Suppose $\Omega\subset B^n_1(0)$ is open, $Q_1,Q_2\in\{1,\dots, Q-1\}$ are such that $Q_1+Q_2=Q$, and $w|_\Omega\equiv w_1+w_2$, for some $w_i\in C^1(\Omega;\A_{Q_i}(\R^k))$. Then each $w_i$ is a $Q_i$-valued $C^1$ harmonic function.
\end{lem}
\begin{proof}
Without loss of generality, it suffices to show that $w_1$ is a $Q_1$-valued $C^1$ harmonic function. Let $z\in \Omega\setminus\B_{w_1}$. Then there exists an open ball $B\subset\Omega$ containing $z$ and $Q_1$ single-valued $C^1$ functions $u_1,\dots,u_{Q_1}:B\to\R^k$ such that $w_1|_B\equiv\sum_{i=1}^{Q_1}\mb{u_i}$. 

\textbf{Claim 1:} $w_2$ is locally harmonic in $B\setminus\B_{w_2}$ and each $u_i$ is harmonic in $B\setminus\B_{w_2}$. 

Fix $y\in B\setminus\B_{w_2}$. Then $y\in\Omega\setminus\B_w$, and there exists an open ball $\widetilde{B}\subset B\setminus\B_{w_2}$ containing $y$, $Q$ single-valued harmonic functions $\widetilde{w}_1,\dots,\widetilde{w}_Q:\widetilde{B}\to\R^k$, and, $Q_2$ single-valued $C^1$ functions $v_1,\dots, v_{Q_2}:\widetilde{B}\to\R^k$ such that $w|_{\widetilde{B}}\equiv\sum_{\ell=1}^Q\mb{\widetilde{w}_\ell}$ and $w_2|_{\widetilde{B}}\equiv\sum_{j=1}^{Q_2}\mb{v_j}$. Let $q\in\{1,\dots, Q\}$ denote the number of distinct functions in the collection $\{\widetilde{w}_1,\dots,\widetilde{w}_Q\}$. Note that, since $w$ is continuous, the function $\widetilde{B}\ni x\mapsto\text{card}(\spt(w(x)))$ is lower semi-continuous. 

Suppose first that $\Omega'$ is a non-empty open connected subset of $\widetilde{B}$ and $Q'\in\{1,\dots,Q\}$ such that $\text{card}(\spt(w))\equiv Q'$ in $\Omega'$. By the weak unique continuation for harmonic functions, and, by the constancy of $\text{card}(\spt(w))$ in $\Omega'$, for any indices $i,j\in\{1,\dots,Q\}$ we have that either $\widetilde{w}_i\equiv\widetilde{w}_j$ in $\widetilde{B}$, or $\widetilde{w}_i(x)\neq\widetilde{w}_j(x)$ for all $x\in\Omega'$. Hence, $Q'=q$, and for each $i\in\{1,\dots, Q_1\},j\in\{1,\dots,Q_2\}$, there exist $\ell_i,\ell_j\in\{1,\dots,Q\}$ such that $u_i|_{\Omega'}\equiv\widetilde{w}_{\ell_i}$ and $v_j|_{\Omega'}\equiv\widetilde{w}_{\ell_j}$. In particular, each function $u_i$ and $v_j$ is harmonic in $\Omega'$. 

We now assume without loss of generality that $q\geq 2$, and define the set
$$
\CK:=\{x\in\widetilde{B}:\widetilde{w}_i(x)=\widetilde{w}_j(x)\text{ and }D\widetilde{w}_i(x)=D\widetilde{w}_j(x)\text{ for some distinct functions }\widetilde{w}_i,\widetilde{w}_j\},
$$
and for every pair of indices $i,j\in\{1,\dots,Q\}$ such that $i<j$, define the (potentially empty) set
$$
\Sigma_{i,j}:=\{x\in\widetilde{B}\setminus{\CK}:\widetilde{w}_i(x)=\widetilde{w}_j(x)\text{ and }D\widetilde{w}_i(x)\neq D\widetilde{w}_j(x)\}.
$$
Note that $\CK$ is relatively closed in $\widetilde{B}$, and each $\Sigma_{i,j}$ is relatively closed in $\widetilde{B}\setminus\CK$. Define the function $m:\widetilde{B}\setminus\CK\to \Z_{\geq0}$ by setting
$$
m(x):=\text{card}(\{(i,j)\in\{1,\dots,Q\}^2: i<j\text{ and }x\in\Sigma_{i,j}\}).
$$
Note that $m$ is upper semi-continuous. We have the following claim:

\textbf{Claim 2:} For every $x\in\widetilde{B}\setminus\CK$ there exists $r_x>0$ such that each $u_i$ and $v_j$ is harmonic in $B^n_{r_x}(x)$. 

We prove Claim 2 by induction on $m(x)$. If $m(x)=0$, then, since $x\not\in\CK$, there exists $r_x>0$ such that $\text{card}(w)\equiv q$ in $B^n_{r_x}(x)$, and thus the claim follows by the above. Fix $m\in\Z_{\geq 1}$ and suppose that the claim has been established whenever $m(x)\leq m-1$. Suppose $z\in\widetilde{B}\setminus\CK$ satisfies $m(z)=m$. Let $\{(i_\ell,j_\ell)\}_{\ell=1}^m$ be the collection of indices such that $i_\ell<j_\ell$ and $z\in\Sigma_{i_\ell,j_\ell}$ for all $\ell\in\{1,\dots,m\}$. Since $m(z)=m$ and $z\not\in\CK$, there exists $r_z>0$ such that $m(x)\leq m-1$ for all $x\in B^n_{r_z}(z)\setminus\cap_{\ell=1}^m\Sigma_{i_\ell,j_\ell}$. Thus, the inductive hypothesis implies that each function $u_i$ and $v_j$ is harmonic in $B^n_{r_z}(z)\setminus\cap_{\ell=1}^m\Sigma_{i_\ell,j_\ell}$. In particular, each function $u_i$ and $v_j$ is harmonic in $B^n_{r_z}(z)\setminus\Sigma_{i_1,j_1}$. Furthermore, provided $r_z$ is sufficiently small, the set $B^n_{r_z}(z)\cap\Sigma_{i_1,j_1}$ is contained in a connected analytic embedded $(n-1)$-dimensional submanifold of $\widetilde{B}$ with no interior boundary and locally finite $\H^{n-1}$ measure. Since each $u_i$ and $v_j$ is $C^1$, it now follows from the divergence theorem and a standard extension argument that each $u_i$ and $v_j$ is harmonic in $B^n_{r_z}(z)$. 

We have thus established Claim 2, which implies that each function $u_i$ and $v_j$ is harmonic in $\widetilde{B}\setminus\CK$. Furthermore, standard unique continuation for harmonic functions gives $\dim_\H\CK\leq n-2$, and it follows from a standard excision argument that each $u_i$ and $v_j$ is harmonic in $\widetilde{B}$, which proves Claim 1. In particular, the weak unique continuation for harmonic functions implies that for all $i\in\{1,\dots, Q_1\},j\in\{1,\dots,Q_2\}$, there exist $\ell_i,\ell_j\in\{1,\dots,Q\}$ such that $u_i|_{\widetilde{B}}\equiv\widetilde{w}_{\ell_i}$ and $v_j|_{\widetilde{B}}\equiv\widetilde{w}_{\ell_j}$. 

Since $Q_2\in\{1,\dots, Q-1\}$ and Hypothesis $(\star)_{Q}$ is true, Claim 1 implies that $\dim_\H(\B_{w_2}\cap B)\leq n-2$. It thus follows from Claim 1 and a standard excision argument that each $u_i$ is harmonic in $B$. Since $z\in \Omega\setminus\B_{w_1}$ was arbitrary, it now follows that $w_1$ is a $Q_1$-valued $C^1$ harmonic function.
\end{proof}
\begin{remark}\label{remark: uniqueness of selection on R_w assuming ind hyp}
Let $Q\geq 2$ and suppose Hypothesis $(\star)_Q$ holds. Suppose $w:B^n_1(0)\to\A_Q(\R^k)$ is a $Q$-valued $C^1$ harmonic function. Lemma \ref{lem: separation harmonic assuming induction hyp} and Proposition \ref{prop: separation principle} imply that $\textnormal{cap}_2(\B_w\setminus\CK_w)=0$. Furthermore, by the argument given in the proof of Lemma \ref{lem: separation harmonic assuming induction hyp}, every point in $\CR_w$ is contained in an open ball in which $w$ has a unique decomposition by $C^1$ functions, up to reordering. In particular, any such functions must be harmonic, and thus analytic.
\end{remark}

We can now establish $W^{2,2}_{\textnormal{loc}}$ regularity for the average-free part:

\begin{lem}\label{lem:  W^2,2 reg for average-free part assuming ind hyp}
Let $Q\geq 1$ and suppose Hypothesis $(\star)_Q$ holds. Suppose $w:B^n_1(0)\to\A_{Q}(\R^k)$ is a $Q$-valued $C^1$ harmonic function. Then $w_f\in W^{2,2}_{\textnormal{loc}}(B^n_1(0);\A_Q(\R^k))$.
\end{lem}
\begin{remark}\label{remark: branch set of average-free part}
If $w:B^n_1(0)\to\A_Q(\R^k)$ is a $Q$-valued $C^1$ harmonic function, then, since $w_a$ is a single-valued $C^1$ function, $\B_w=\B_{w_f}$ and $w_f$ is a $Q$-valued $C^1$ harmonic function. 
\end{remark}
\begin{proof}
Assume without loss of generality that $Q\geq 2$. For $\eps\in(0,1)$, let $(\CK_{w_f})_\eps$ denote the $\eps$-neighbourhood of $\CK_{w_f}$. Since $Dw_f$ is locally harmonic in $B^n_1(0)\setminus\B_{w_f}$, we have $\Delta|Dw_f|^2=2|D^2w_f|^2$ on $B^n_1(0)\setminus\B_{w_f}$. Fix $\delta\in(0,1)$ and let $\gamma_\delta$ be a smooth nonnegative convex function on $\R$ such that $\gamma_\delta\equiv0$ in a neighbourhood of $0$ and $\gamma_\delta'(t)\equiv1$ on $[\delta,\infty)$. Let $\sigma\in(0,1)$ be such that $B^n_\rho(y)\Subset B^n_\sigma(0)$. Since $w_f\in C^1$,  there exists $\eps\in(0,1)$ such that 
$$
\spt(\gamma_\delta(|Dw_f|^2))\cap B^n_\sigma(0)\subseteq B^n_\sigma(0)\setminus(\CK_{w_f})_\eps.
$$
We thus have
\begin{equation}\label{eq: sub-harmonic equation for squared derivative assuming ind hyp}
\Delta\gamma_\delta(|Dw_f|^2)\geq 2\gamma_\delta'(|Dw_f|^2)|D^2w_f|^2,\quad\text{weakly in }B^n_\sigma(0)\setminus(\B_{w_f}\setminus(\CK_{w_f})_{\eps/2}).
\end{equation}
By Hypothesis $(\star)_Q$, Proposition \ref{prop: separation principle} and Lemma \ref{lem: separation harmonic assuming induction hyp}, we have that $\gamma_\delta(|Dw_f|^2)\in W^{1,2}_{\textnormal{loc}}(B^n_\sigma(0))$. Also, Remark \ref{remark: uniqueness of selection on R_w assuming ind hyp} implies that $\textnormal{cap}_2(\B_{w_f}\setminus(\CK_{w_f})_{\eps/2})=0$, and it thus follows from a standard excision argument that \eqref{eq: sub-harmonic equation for squared derivative assuming ind hyp} holds, weakly in $B^n_\sigma(0)$. Multiplying \eqref{eq: sub-harmonic equation for squared derivative assuming ind hyp} by a function $\zeta\in C^2_c(B^n_\rho(y))$ satisfying $0\leq\zeta\leq 1$ on $B^n_\rho(y)$, $\zeta\equiv 1$ on $B^n_{\rho/2}(y)$, $\zeta\equiv0$ outside $B^n_{3\rho/4}(y)$ and $|D^2\zeta|\leq C(n)\rho^{-2}$, integrating by parts, and letting $\delta\downarrow0$ gives
\begin{equation*}
\rho^{2-n}\int_{B^n_{\rho/2}(y)}|D^2w_f|^2\leq C(n)\rho^{-n}\int_{B^n_{3\rho/4}(y)}|Dw_f|^2,
\end{equation*}
It now follows from standard arguments that $w_f\in W^{2,2}_{\textnormal{loc}}$  (cf.~\cite[Lemma 3.5]{LucheseSchauder}).
\end{proof}

We have the following consequence of Lemma \ref{lem:  W^2,2 reg for average-free part assuming ind hyp} and Proposition \ref{prop: W^2,2 Dir stationary is C^1,alpha}:
\begin{lem}\label{lem: Dir stationarity of average-free part assuming induction hyp}
Let $Q\geq 1$ and suppose Hypothesis $(\star)_{Q}$ holds. If $w:B^n_1(0)\to\A_{Q}(\R^k)$ is a $Q$-valued $C^1$ harmonic function, then $w_f$ is Dirichlet stationary. In particular, $Dw_f$ is locally Dirichlet minimising.
\end{lem}
\begin{proof}
We first show that ${w_f}$ is stationary for outer variations. For $\delta>0$, let $\rho_\delta$ be a smooth function on $\R$ such that $\rho_\delta(-t)=\rho_\delta(t)$ for all $t$, $\rho_\delta(t)\equiv 1$ for all $t\geq\delta$, $\rho_\delta\equiv0$ in a neighbourhood of $0$, and $|\rho'_\delta|\leq 4/\delta$. Fix $\psi\in C^1_c(\R^k\times B^n_1(0);\R^k)$, and note that
\begin{equation*}
\int\sum_{\ell=1}^QDw_f^\ell(x):D(\psi(x,w^\ell_f(x)))=\lim_{\delta\downarrow0}\int\rho_\delta(|w_f|^2+|Dw_f|^2)\sum_{\ell=1}^QDw_f^\ell(x):D(\psi(x,w^\ell_f(x))).
\end{equation*}
Since $w_f$ is average-free, integrating by parts gives
\begin{align*}
\int\rho_\delta(|w_f|^2+|Dw_f|^2)\sum_{\ell=1}^QDw_f^\ell(x):D(\psi(x,w^\ell_f(x)))&=\\
&\hspace{-7.5cm}-\int\rho_\delta(|w_f|^2+|Dw_f|^2)\sum_{\ell=1}^Q\Delta w^\ell_f\psi(x,w_f^\ell(x))\\
&\hspace{-6.5cm}-2\int\rho'_\delta(|w_f|^2+|Dw_f|^2)\sum_{\ell=1}^Q(w_fDw_f+Dw_fD^2w_f):(Dw^\ell_f\cdot(\psi(x,w_f^\ell(x))-\psi(x,0))).
\end{align*}
Since $\sum_{\ell=1}^Q\Delta w^\ell_f\psi(x,w_f^\ell(x))=0$ $\H^n$-a.e. on $\{|w_f|^2+|Dw_f|^2>0\}$, we can apply Lemma \ref{lem:  W^2,2 reg for average-free part assuming ind hyp} to get
$$
\biggl|\int\rho_\delta(|w_f|^2+|Dw_f|^2)\sum_{\ell=1}^QDw_f^\ell(x):D(\psi(x,w^\ell_f(x)))\biggr|\leq\ C\int_{\pi_0(\spt\psi)}(\delta+\delta^{1/2}|D^2w_f|)\leq C\delta^{1/2},
$$
where $C\in(0,\infty)$ is independent of $\delta$ and $\pi_0:\R^k\times\R^n\to\R^n$ is the projection map. Taking $\delta\downarrow0$, we see that $w_f$ is stationary for outer variations. It now follows from Lemma \ref{lem:  W^2,2 reg for average-free part assuming ind hyp}, \cite[Theorem 3.11]{LucheseSchauder} and Proposition \ref{prop: W^2,2 Dir stationary is C^1,alpha} that $w_f$ is Dirichlet stationary and $Dw_f$ is locally Dirichlet minimising. This completes the proof.
\end{proof}

Finally, we have the following:
\begin{lem}\label{lem: branch set contained in that of derivative assuming ind hyp}
Let $Q\geq 1$ and suppose Hypothesis $(\star)_Q$ holds. If $w:B^n_1(0)\to\A_Q(\R^k)$ is a $Q$-valued $C^1$ harmonic function, then $\B_{w_f}=\B_{Dw_f}$.
\end{lem}
\begin{proof}
By Remark \ref{remark: branch set of average-free part}, we can assume without loss of generality that $w$ is average-free. Furthermore, Remark \ref{remark: uniqueness of selection on R_w assuming ind hyp} implies that $\B_{Dw}\subseteq\B_w$, and it thus remains to establish the reverse inclusion. We proceed by induction on $Q$. If $Q=1$, then the result is clear. Suppose $q\geq 2$ and that the result is proved whenever $Q\leq q-1$. Suppose $B\subseteq B^n_1(0)$ is an open ball and $H_1,\dots,H_{q}:B\to\R^{n\times k}$ are $C^1$ functions such that $Dw\equiv\sum_{j=1}^{q}\mb{H_j}$ in $B$. By Proposition \ref{prop: separation principle}, Lemma \ref{lem: separation harmonic assuming induction hyp} and the induction hypothesis, it suffices to assume that $B\cap\CK_w\neq\emptyset$. Define $d:B\to\R$ by setting $d(x):=\frac{1}{2}\min\{\dist(x,\partial B),\dist(x,\CK_w)\}$. Then for each $x\in B\setminus\CK_w$, the set $B^n_{d(x)}(x)$ is non-empty and simply connected.

\textbf{Claim:} For every $x\in B\setminus\CK_w$ there exists $r(x)\in(0,d(x))$ and $q$ analytic functions $w^x_1,\dots,w^x_{q}:B^n_{r(x)}(x)\to\R^k$ such that
$$
w\equiv\sum_{j=1}^{q}\mb{w^x_j}\quad\text{ and }\quad Dw^x_j\equiv H_j,\quad\text{ in }B^n_{r(x)}(x).
$$
\begin{proof}[Proof of Claim]
Fix $x\in B\setminus\CK_w$. By Proposition \ref{prop: separation principle} and Lemma \ref{lem: separation harmonic assuming induction hyp}, there exist $r(x)\in(0,d(x))$ and $Q_1,Q_2\in\{1,\dots,Q-1\}$, such that for each $i\in\{1,2\}$, there exists a $Q_i$-valued $C^1$ harmonic function $w_i:B^n_{r(x)}(x)\to\A_{Q_i}(\R^k)$ such that $w\equiv w_1+w_2$ in $B^n_{r(x)}(x)$. Also, since $w$ is average-free, Lemma \ref{lem: Dir stationarity of average-free part assuming induction hyp} implies that $Dw$ is locally Dirichlet minimising. In particular, each $H_i$ is an analytic single-valued harmonic function and $Dw|_B$ is a $Q$-valued $C^1$ harmonic function. It thus follows from the induction hypothesis and Remark \ref{remark: uniqueness of selection on R_w assuming ind hyp} (applied with $Dw|_B$) that it suffices to show that 
\begin{equation}\label{eq: claim branch equals that of derivative 1}
\B_{Dw_i}\text{ is empty for each }i\in\{1,2\}.
\end{equation}
Now, since $Dw$ is locally Dirichlet minimising, for any indexing set $I\subseteq\{1,\dots,Q\}$, the $|I|$-valued function $\sum_{i\in I}\mb{H_i}$ is locally Dirichlet minimising. Furthermore, Lemma \ref{lem: Dir stationarity of average-free part assuming induction hyp} implies that for each $i\in\{1,2\}$, the function $(Dw_i)_f$ is locally Dirichlet minimising. Since Hypothesis $(\star)_Q$ holds, we have that $\dim_\H\B_{w_i}\leq n-2$ for each $i\in\{1,2\}$, and it thus follows from \cite[Lemma 3.23(b)]{DLS11Qvalued} that each $Dw_i$ is locally Dirichlet minimising. Fix $y\in\CR_{w_1}\cap\CR_{w_2}$ (which exists by Proposition \ref{prop: separation principle}). Applying Remark \ref{remark: uniqueness of selection on R_w assuming ind hyp} with $Dw$, there exists $r_y>0$ and a permutation $\sigma_y\in S_Q$ such that 
\begin{equation}\label{eq: claim branch equals that of derivative 2}
Dw_1\equiv\sum_{\ell=1}^{Q_1}\mb{H_{\sigma_y(\ell)}}\quad\text{ and }\quad Dw_2\equiv\sum_{\ell=Q_1+1}^{Q}\mb{H_{\sigma_y(\ell)}},\quad\text{ in }B^n_{r_y}(y).
\end{equation} 
It now follows from \cite[Lemma 6.9]{DLMSVDirRectifiability} that \eqref{eq: claim branch equals that of derivative 2} holds in $B^n_{r(x)}(x)$, which establishes \eqref{eq: claim branch equals that of derivative 1} and proves the claim.
\end{proof}

Fix $j\in\{1,\dots,Q\}$ and define $w_j:B\to\R^k$ by setting
$$
w_j(x):=
\begin{cases}
w^x_j(x),\text{ if }x\in B\setminus\CK_w,\\
w_a(x),\text{ if }x\in\CK_w.
\end{cases}
$$
We note that $w_j$ is well-defined. Indeed, fix $x\in B\setminus\CK_w$, and let $U_x\subset B\setminus\CK_w$ be an open and simply connected set such that $x\in U_x$ and $\partial U_x\cap B\cap\CK_w\neq\emptyset$. By Remark \ref{remark: uniqueness of selection on R_w assuming ind hyp}, $w$ admits a unique decomposition into analytic functions in $U_x$, which, by the uniqueness of the decomposition into analytic functions in $B^n_{r(x)}(x)$, we can also label by $w^x_1,\dots,w^x_{q}:U_x\to\R^k$. Suppose that $j_1,j_2\in\{1,\dots,q\}$ such that $Dw^x_{j_1}\equiv Dw^x_{j_2}$ in $B^n_{r(x)}(x)$. It then follows by the unique continuation for analytic functions that $Dw^x_{j_1}\equiv Dw^x_{j_2}$ in $U_x$. Since $U_x$ is connected, $w^x_{j_1}\equiv w^x_{j_2}+C$ for some $C\in\R^k$. Since $\partial U_x\cap\CK_w\neq\emptyset$, it follows from continuity of $w$ that $C=0$. Thus $w^x_{j_1}(x)=w^x_{j_2}(x)$, and $w_j$ is well-defined. Finally, $w_j$ is $C^1$ and $w\equiv\sum_{i=1}^Q\mb{w_i}$ in $B$. This completes the proof.
\end{proof}

We can now prove Proposition \ref{prop: ind hyp always true}:
\begin{proof}[Proof of Proposition \ref{prop: ind hyp always true}]
Hypotheses $(\star)_1$ and $(\star)_2$ are clearly true. Suppose $Q\geq 2$ and Hypothesis $(\star)_{Q}$ holds. Suppose $w:B^n_1(0)\to\A_Q(\R^k)$ is a $Q$-valued $C^1$ harmonic function. By Remark \ref{remark: branch set of average-free part}, Lemmas \ref{lem: Dir stationarity of average-free part assuming induction hyp}, \ref{lem: branch set contained in that of derivative assuming ind hyp}, and Proposition \ref{prop: W^2,2 Dir stationary is C^1,alpha}, $\B_w$ coincides with the branch set of a Dirichlet minimiser. Thus, \cite{DLMSVDirRectifiability,KW17} imply that $\textnormal{cap}_2(\B_w)=0$. In particular $w_a$ is a harmonic single-valued function, and it thus follows from Lemma \ref{lem:  W^2,2 reg for average-free part assuming ind hyp} that $w\in W^{2,2}_{\text{loc}}(B^n_1(0);\A_Q(\R^k))$. Since $w$ was arbitrary, it follows that Hypothesis $(\star)_{Q+1}$ holds. Proposition \ref{prop: ind hyp always true} now follows by induction on $Q$.
\end{proof}

Theorem \ref{thm: derivative is Dir min} now follows from Proposition \ref{prop: ind hyp always true}, Lemmas \ref{lem: Dir stationarity of average-free part assuming induction hyp}, \ref{lem: branch set contained in that of derivative assuming ind hyp} and \cite[Lemma 3.23(b)]{DLS11Qvalued}. Theorems \ref{thm: main-reg} and \ref{thm: rectifiability of branch set} follow from Theorem \ref{thm: derivative is Dir min}, \cite[Theorem 0.12]{DLS11Qvalued}, \cite[Theorem 2.13]{Almgren00} and \cite{DLMSVDirRectifiability,KW17}. Theorem \ref{thm: Minkowski-type estimates in two-valued setting} follows from Theorem \ref{thm: derivative is Dir min}, Proposition \ref{prop: separation principle} and \cite[Theorem 1.2]{DLMSVDirRectifiability}.

\subsection{Proofs of Theorems \ref{thm: weak unique continuation} and \ref{thm: quantitative unique continuation}}
We can now prove Theorem \ref{thm: weak unique continuation}:
\begin{proof}[Proof of Theorem \ref{thm: weak unique continuation}]
We proceed by induction on $Q$. The result for $Q=1$ follows from the standard weak unique continuation for single-valued harmonic functions. Let $q\geq 2$ and suppose that the result has been established whenever $Q\leq q-1$. Suppose $w,v:B^n_1(0)\to\A_q(\R^k)$ are two $q$-valued $C^1$ harmonic functions which agree in a non-empty open subset of $B^n_1(0)$. By Theorem \ref{thm: rectifiability of branch set}, we can assume without loss of generality that $w$ and $v$ are average-free and non-vanishing. Define $\CK:=\CK_w\cup\CK_v$. By Theorem \ref{thm: rectifiability of branch set} and standard facts about harmonic functions, we have $\H^{n-1}(\CK)=0$. Thus, $\Omega:=B^n_1(0)\setminus\CK$ is open and connected. Since $\Omega$ is dense, it suffices to show that $w\equiv v$ in $\Omega$. Let $E$ be the closure in $\Omega$ of the interior of the set $\{j_1w=j_1v\}$. Then $E$ is non-empty and closed. Fix $x\in E$. Then $x\in \Omega$ and there exist $q_1,q_2\in\{1,\dots,q-1\}$, and for each $i\in\{1,2\}$, there exist $T_i\in\A_{q_i}(\R^k\times\R^{n\times k})$ such that $j_1w(x)=j_1v(x)=T_1+T_2$ and $\spt(T_1)\cap\spt(T_2)=\emptyset$. Hence, there exists $\delta>0$ such that if $S_i\in\A_{q_i}(\R^k\times\R^{n\times k})$ satisfy $\max\{\G(T_1,S_1),\G(T_2,S_2)\}\leq\delta$, then $\spt(S_1)\cap\spt(S_2)=\emptyset$. Since $w$ and $v$ are $C^1$, it follows that there exists an open ball $B\subset B^n_1(0)$ containing $x$, and for each $i\in\{1,2\}$, there exist two $q_i$-valued $C^1$ functions $w_i,v_i:B\to\A_{q_i}(\R^k)$ with $j_1w_i(x)=j_1v_i(x)=T_i$, such that $w|_B\equiv w_1+w_2$ and $v|_B\equiv v_1+v_2$. Hence, Proposition \ref{prop: ind hyp always true} and Lemma \ref{lem: separation harmonic assuming induction hyp} imply that each $w_i$ and $v_i$ is a $q_i$-valued $C^1$ harmonic function on $B$. Furthermore, since $x\in E$, there exists a non-empty open set $U\subseteq B$ such that $w_i|_U\equiv v_i|_U$ for each $i\in\{1,2\}$. By the induction hypothesis, $w\equiv v$ in $B$ and hence $B\subseteq E$. Since $x\in E$ was arbitrary, we conclude that $E$ is open. Finally, since $\Omega$ is connected, we have that $\Omega=E$, which completes the proof.
\end{proof}
Arguing analogously to the proof of Theorem \ref{thm: weak unique continuation}, using \cite[Lemma 5.3 and Theorem 5.7]{LucheseSchauder} in place of Proposition \ref{prop: ind hyp always true}, Lemma \ref{lem: separation harmonic assuming induction hyp} and Theorem \ref{thm: rectifiability of branch set}, we obtain the following:
\begin{thm}\label{thm: weak unique continuation for stationary C^1,alpha graphs}
Let $n,Q\geq2$, $k\geq 1$ and $\al\in(0,1)$. Suppose $w,v\in C^{1,\al}(B^n_1(0);\A_Q(\R^k))$ are such that the integral varifolds $\bf{v}(w)$ and $\bf{v}(v)$ are stationary in $\R^k\times B^n_1(0)$. If $U\subseteq B^n_1(0)$ is non-empty and open such that $w\equiv v$ in $U$, then $w\equiv v$ in $B^n_1(0)$.
\end{thm}
We now prove Theorem \ref{thm: quantitative unique continuation}:
\begin{proof}[Proof of Theorem \ref{thm: quantitative unique continuation}]
Suppose $w,v:B^n_1(0)\to\A_Q(\R^k)$ are $Q$-valued $C^1$ harmonic functions. Let $\B:=\B_w\cup\B_v$ and set $\Omega:=B^n_1(0)\setminus\B$. By Theorem \ref{thm: rectifiability of branch set}, we have
\begin{equation}\label{eq: dim bound of B}
\dim_\H\B\leq n-2.
\end{equation}
We establish (1), the proof of (2) is analogous. Define $
A:=\{w=v\}$, and suppose that $\dim_\H A>n-1$. Then \eqref{eq: dim bound of B} implies that there exists an open ball $B\subset\Omega$ such that $\dim_\H(A\cap B)>n-1$. Since $B$ is simply connected, there exist $2Q$ single-valued harmonic functions $w_1,\dots,w_Q,v_1,\dots,v_Q:B\to\R^k$ such that $w|_B\equiv\sum_{i=1}^Q\mb{w_i}$ and $v|_B\equiv\sum_{i=1}^Q\mb{v_i}$. For each permutation $\sigma\in S_Q$, define 
$$ \Sigma_\sigma:=\cap_{i=1}^Q\{w_i=v_{\sigma(i)}\}.
$$
Then $A\cap B=\cup_{\sigma\in S_Q}\Sigma_\sigma$, and thus there exists $\sigma_0\in S_Q$ such that $\dim_\H\Sigma_{\sigma_0}>n-1$. By standard results for harmonic functions, we have that $\Sigma_{\sigma_0}=B$. Hence, $w\equiv v$ in $B$ and it follows from Theorem \ref{thm: weak unique continuation} that $w\equiv v$ in $B^n_1(0)$. This completes the proof.
\end{proof}

\subsection{Proofs of Theorem \ref{thm: real part of Dir min} and Corollary \ref{cor: Full Series Expansion in Disc}}
We now prove Theorem \ref{thm: real part of Dir min}:

\textit{Proof of Theorem \ref{thm: real part of Dir min}.} By Theorem \ref{thm: rectifiability of branch set}, after translating and rescaling, it suffices to prove the following claim:
\begin{claim}\label{claim: real part of Dir min}
Suppose $w\in L^\infty_{\textnormal{loc}}(\D;\A_Q(\R^k))$ is locally harmonic in $\D\setminus\{0\}$. Suppose further that $w$ does not decompose in $\D$ as the sum of two $L^\infty_{\textnormal{loc}}$ ($<Q$)-valued functions which are locally harmonic in $\D\setminus\{0\}$. Then there exists a holomorphic function $F:\D\to\C^k$ such that the function $g:\D\to\A_Q(\C^k)$ defined by
\begin{equation}\label{eq: real part of Dir min}
g(z)=\sum_{\zeta^Q=z}\mb{F(\zeta)}
\end{equation}
is locally Dirichlet minimising, and $w\equiv \Re (g)$ in $\D$.
\end{claim}
\begin{proof}[Proof of Claim \ref{claim: real part of Dir min}] Let $\Omega:=\D\setminus\Gamma$, where $\Gamma:=\{z\in\D:\Re(z)\geq0,\Im(z)=0\}$. We also set $D^\pm:=\{z\in\D:\pm\Im(z)>0\}$. Since $\Omega$ is simply connected, there exist $Q$ single-valued harmonic functions $w_1,\dots,w_Q:\Omega\to\R^k$ such that $w|_\Omega\equiv\sum_{i=1}^Q\mb{w_i}$. Set $w_{Q+1}:=w_1$. Then, as $w$ does not decompose in $\D$ as the sum of two $L^\infty_{\textnormal{loc}}$ ($<Q$)-valued functions which are locally harmonic in $\D\setminus\{0\}$, after relabelling, we have
\begin{equation*}
\lim_{z\in\D^-:z\to\omega}w_i(z)=\lim_{z\in\D^+:z\to\omega}w_{i+1}(z),\quad\text{ for all }\omega\in\Gamma.
\end{equation*}
For each $i\in\{1,\dots,Q\}$, denote the open set $U_i:=\{\zeta\in\D\setminus\{0\}:2\pi(i-1)/Q<\arg(\zeta)<2\pi i/Q\}$. Define $u:\cup_{i=1}^QU_i\to\R^k$ by setting $u(\zeta):=w_i(\zeta^Q)$ for all $\zeta\in U_i$, and extend $u$ continuously to $\D$. It follows from the 2-dimensional conformal invariance of Dirichlet energy and standard removability that $u$ is harmonic and satisfies
$$
w(z)=\sum_{\zeta^Q=z}\mb{u(\zeta)},\quad\text{ for all }z\in\D.
$$
Let $v:\D\to\R^k$ be a harmonic function such that $F:=u+iv$ is holomorphic. Define $g:\D\to\A_Q(\C^k)$ as in \eqref{eq: real part of Dir min}. Then by \cite[Theorem 2.20]{Almgren00} (cf.~\cite[Theorem 0.1]{Spadaro2010ComplexVarieties}), $g$ is Dirichlet minimising. This establishes Claim \ref{claim: real part of Dir min}, and completes the proof.
\end{proof}

We can now prove Corollary \ref{cor: Full Series Expansion in Disc}:

\begin{proof}[Proof of Corollary \ref{cor: Full Series Expansion in Disc}]
Let $F_1,\dots,F_N:\D\to\C^k$ be the holomorphic functions given by Theorem \ref{thm: real part of Dir min}. Fix $j\in\{1,\dots,N\}$ and write
$$
F_j(\zeta)=\sum_{m=0}^{\infty}a_{j,m}\zeta^m,\quad\text{ for }\zeta\in\D.
$$
Since $Dw$ is locally bounded, the chain rule and the Cauchy–Riemann equations imply that $\zeta^{1-Q_j}F_j'(\zeta)$ is bounded in $\D\setminus\{0\}$. It thus follows from the Cauchy integral formula that $a_{j,m}=0$ for all $m:1\leq m<Q_j$. Finally, $h_j(z):=\operatorname{Re}
(\sum_{\ell=0}^{\infty}a_{j,\ell Q_j}z^\ell)$ is the real part of a holomorphic function and is thus harmonic. This completes the proof.
\end{proof}

\subsection{Proof of Theorem \ref{thm: reg of C^1 outer stat}}\label{sec: proof of reg of C^1 outer stat}

We can now prove Theorem \ref{thm: reg of C^1 outer stat}:
\begin{proof}[Proof of Theorem \ref{thm: reg of C^1 outer stat}]
We proceed by induction on $Q$. The result is clear when $Q=1$. Fix $q\geq 2$ and assume that the result has been established whenever $Q\leq q-1$. Suppose $w\in C^1(B^n_1(0);\A_q(\R^k))$ is stationary for outer variations. Since $w_a$ is a single-valued harmonic function, we can assume without loss of generality that $w$ is average-free. By Proposition \ref{prop: separation principle}, \cite[Lemma 3.2]{LucheseSchauder} and the induction hypothesis, $w$ is a $q$-valued $C^1$ harmonic function in $B^n_1(0)\setminus\CK_w$. In particular, Theorem \ref{thm: rectifiability of branch set} implies that $\textnormal{cap}_2(\B_w\setminus\CK_w)=0$. We can now argue as in the proof of Lemma \ref{lem:  W^2,2 reg for average-free part assuming ind hyp}, using \cite[Lemma 3.2]{LucheseSchauder} in place of Lemma \ref{lem: separation harmonic assuming induction hyp} to see that $w\in W^{2,2}_{\textnormal{loc}}(B^n_1(0);\A_q(\R^k))$. Hence, Proposition \ref{prop: W^2,2 Dir stationary is C^1,alpha} implies that $Dw$ is Dirichlet minimising. Note that if $Dw\equiv q\mb{0}$, then \cite[Proposition 2.12]{DLS11Qvalued} implies that $w\equiv\textnormal{const}$ which establishes the claim. We can thus assume without loss of generality that $Dw\not\equiv q\mb{0}$. Hence, since $Dw$ is average-free, it follows from \cite[Proposition 3.22]{DLS11Qvalued} that we have $\dim_\H\{Dw=q\mb{0}\}\leq n-2$. Since $\CK_w\subseteq \{Dw=q\mb{0}\}$, we can now conclude that $\dim_\H\B_w\leq n-2$. It now follows from Proposition \ref{prop: ind hyp always true}, Remark \ref{remark: uniqueness of selection on R_w assuming ind hyp} and the induction hypothesis that $w$ is a $q$-valued $C^1$ harmonic function. This completes the proof.
\end{proof}

\subsection{Proof of Theorem \ref{thm: Uniqueness of Tangent Maps}}\label{sec: Proof of Uniqueness of Tangent Maps}
It follows from an induction argument on $Q$, using Propositions \ref{prop: separation principle}, \ref{prop: ind hyp always true} and Lemma \ref{lem: separation harmonic assuming induction hyp} that Theorem \ref{thm: Uniqueness of Tangent Maps} follows from the following:
\begin{thm}\label{thm: Density Q uniqueness of tangent maps}
Suppose $w:B^n_1(0)\to\A_Q(\R^k)$ is a $Q$-valued $C^1$ harmonic function such that $w\not\equiv Q\mb{w_a}$ in $B^n_1(0)$. For $\H^{n-2}$-a.e. $z\in\CK_w$, there exists $\rho_z\in(0,1-|z|)$, an integer $N_z\geq 0$, relatively prime positive integers $q_z\geq 1,k_z\geq q_z+1$ with $N_zq_z\leq Q$, a unique non-zero, homogeneous, average-free $Q$-valued $C^1$ harmonic function $\phi^{(z)}:\R^n\to\A_Q(\R^k)$ given by
\begin{equation}\label{eq: cylindrical form of tangent-map}
\phi^{(z)}(x)=\sum_{j=1}^Q\mb{\phi_j^{(z)}}=(Q-N_zq_z)\mb{0}+\sum_{\ell=1}^{N_z}\textnormal{Re}(c_\ell^{(z)}(x_1+ix_2)^{k_z/q_z})
\end{equation}
for some $c_\ell^{(z)}\in\C^k\setminus\{0\}$, and an orthogonal rotation $Q_z$ of $\R^n$ such that
$$
w(z+x)=\sum_{j=1}^Q\mb{w_a(z+x)+\phi_j^{(z)}(Q_zx)+\eps_j^{(z)}(x)},\quad\text{ for all }x\in B^n_{\rho_z}(0),
$$
where, for each $x\in B^n_{\rho_z}(0)$, $\eps^{(z)}_j(x)\in\R^k$ such that $\sum_{j=1}^Q\eps_j^{(z)}(x)=0$ and
$$
\sigma^{-n}\int_{B^n_\sigma(0)}\sum_{j=1}^Q|\eps_j^{(z)}|^2= o(\sigma^{2k_z/q_z}),\quad\text{ for all }\sigma\in(0,\rho_z).
$$
\end{thm}

Before proceeding to the proof of Theorem \ref{thm: Density Q uniqueness of tangent maps}, we need to record a few facts about $Q$-valued $C^1$ harmonic functions that we can now prove. Let $w:B^n_1(0)\to\A_Q(\R^k)$ be a $Q$-valued $C^1$ harmonic function. By the results of Section \ref{sec: Proofs of reg thms for Q-valued C^1 harmonic functions}, we have that $w\in C^{1,\al}(B^n_1(0);\A_Q(\R^k))\cap W^{2,2}_{\textnormal{loc}}(B^n_1(0);\A_Q(\R^k))$ is Dirichlet stationary and $Dw$ is Dirichlet minimising. In particular, by \cite[Theorem 2.13]{Almgren00}, for every ball $B^n_\rho(y)\Subset B^n_1(0)$, $w$ satisfies the estimate
\begin{equation}\label{eq: C^1,al+W^2,2 estimates}
|w|^2_{0;B^n_{\rho/2}(y)}+\rho^2|Dw|^2_{0;B^n_{\rho/2}(y)}+\rho^{2+2\al}[Dw]^2_{\al;B^n_{\rho/2}(y)}+\rho^{4-n}\int_{B^n_{\rho/2}(y)}|D^2w|^2\leq C\rho^{-n}\int_{B^n_\rho(y)}|w|^2,
\end{equation}
where $C=C(n,Q)\in(0,\infty)$ and $\al=\al(n,Q)\in(0,1)$. Suppose without loss of generality that $w$ is average-free and $w\not\equiv Q\mb{0}$ in $B^n_1(0)$. Fix $z\in\CK_w$ and for $\sigma\in(0,1-|z|)$, define $w_{z,\sigma}:B^n_{\sigma^{-1}(1-|z|)}(0)\to\A_Q(\R^k)$ by setting
$$
w_{z,\sigma}(x):=(\sigma^{-n/2}\smallnorm{w}_{L^2(B^n_\sigma(z))})^{-1}w(z+\sigma x).
$$
Then by \eqref{eq: C^1,al+W^2,2 estimates} and \cite[Lemmas 2.5 and 2.10]{LucheseSchauder}, for any sequence $\sigma_j\downarrow0$, there exists a function $v\in C^{1,\al}(B^n_1(0);\A_Q(\R^k))\cap W^{2,2}_{\textnormal{loc}}(B^n_1(0);\A_Q(\R^k))$ which is stationary for outer variations, such that, after taking a subsequence and relabelling, $(w_{z,\sigma_j})_j$ converges to $v$ in the local $Q$-valued $C^1$ topology. In particular, Theorem \ref{thm: reg of C^1 outer stat} implies that $v$ is a $Q$-valued $C^1$ harmonic function. A $Q$-valued $C^1$ harmonic function $v$ arising as above is called a \textit{tangent map} to $w$ at $z$.

In fact, $v$ is Dirichlet stationary (cf.~\cite[Theorem 3.11]{LucheseSchauder} or Section \ref{sec: Proofs of reg thms for Q-valued C^1 harmonic functions}), and it follows from standard frequency function techniques (cf.~\cite[Chapter 2]{Almgren00} or \cite[Chapter 3]{DLS11Qvalued}) that $v$ is homogeneous of degree $N_w(z)$, where $N_w(z):=\lim_{\rho\downarrow0}N_{w,z}(\rho)$ and $N_{w,z}:(0,1-|z|)\to \R$ is the monotone non-decreasing function defined by
$$
N_{w,z}(\rho):=\frac{\rho^{2-n}\int_{B^n_\rho(z)}|Dw|^2}{\rho^{1-n}\int_{\partial B^n_\rho(z)}|w|^2}.
$$
Moreover, the set
$$
S(v):=\{y\in\R^n:v(x+y)=v(x)\text{ for all }x\in\R^n\}
$$
is a proper linear subspace of $\R^n$. In fact, since $w\in C^{1,\al}$ and $z\in\CK_w$, we have that $N_w(z)\geq1+\al>1$, and thus $\dim S(v)\leq n-2$. We denote by $\Sigma_w$ the set of points of $\CK_w$ at which there exists a tangent map to $w$ whose spine dimension is equal to $n-2$. It now follows from standard geometric measure theory techniques (cf.~\cite[Section 3.4]{Simon96}) that $\H^{n-2}(\CK_w\setminus\Sigma_w)=0$.

We can now prove Theorem \ref{thm: Density Q uniqueness of tangent maps}:
\begin{proof}[Proof of Theorem \ref{thm: Density Q uniqueness of tangent maps}]
By Corollary \ref{cor: Full Series Expansion in Disc}, we can assume without loss of generality that $n\geq 3$. By Section \ref{sec: Proofs of reg thms for Q-valued C^1 harmonic functions}, we can assume without loss of generality that $w$ is average-free. Furthermore, $Dw$ is Dirichlet minimising, and thus \cite[Theorem A]{KW17} implies that $Dw$ has a unique tangent map at $\H^{n-2}$-a.e. point of $\CK_w$. Let $\Sigma'_w$ denote the set of points of $z\in\Sigma_w$ at which $Dw$ has a unique tangent map. Note that $\H^{n-2}(\Sigma_w\setminus\Sigma'_w)=0$. We have a first claim:

\begin{claim}\label{claim: uniqueness of tangent maps}
If $z\in\Sigma'_w$, then $w$ has a unique tangent map at $z$.
\end{claim}
\begin{proof}[Proof of Claim \ref{claim: uniqueness of tangent maps}]
Fix $z\in\Sigma'_w$, and suppose that $u,v:B^n_1(0)\to\A_Q(\R^k)$ are tangent maps to $w$ at $z$. By the monotonicity of $N_{w,z}$ and the doubling estimates given by the frequency function  of $w$, for every $0<\sigma<\rho<(1-|z|)/2$ we have
$$
N_w(z)(\rho/\sigma)^{2N_w(z)}\sigma^{1-n}\int_{\partial B^n_\sigma(z)}|w|^2\leq \rho^{2-n}\int_{B^n_\rho(z)}|Dw|^2\leq N_{w,z}(\rho)(\rho/\sigma)^{2N_{w,z}(\rho)}\sigma^{1-n}\int_{\partial B^n_\sigma(z)}|w|^2.
$$
Integrating in $\sigma$ gives that for every $\rho\in(0,(1-|z|)/2)$ we have 
$$
N_w(z)(n+2N_w(z))\rho^{-n}\int_{B^n_\rho(z)}|w|^2\leq\rho^{2-n}\int_{B^n_\rho(z)}|Dw|^2\leq N_{w,z}(\rho)(n+2N_{w,z}(\rho))\rho^{-n}\int_{B^n_\rho(z)}|w|^2.
$$
In particular, since $N_{w,z}(\rho)\geq N_w(z)\geq 1+\al>0$, there exists $c_z>0$ such that both $c_zDu$ and $c_zDv$ are tangent maps to $Dw$ at $z$. By definition of $\Sigma'_w$, we must have that $Du=Dv$. Since $u$ and $v$ are $C^1$ and homogeneous of the same degree, we have that $u=v$. This completes the proof.
\end{proof}
The result now follows from \cite[Theorem A]{KW17}.
\end{proof}

\subsection{Proofs of Theorems \ref{thm: Lip+W^2,2 Stationary + Lip small implies C^1,alpha}, \ref{thm: Bel25 Theorem 5 reg improvement} and \ref{thm: C^1+W^2,2 implies C^1,alpha}}
In this section, we give proofs of Theorems \ref{thm: Lip+W^2,2 Stationary + Lip small implies C^1,alpha}, \ref{thm: Bel25 Theorem 5 reg improvement} and \ref{thm: C^1+W^2,2 implies C^1,alpha}. We begin by adapting the methods of Section \ref{sec: W^2,2 outer stationary} to prove Theorem \ref{thm: Lip+W^2,2 Stationary + Lip small implies C^1,alpha}:
\begin{proof}[Proof of Theorem \ref{thm: Lip+W^2,2 Stationary + Lip small implies C^1,alpha}]
Define $g:\R^{n\times k}\to\R^{n\times n}$ by setting
\begin{equation}\label{eq: defn of g}
g(X):=(I_n+X^TX)^{-1}.
\end{equation}
Set $F:=Du$ and write $F=(F_1,\dots,F_n)$ where each $F_i\in W^{1,2}_{\textnormal{loc}}(B^n_1(0);\A_Q(\R^k))$. Arguing as in the proof of \cite[Lemma 3.10]{LucheseSchauder}, for every $\ell,\kappa$, we have
\begin{equation}\label{eq: Q-valued MSS PDE}
g(F^\ell):DF^{\ell,\kappa}=0\quad\H^n\text{-a.e. in }B^n_1(0).
\end{equation}  
For $P=(P^1,\dots,P^k)\in\R^{n^2\times k}$, define
$$
N(P):=\sum_{\kappa=1}^k[(\textnormal{tr}P^\kappa)^2-\textnormal{tr}((P^\kappa)^2)]=2\sum_{\kappa=1}^k\sum_{i<j}(P^\kappa_{ii}P^\kappa_{jj}-P^\kappa_{ij}P^\kappa_{ji}).
$$
Note that, arguing as in the proof of Proposition \ref{prop: W^2,2 Dir stationary is C^1,alpha}, for every ball $B^n_r(y)\Subset B^n_1(0)$ and every $G\in W^{1,2}(B^n_r(y);\A_Q(\R^{n\times k}))$ satisfying $G=F$ on $\partial B^n_r(y)$, we have
$$
\int_{B^n_r(y)}\sum_{\ell=1}^Q N(DG^\ell)=\int_{B^n_r(y)}\sum_{\ell=1}^Q N(DF^\ell).
$$
Define $R:\R^{n\times k}\times\R^{n^2\times k}\to\R$ by setting
$$
R(X,P):=\sum_{\kappa=1}^k(g(X):P^\kappa)^2+\sum_{\kappa=1}^k\sum_{i<j}(P^\kappa_{ij}-P^\kappa_{ji})^2.
$$
Since $u\in W^{2,2}_{\textnormal{loc}}(B^n_1(0);\A_Q(\R^k))$, \cite[Lemma 3.10]{LucheseSchauder} implies that each $DF^\ell$ is symmetric. Combining with \eqref{eq: Q-valued MSS PDE} gives $R(F^\ell,DF^\ell)=0$ $\H^n$-a.e. in $B^n_1(0)$. Note that
\begin{equation}\label{eq: martix symmetry property}
\sum_{\kappa=1}^k\sum_{i<j}(P^\kappa_{ij}-P^\kappa_{ji})^2=|P|^2-\sum_{\kappa=1}^k\textnormal{tr}((P^\kappa)^2),\quad\text{ for all }P=(P^1,\dots,P^k)\in\R^{n^2\times k}.
\end{equation}
Setting $J(X,P):=R(X,P)-N(P)$, we thus have
\begin{equation}\label{eq: equation for energy}
J(X,P)=|P|^2+\sum_{\kappa=1}^k[(g(X):P^\kappa)^2-(\textnormal{tr}P^\kappa)^2].
\end{equation}
Fix an open ball $B^n_R(y)\Subset B^n_1(0)$. For all $r\in(0,R)$ such that $F|_{\partial B^n_r(y)}\in W^{1,2}(\partial B^n_r(y);\A_Q(\R^{n\times k}))$, define
$$
\mathscr{J}(G;B^n_r(y)):=\int_{B^n_r(y)}\sum_{\ell=1}^Q J(G^\ell,DG^\ell),\quad\text{ for }G\in W^{1,2}(B^n_r(y);\A_Q(\R^{n\times k})).
$$
If $G=F$ on $\partial B^n_r(y)$, then the above gives
\begin{equation*}
\mathscr{J}(G;B^n_r(y))-\mathscr{J}(F;B^n_r(y))=\int_{B^n_r(y)}\sum_{\ell=1}^Q R(G^\ell,DG^\ell)\geq 0.
\end{equation*}
Thus $F$ is a local minimiser of $\mathscr{J}$. Fix $\eps\in(0,1/2)$ to be determined depending on $n$ and $Q$, and suppose that $L=\Lip(u)\leq\eps/\sqrt{n}$. Then for all $\ell$, we have that $|F^\ell|\leq \eps$ $\H^n$-a.e. in $B^n_1(0)$. For $X\in B^{n\times k}_\eps(0)$ and all $P^\kappa\in\R^{n\times n}$, we have
\begin{equation*}
|(g(X):P^\kappa)^2-(\textnormal{tr}P^\kappa)^2|=|(g(X)-I_n):P^\kappa||(g(X)+I_n):P^\kappa|\leq C(n)\eps^2|P^\kappa|^2.
\end{equation*}
In particular, by \eqref{eq: equation for energy}, for all $X\in B^{n\times k}_\eps(0)$ and all $P\in\R^{n^2\times k}$, we have
$$
(1-C(n)\eps^2)|P|^2\leq J(X,P)\leq (1+C(n)\eps^2)|P|^2.
$$
Thus, for any $G\in W^{1,2}(B^n_r(y);\A_Q(B^{n\times k}_\eps(0)))$ satisfying $G=F$ on $\partial B^n_r(y)$, it follows from the above established local minimisation property of $F$ that
\begin{equation}\label{eq: codim-one almost Dir min}
\int_{B^n_r(y)}|DF|^2\leq\Lambda\int_{B^n_r(y)}|DG|^2,\quad\text{ where }\Lambda:=(1+C(n)\eps^2)/(1-C(n)\eps^2).
\end{equation}
Let $C_\star=C_\star(n,Q)\in(0,\infty)$ be the constant given by \cite[Proposition 3.10]{DLS11Qvalued}. Then $C_\star(2,Q)=Q$, and $C_\star(n,Q)<1/(n-2)$ for all $n\geq 3$. If $n=2$, set $\eps\in(0,1/2)$ such that $\Lambda=2$. If $n\geq3$, choose $\eps=\eps(n,Q)\in(0,1)$ such that $\Lambda C_\star<1/(n-2)$. For $\H^1$-a.e. $r$ with $B^n_r(y)\Subset B^n_1(0)$, we have $F|_{\partial B^n_r(y)}\in W^{1,2}(\partial B^n_r(y);\A_Q(\R^{n\times k}))$. Fix such an $r$ and let $G$ be a Dirichlet-minimising extension of $F|_{\partial B^n_r(y)}$. Note that we can ensure that $G$ takes values in $\A_Q(B^{n\times k}_\eps(0))$, as composing each component of $G$ with the nearest-point projection onto $B^{n\times k}_\eps(0)$ preserves the trace and does not increase Dirichlet energy. By \eqref{eq: codim-one almost Dir min} and \cite[(3.19)]{DLS11Qvalued}, we have
\begin{equation}
\int_{B^n_r(y)}|DF|^2\leq C_\star\Lambda\,r\int_{\partial B^n_r(y)}|DF|^2,\quad\text{ for }\H^1\text{-a.e. }r\in(0,R).
\end{equation}
Arguing as in the proof of \cite[Theorem 3.9]{DLS11Qvalued} now gives that $F\in C^{0,\alpha}(B^n_{1/2}(0))$, where $\al=\al(n,Q)\in(0,1)$. Finally, Proposition \ref{prop: W^2,2 + Dw cont implies C^1} gives $u\in C^{1,\al}(B^n_{1/2}(0);\A_Q(\R^k))$, which completes the proof.
\end{proof}
Theorem \ref{thm: Bel25 Theorem 5 reg improvement} follows from \cite[Theorem 5]{Bel25} and the following consequence of Theorem \ref{thm: Lip+W^2,2 Stationary + Lip small implies C^1,alpha}:
\begin{prop}\label{prop: proof of Bel25 Theorem 5 reg improvement}
Fix $Q\geq 1$ and $n\geq 2$. There exists $\eps=\eps(n,Q)\in(0,1)$ such that the following holds. Suppose $u:B^n_1(0)\to\A_Q(\R)$ is Lipschitz, and that the integral varifold $V:=\bf{v}(u)$ is stationary in $\R\times B^n_1(0)$. Suppose further that there exists a closed set $\Sigma\subset \R\times B^n_1(0)$ with $\textnormal{cap}_2(\Sigma\cap\spt\smallnorm{V})=0$ such that $V$ is represented in $(\R\times B^n_1(0))\setminus\Sigma$ by a properly immersed, two-sided, stable minimal hypersurface. If $\Lip(u)\leq\eps$, then $u\in C^{1,\al}(B^n_{1/2}(0);\A_Q(\R))$ for some $\al=\al(n,Q)\in(0,1]$ with $\al(2,Q)=1/Q$ and $\al(n,2)=1/2$.
\end{prop}
\begin{proof}
Set $M:=\spt\smallnorm{V}\setminus\Sigma$. By assumption, $M$ is a properly immersed, two-sided, stable minimal hypersurface in $(\R\times B^n_1(0))\setminus\Sigma$. In particular, $M$ satisfies the stability inequality in $(\R\times B^n_1(0))\setminus\Sigma$. Since $\textnormal{cap}_2(\Sigma)=0$, the stability inequality for $M$ can be extended across $\Sigma$, and we thus have
\begin{equation}\label{eq: stability ineq}
\int|A_M|^2\phi^2\leq\int|\nabla^M\phi|^2,\quad\text{ for all }\phi\in C^1_c(\R\times B^n_1(0)).
\end{equation}
Let $\pi:\R\times B^n_1(0)\to B^n_1(0)$ be the projection map. Set $L:=\Lip(u)$. Since $\pi$ is Lipschitz, we have that $\H^n(\pi(\Sigma\cap\spt\smallnorm{V}))=0$, and it follows from \eqref{eq: stability ineq}, the area formula and a standard cutoff argument that for every compact set $K\subset B^n_1(0)$, we have
\begin{equation}\label{eq: SMH W^2,2 estimate}
\int|D^2u|^2\phi^2\leq C(n,L,M,K)\int(|\phi|^2+|D\phi|^2),\quad\text{ for all }\phi\in C^1_c(K).
\end{equation}
Thus, there exists a function $H\in L^2_{\textnormal{loc}}(B^n_1(0);\A_Q(\R^n\times\R^{n\times n}))$ such that $Dj_1u=H$ in $B^n_1(0)\setminus\pi(\Sigma)$. Since $\H^{n-1}(\pi(\Sigma))=0$, it now follows from \cite[Theorem 2.4]{DLS11Qvalued} that $u\in W^{2,2}_{\textnormal{loc}}(B^n_1(0);\A_Q(\R))$. The result now follows from Theorem \ref{thm: Lip+W^2,2 Stationary + Lip small implies C^1,alpha}.
\end{proof}

Finally, we prove Theorem \ref{thm: C^1+W^2,2 implies C^1,alpha}:

\begin{proof}[Proof of Theorem \ref{thm: C^1+W^2,2 implies C^1,alpha}]
Let $g:\R^{n\times k}\to\R^{n\times n}$ be defined as in \eqref{eq: defn of g}. Set $F:=Du$ and write $F=(F_1,\dots,F_n)$ where each $F_i\in W^{1,2}_{\textnormal{loc}}(B^n_1(0);\A_Q(\R^k))$. Arguing as in the proof of Theorem \ref{thm: Lip+W^2,2 Stationary + Lip small implies C^1,alpha}, we see that $DF^\ell$ is symmetric and \eqref{eq: Q-valued MSS PDE} holds. Fix $x_0\in B^n_1(0)$, and write $F(x_0)=\sum_{j=1}^Jq_j\mb{P_j}$, for some distinct $P_1,\dots,P_J\in\R^{n\times k}$. Since $F$ is continuous, there exists an open neighbourhood $U\Subset B^n_1(0)$ of $x_0$ in which $F$ decomposes as the sum $F=\sum_{j=1}^J\widetilde{F}_j$, where each $\widetilde{F}_j\in W^{1,2}(U;\A_{q_j}(\R^{n\times k}))$ and satisfies $\widetilde{F}_j(x_0)=q_j\mb{P_j}$. Fix $j\in\{1,\dots J\}$ and set $q:=q_j$ and $P_0=P_j$. Then $S:=g(P_0)^{1/2}$ is a symmetric positive-definite matrix. Set $V:=S^{-1}(U-x_0)$ and define $\widehat{F}:V\to\A_q(\R^{n\times k})$ by 
$$
\widehat{F}(y):=\sum_{\ell=1}^{q}\mb{\widetilde{F}_j^\ell(x_0+Sy)S}.
$$
Define $A:\R^{n\times k}\to\R^{n\times n}$ by $A(X):=S^{-1}g(XS^{-1})S^{-1}$ and set $X_0:=P_0S$ so that $A(X_0)=I_n$. Then $D\widehat{F}^{\ell,\kappa}=SDF^{\ell,\kappa}S$  for all $\ell,\kappa$, so $D\widehat{F}$ is symmetric. Also, \eqref{eq: Q-valued MSS PDE} implies that
\begin{equation}\label{eq: transformed PDE}
A(\widehat{F}^\ell):D\widehat{F}^{\ell,\kappa}=0\quad\H^n\text{-a.e. in }V.
\end{equation}
Define $\widehat{R}:\R^{n\times k}\times\R^{n^2\times k}\to\R$ by setting
$$
\widehat{R}(X,P):=\sum_{\kappa=1}^k(A(X):P^\kappa)^2+\sum_{\kappa=1}^k\sum_{i<j}(P^\kappa_{ij}-P^\kappa_{ji})^2,
$$
By \eqref{eq: transformed PDE}, we have that $\widehat{R}(\widehat{F}^\ell,D\widehat{F}^\ell)=0$ $\H^n$-a.e. in $V$. Setting $\widehat{J}(X,P):=\widehat{R}(X,P)-N(P)$, it follows from \eqref{eq: martix symmetry property} that
\begin{equation}\label{eq: equation for transformed energy}
\widehat{J}(X,P)=|P|^2+\sum_{\kappa=1}^k[(A(X):P^\kappa)^2-(\textnormal{tr}P^\kappa)^2].
\end{equation}
Fix an open ball $B^n_R(y)\Subset V$. For every $r\in(0,R)$ such that $\widehat{F}|_{\partial B^n_r(y)}\in W^{1,2}(\partial B^n_r(y);\A_q(\R^{n\times k}))$, define the functional
$$
\widehat{\mathscr{J}}(G;B^n_r(y)):=\int_{B^n_r(y)}\sum_{\ell=1}^q \widehat{J}(G^\ell,DG^\ell),\quad\text{ for }G\in W^{1,2}(B^n_r(y);\A_q(\R^{n\times k})).
$$
If $G=\widehat{F}$ on $\partial B^n_r(y)$, then the above gives
\begin{equation*}
\widehat{\mathscr{J}}(G;B^n_r(y))-\widehat{\mathscr{J}}(\widehat{F};B^n_r(y))=\int_{B^n_r(y)}\sum_{\ell=1}^q \widehat{R}(G^\ell,DG^\ell)\geq 0.
\end{equation*}
Thus $\widehat{F}$ is a local minimiser of $\widehat{\mathscr{J}}$. Let $\eps\in(0,1)$, to be determined depending on $n$ and $Q$. Since $A(X_0)=I_n$, there exists an open ball $B=B^{n\times k}_\rho(X_0)\subset\R^{n\times k}$ such that
\begin{equation}\label{eq: choice of delta}
\delta(2\sqrt n+\delta)\leq\eps,\quad\text{ where }\quad \delta:=\sup_{X\in B}|A(X)-I_n|.
\end{equation}
For every $X\in B$ and every $P^\kappa\in\R^{n\times n}$, we have
\begin{equation*}
|(A(X):P^\kappa)^2-(\textnormal{tr}P^\kappa)^2|=|(A(X)-I_n):P^\kappa||(A(X)+I_n):P^\kappa|\leq\delta(2\sqrt n+\delta)|P^\kappa|^2.
\end{equation*}
Thus \eqref{eq: equation for transformed energy} and \eqref{eq: choice of delta} give
\begin{equation}\label{eq: closeness of integrand to square norm}
(1-\eps)|P|^2\leq \widehat{J}(X,P)\leq(1+\eps)|P|^2,\quad\text{ for all }X\in B.
\end{equation}
Since $\widehat{F}(0)=q\mb{X_0}$ and $\widehat{F}$ is continuous, after replacing $U$ with a smaller open neighbourhood, we can assume without loss of generality that every value of $\widehat{F}$ lies in $B$. Thus, \eqref{eq: closeness of integrand to square norm} implies that for any $G\in W^{1,2}(B^n_r(y);\A_q(B))$ satisfying $G=\widehat{F}$ on $\partial B^n_r(y)$, we have
\begin{equation*}
\int_{B^n_r(y)}|D\widehat{F}|^2\leq\Lambda\int_{B^n_r(y)}|DG|^2,\quad\text{ where }\Lambda:=(1+\eps)/(1-\eps).
\end{equation*}
We can now argue as in the final part of the proof of Theorem \ref{thm: Lip+W^2,2 Stationary + Lip small implies C^1,alpha} to see that $\widehat{F}\in C^{0,\alpha}(V)$, for some $\al=\al(n,Q)\in(0,1)$. In particular, each $\widetilde{F}_j$ is $C^{0,\alpha}$ and hence, $F=Du$ is $C^{0,\alpha}$. The result now follows from Proposition \ref{prop: W^2,2 + Dw cont implies C^1}.
\end{proof}

\subsection{Proof of Theorem \ref{thm: Compactness Branched SMH n=2}}
In this section, we prove Theorem \ref{thm: Compactness Branched SMH n=2}. Define $g:\R^2\to\R^{2\times 2}$ by setting $g(y):=(I_2+y\otimes y)^{-1}$. The following is the main result of the section:
\begin{thm}\label{thm: graph of derivative is calibrated in 2d}
There exists a metric $h$ on $\R^2\times B^2_1(0)$ and a calibration $\omega$ on $(\R^2\times B^2_1(0),h)$ such that the following holds. Suppose $G\in W^{1,2}_{\textnormal{loc}}(B^2_1(0);\A_Q(\R^2))$ is continuous such that for every $\ell\in\{1,\dots,Q\}$ and $i,j\in\{1,2\}$ we have
\begin{equation}\label{eq: DG symmetric + PDE}
D_iG_j^\ell(x)=D_jG_i^\ell(x)\quad\textnormal{ and }\quad g(G^\ell(x)):DG^\ell(x)=0,
\end{equation}
for $\H^2$-a.e. $x\in B^2_1(0)$. Then, $\bf{T}_G$ is calibrated in $(\R^2\times B^2_1(0),h)$ by $\omega$ and $\B_G$ is discrete.
\end{thm}
\begin{remark}\label{remark: Obstruction to Calibrations in Higher Dimensions}
In particular, $h$ and $\omega$ are independent of $G$. This property does not generalise to higher dimensions, even in the $Q=1$ case. We discuss this further in Appendix \ref{app: Obstruction to Calibrations in Higher Dimensions}. We also note that the proof of Theorem \ref{thm: graph of derivative is calibrated in 2d} below shows that when $n=2$ and $k=1$, the graph of any $Q$-valued $W^{1,2}_{\textnormal{loc}}$ function satisfying \eqref{eq: Hodge structure} is calibrated in $\R^2\times B^2_1(0)$ with the standard Euclidean metric.
\end{remark}
\begin{proof}[Proof of Theorem \ref{thm: graph of derivative is calibrated in 2d}]
Define $H:\R^2\to\R^{2\times 2}$ by setting
$$
H(y):=\frac{I_2+y\otimes y}{\sqrt{1+|y|^2}}.
$$
Then $H$ satisfies $\det H(y)=1$ and $H(y)^{-1}=\sqrt{1+|y|^2}g(y)$ for all $y\in\R^2$. For $(y,x)\in \R^2\times B^2_1(0)$, define $h_{(y,x)}:\R^{2\times 2}\times\R^{2\times 2}\to\R$ by setting
$$
h_{(y,x)}((v,w),(\eta,\xi)):=v^TH(y)^{-1}\eta+w^TH(y)\xi.
$$
Then $h:=(h_{(y,x)})_{(y,x)\in\R^2\times B^2_1(0)}$ defines a smooth metric on $\R^2\times B^2_1(0)$. Let $\omega$ be the closed $2$-form defined by
$$
\omega:=-dy_1\wedge dy_2+dx_1\wedge dx_2.
$$
We claim that $\omega$ defines a calibration on $(\R^2\times B^2_1(0),h)$. To see this, let $R:\R^2\to\R^2$ be the rotation map by $\pi/2$, and for $(y,x)\in\R^2\times B^2_1(0)$, define $J_{(y,x)}:\R^{2\times 2}\to\R^{2\times 2}$ by setting
$$
J_{(y,x)}(v,w):=(-H(y)Rv,H(y)^{-1}Rw).
$$
Since $H$ is symmetric and has unit determinant, we have $HRH=R$. Hence, $J_{(y,x)}$ satisfies
\begin{equation}\label{eq: almost complex J}
J^2_{(y,x)}=-I_4, \quad h(J_{(y,x)}(v,w),J_{(y,x)}(\eta,\xi))=h((v,w),(\eta,\xi)),
\end{equation}
and
\begin{equation}\label{eq: compatibility of omega}
\langle\omega,(v,w)\wedge(\eta,\xi)\rangle=h(J_{(y,x)}(v,w),(\eta,\xi)),
\end{equation}
for all $(v,w),(\eta,\xi)\in\R^2\times\R^2$. In particular, $\omega$ defines a calibration on $\R^2\times B^2_1(0)$. We wish to show that $\omega$ calibrates $\bf{T}_G$. To prove the claim, we first note the following. Suppose now that $P\in\R^{2\times 2}$ such that $P^T=P$ and $g(y):P=0$ for some $y\in\R^2$. Then we have $H(y)^{-1}:P=0$. In particular, the matrix $S:=H(y)^{-1/2}PH(y)^{-1/2}$ is a symmetric and trace free $2\times 2$ matrix. Thus, $S^2=-\det(S)I_2$. Since $\det(H(y))=1$, we have $\det(S)=\det(P)$, and thus
$$
P^TH(y)^{-1}P=-\det(P)H(y).
$$
Setting $\zeta(P):=(Pe_1,e_1)\wedge(Pe_2,e_2)$, it follows from the above that
\begin{equation}\label{eq: DG symmetric + PDE 1}
|\zeta(P)|_h=1-\det(P)=\langle \omega,\zeta(P)\rangle,
\end{equation}
for every such $P$. We can now show that $\omega$ calibrates $\bf{T}_G$. Indeed, let $\mu$ be the measure on $\R^2\times B^2_1(0)$ defined by 
$$
\mu(A):=\int\sum_{\ell=1}^Q\mathbf{1}_A(G^\ell(x),x)|\zeta(DG^\ell(x))|_h\,\d x.
$$
Since $G$ is locally bounded, and $DG\in L^2_{\textnormal{loc}}(B^2_1(0);\A_Q(\R^2))$, $\mu$ is locally finite. Note that $\smallnorm{\bf{T}_G}_h\leq\mu$. Furthermore, by \eqref{eq: DG symmetric + PDE} and \eqref{eq: DG symmetric + PDE 1}, for any $\psi\in C^1_c(\R^2\times B^2_1(0))$, we have
$$
\bf{T}_G(\psi\omega)=\int \psi\,\d\mu,
$$
and thus $\mu\leq\smallnorm{\bf{T}_G}_h$. In particular, we have $\smallnorm{\bf{T}_G}_h=\mu$. Furthermore, if $\tau_h$ denotes the orientation form of $\bf{T}_G$, then by the above, we have
$$
\int\psi(1-\langle\omega,\tau_h\rangle)\,\d\smallnorm{\bf{T}_G}_h=0,\quad\text{ for all }\psi\in C^1_c(\R^2\times B^2_1(0)).
$$
Hence, $\langle \omega,\tau_h\rangle=1$ $\smallnorm{\bf{T}_G}_h$-a.e. on $\spt\smallnorm{\bf{T}_G}_h$, and thus $\omega$ calibrates $\bf{T}_G$.

It remains to show that $\B_G$ is discrete. For this, we will need the following claim:

\begin{claim}\label{claim: local parametrisation of graph of derivative}
Fix $X\in\spt\smallnorm{\bf{T}_G}$. There exists $\rho>0$, an integer $N\geq 1$, such that 
\begin{equation}\label{eq: local parametrisation of graph of derivative}
\bf{T}_G\res B^4_\rho(X)=\sum_{j=1}^Nm_j\cdot(\Psi_j)_{\#}\mb{D_j}\res B^4_\rho(X),
\end{equation}
where, for each $j\in\{1,\dots,N\}$, $m_j\geq 1$ is an integer, $D_j\subset \R^2$ is an open disk, and $\Psi_j=\Psi_j(s,t):D_j\to\R^2\times B^2_1(0)$ is a smooth proper map which is a homeomorphism onto its image and satisfies $\Psi_j(0)=X$ and
\begin{equation}\label{eq: PDE for parametrisations}
(\Psi_j)_t=J_{\Psi_j}(\Psi_j)_s\quad \text{ in }D_j.
\end{equation}
\end{claim}
\begin{proof}[Proof of Claim \ref{claim: local parametrisation of graph of derivative}]
By the Nash embedding theorem, there exists an isometric embedding $\iota:(\R^2\times B^2_1(0),h)\to\R^d$, for some integer $d\geq 1$. Since $\iota$ is an isometry, $\iota_{\#}\bf{T}_G$ is area minimising within the $4$-dimensional submanifold $\iota(\R^2\times B^2_1(0))$ of $\R^d$. After rescaling in $\R^d$, we can ensure that $\iota$ is chosen so that there exists $\rho_0>0$ such that $\Sigma:=\iota(B^4_{\rho_0}(X))$ satisfies \cite[Assumption 1.1]{DLSS20}. Furthermore, arguing as in \cite[p.~128]{DLSS20}, we obtain $\rho_1>0$, an integer $N\geq 1$ and integer rectifiable currents $T_1,\dots,T_N$ in $B^d_{\rho_1}(\iota(X))$ with $\spt\smallnorm{T_j}\subset B^d_{\rho_1}(\iota(X))\cap\Sigma$ such that each $T_j$ satisfies \cite[Assumption 1.2]{DLSS20} and
$$
\iota_{\#}\bf{T}_G\res B^d_{\rho_1}(\iota(X))=\sum_{j=1}^NT_j\res B^d_{\rho_1}(\iota(X)).
$$
It now follows from the argument of \cite[Section 1.5]{DLSS20} that there exists $\rho_2>0$ such that for each $j\in\{1,\dots,N\}$, there exists integers $m_j,q_j\geq 1$ such that $T_j\res B^d_{\rho_2}(\iota(X))$ is an $m_j$-multiple of a $q_j$-branching in the sense of \cite[Definition 1.4]{DLSS20}. Furthermore, by \cite[Proposition 2.4]{DLSS20}, we can ensure that each such $q_j$-branching is represented by a conformal map. Pre-composing each such conformal $q_j$-branching with the map $\zeta\mapsto(\zeta^{q_j},\zeta)$ and pulling back by $\iota$, we obtain $r,\rho>0$ and an integer $N\geq 1$, such that \eqref{eq: local parametrisation of graph of derivative} holds, where for each $j\in\{1,\dots,N\}$,  $D_j\subset\R^2$ is an open disc containing the origin and $\Psi_j:D_j\to B^4_\rho(X)$ satisfies:
\begin{itemize}
\item[(i)] $\Psi_j$ is a homeomorphism onto its image;
\item[(ii)] $\Psi_j$ is $C^3$ in $D_j\setminus\{0\}$, and we have 
\begin{equation}\label{eq: gradient bound for conformal parametrisation}
|D\Psi_j(\zeta)|\leq C_j|\zeta|^{q_j-1},\quad\text{ for all }\zeta\in D_j\setminus\{0\},
\end{equation}
where $C_j\in (0,\infty)$ is independent of $\zeta$;
\item[(iii)] $\Psi_j$ is conformal on $D_j\setminus\{0\}$.
\end{itemize}
Fix $j\in\{1,\dots,N\}$ and set $D:=D_j, \Psi:=\Psi_j$. Then by the above, $\Psi|_{D\setminus\{0\}}$ is locally a minimal conformal parametrisation in $(\R^2\times B^2_1(0),h)$. In particular, by \cite[p.~116, 126]{EellsSampson64Harmonic}, $\Psi$ satisfies
\begin{equation}\label{eq: PDE for conformal parametrisation}
\Delta\Psi^\nu=-\Gamma^\nu_{\kappa\gamma}(\Psi)(\Psi^\kappa_s\Psi^\gamma_s+\Psi^\kappa_t\Psi^\gamma_t)=:\Phi^\nu,\quad\text{ in }D\setminus\{0\},
\end{equation}
for every $\nu\in\{1,\dots,4\}$, where $\Gamma^\nu_{\kappa\gamma}$ are the Christoffel symbols of $h$. By \eqref{eq: gradient bound for conformal parametrisation}, we have that $D\Psi,\Phi^\nu\in L^\infty(D)$. Since the singleton $\{0\}$ has zero $1$-capacity in $D$, it follows from a standard capacity argument that \eqref{eq: PDE for conformal parametrisation} holds weakly in $D$. Hence, a standard difference quotient argument gives that $\Psi\in W^{2,p}_{\textnormal{loc}}(D)$ for all $p\in(1,\infty)$. In particular, the Sobolev embedding theorem implies that $\Psi\in C^{1,\beta}(D)$, for all $\beta\in(0,1)$. It now follows from standard Schauder theory that $\Psi$ is smooth in $D$. Finally, since $\Psi$ is conformal in $D\setminus\{0\}$, we have
$$
h(\Psi_s,\Psi_s)=h(\Psi_t,\Psi_t)>0\quad\text{ and }\quad h(\Psi_s,\Psi_t)=0,\quad\text{ in }D\setminus\{0\}.
$$
Furthermore, since $\omega$ calibrates $\bf{T}_G$, we have
$$
\langle\omega,\Psi_s\wedge\Psi_t\rangle=|\Psi_s\wedge\Psi_t|_h=h(\Psi_s,\Psi_s),\quad\text{ in }D\setminus\{0\}.
$$
By \eqref{eq: compatibility of omega}, this implies that $\Psi_t=J_\Psi\Psi_s$ in $D\setminus\{0\}$. Since $\Psi$ is smooth, this establishes \eqref{eq: PDE for parametrisations}, and completes the proof of the claim.
\end{proof}
Fix $X\in\spt\smallnorm{\bf{T}_G}$ and let $\Psi:D\to\R^2\times B^2_1(0)$ be a smooth map given by Claim \ref{claim: local parametrisation of graph of derivative}. Define $F:=\pi\circ\Psi$ and $Y:=\pi^\perp\circ\Psi$, where $\pi:\R^2\times B^2_1(0)\to\{0\}^2\times B^2_1(0)$ is the projection map and $\pi^\perp:=I_4-\pi$. We claim that the critical points of $F$ are isolated. By \eqref{eq: PDE for parametrisations}, $F$ satisfies 
\begin{equation}\label{eq: PDE for F}
F_t=H(Y)^{-1}RF_s,\quad\text{ in }D.
\end{equation}
Let $a,b,c:D\to\R$ be smooth functions such that $a(s,t)>0, a(s,t)\cdot c(s,t)-b(s,t)^2=1$ and 
$$
H(Y(s,t))^{-1}=
\begin{pmatrix}
a(s,t)&b(s,t)\\
b(s,t)&c(s,t)
\end{pmatrix},\quad\text{ for all }(s,t)\in D.
$$
Write $F=(u,v)$ for some smooth functions $u,v:D\to\R$. By \eqref{eq: PDE for F}, we have that
$$
u_t=bu_s-av_s\quad\text{ and }\quad v_t=cu_s-bv_s,\quad\text{ in }D.
$$
Solving for $Dv$ gives
\begin{equation}\label{eq: u & v coupling}
v_s=a^{-1}(bu_s-u_t)\quad\text{ and }\quad v_t=a^{-1}(u_s+bu_t),\quad\text{ in }D.
\end{equation}
In particular, $u$ cannot be constant. Indeed, if $u$ were constant, then \eqref{eq: u & v coupling} would imply that $F$ is constant. Since $\spt\smallnorm{\bf{T}_G}=\graph(G)$, constancy of $F$ would imply that $\Psi(D)$ is finite, which, since $\Psi$ is proper, is a contradiction. Now, by \eqref{eq: u & v coupling}, we have that
\begin{equation}\label{eq: critical points of F}
\det(DF)=a^{-1}|Du|^2,\quad\text{ in }D,
\end{equation}
and thus the critical points of $F$ are the critical points of $u$. Furthermore, since $D^2v$ is symmetric, it follows from \eqref{eq: u & v coupling} that
$$
D_s(a^{-1}(u_s+bu_t))+D_t(a^{-1}(u_t-bu_s))=0,\quad\text{ in }D.
$$
Expanding and multiplying through by $a$, we see that $u$ solves
$$
\Delta u+B\cdot\nabla u=0,\quad\text{ in }D,
$$
for some smooth function $B:D\to\R^2$. Since $u$ is non-constant, we can apply the Hartman-Wintner theorem (cf.~\cite[Corollary 1]{HartmanWintner1953}) to see that the critical points of $u$ are isolated. By \eqref{eq: critical points of F}, this implies that the critical points of $F$ are isolated, which proves the claim. Finally, if $\textnormal{rank}(D\Psi(s,t))<2$ for some $(s,t)\in D$, then \eqref{eq: almost complex J} and \eqref{eq: PDE for parametrisations} imply that $D\Psi(s,t)=0$, and thus $(s,t)$ is a critical point of $F$. Let $\Sigma_{\Psi}$ denote the critical points of $\pi|_{\reg(\Psi_{\#}\mb{D}\res B^4_\rho(X))}$. Since $\Psi$ is a proper map, it now follows that $\sing(\Psi_{\#}\mb{D}\res B^4_\rho(X))\cup\Sigma_\Psi$ is discrete. 

Let $\Sigma_G$ denote the set of critical points of $\pi|_{\reg(\bf{T}_G)}$. By the above, the set $\sing(\bf{T}_G)\cup\Sigma_G$ is discrete. Since $\pi|_{\spt\smallnorm{\bf{T}_G}}$ is proper, it follows that the set $\pi(\Sigma_G\cup\sing(\bf{T}_G))$ is also discrete. Finally, the inverse function theorem implies that $\B_G\subseteq\pi(\Sigma_G\cup\sing(\bf{T}_G))$. In particular, $\B_G$ is discrete, which completes the proof.
\end{proof}

We will also need the following:
\begin{prop}\label{prop: branch set of derivative of stationary graph}
Fix $Q, n\geq 2, k\geq 1$ and $\al\in(0,1)$. Suppose that $u\in C^{1,\al}(B^n_1(0);\A_Q(\R^k))$ such that the integral varifold $\bf{v}(u)$ is stationary in $\R^k\times B^n_1(0)$. Then we have $\B_u=\B_{Du}$.
\end{prop}

In order to prove Proposition \ref{prop: branch set of derivative of stationary graph}, we will need a further technical lemma:
\begin{lem}\label{lem: C^1 implies analytic for selections of deriviative of stat graph}
Fix $n,Q\geq 2$ and $k\geq 1$. Suppose that $u\in C^1(B^n_1(0);\A_Q(\R^k))$ such that the integral varifold $\bf{v}(u)$ is stationary in $\R^k\times B^n_1(0)$. Suppose that $B\subset B^n_1(0)$ is an open ball and $H_1,\dots, H_q:B\to\R^{n\times k}$ are single-valued $C^1$ functions such that $Du\equiv\sum_{i=1}^q\mb{H_i}$ in $B$. Then each $DH_i$ is symmetric and satisfies $g(H_i):DH_i=0$ in $B$, where $g(X):=(I_n+X^TX)^{-1}$. In particular, each $H_i$ is analytic.
\end{lem}
\begin{proof}
We argue by induction on $Q$. The result follows by standard elliptic regularity when $Q=1$ (cf.~\cite{Mor66MultipleIntegrals}). Let $q\geq 2$ and suppose that the result has been established whenever $Q\leq q-1$. Suppose that $u\in C^1(B^n_1(0);\A_q(\R^k))$ such that the integral varifold $\bf{v}(u)$ is stationary in $\R^k\times B^n_1(0)$. Suppose that $B\subset B^n_1(0)$ is an open ball and $H_1,\dots, H_q:B\to\R^{n\times k}$ are single-valued $C^1$ functions such that $Du\equiv\sum_{i=1}^q\mb{H_i}$ in $B$. By the induction hypothesis and Proposition \ref{prop: separation principle}, each $DH_i$ is symmetric in $B\setminus\CK_u$, and satisfies $g(H_i):DH_i=0$ in $B\setminus\CK_u$. Furthermore, each $H_i$ is continuous, and satisfies $H_i=Du_a$ and $DH_i=D^2u_a$ $\H^n$-a.e. on $\CK_u$. Thus, each $DH_i$ is symmetric and satisfies $g(H_i):DH_i=0$ in $B$. In particular, since $B$ is simply connected, it follows from the Poincaré lemma that for each $i\in\{1,\dots,q\}$, there exists a function $u_i\in C^2(B)$ such that $Du_i\equiv H_i$ in $B$. Furthermore, each $u_i$ is a $C^2$ single-valued solution of the minimal surface system, and hence is analytic (cf.~\cite{Mor66MultipleIntegrals}). This completes the proof.
\end{proof}
We can now prove Proposition \ref{prop: branch set of derivative of stationary graph}:
\begin{proof}[Proof of Proposition \ref{prop: branch set of derivative of stationary graph}]
Fix $\al\in(0,1)$. We argue by induction on $Q$. The result follows from \cite{Mor66MultipleIntegrals} when $Q=1$. Let $q\geq 2$ and suppose that the result has been established whenever $Q\leq q-1$. Suppose $u\in C^{1,\al}(B^n_1(0);\A_q(\R^k))$ and that the integral varifold $\bf{v}(u)$ is stationary in $\R^k\times B^n_1(0)$. We can now use Lemma \ref{lem: C^1 implies analytic for selections of deriviative of stat graph} to argue analogously to the proof of Lemma \ref{lem: branch set contained in that of derivative assuming ind hyp}, using \cite[Lemma 5.3]{LucheseSchauder} in place of Lemma \ref{lem: separation harmonic assuming induction hyp}, and Theorem \ref{thm: weak unique continuation for stationary C^1,alpha graphs} in place of \cite[Lemma 6.9]{DLMSVDirRectifiability}. This completes the proof. 
\end{proof}

Finally, Theorem \ref{thm: Compactness Branched SMH n=2} now follows from Theorem \ref{thm: Bel25 Theorem 5 reg improvement}, \cite[Theorem 5]{Bel25}, \cite[Theorems 1.3 \& 1.9]{WZ26StableMinHypNonBranchedSingSet}, Proposition \ref{prop: proof of Bel25 Theorem 5 reg improvement} (cf.~\eqref{eq: SMH W^2,2 estimate}), \cite[Lemma 2.10]{LucheseSchauder}, Arzelà-Ascoli for multi-valued Lipschitz functions and the following:
\begin{prop}
Fix $Q\geq 1$. There exists $\eps=\eps(Q)\in(0,1)$ such that the following holds. Suppose that $u\in W^{2,2}_{\textnormal{loc}}(B^2_1(0);\A_Q(\R))$ is Lipschitz such that the integral varifold $\bf{v}(u)$ is stationary in $\R\times B^2_1(0)$. If $\Lip(u)\leq\eps$, then $\B_u$ is discrete.
\end{prop}
\begin{proof}
By Theorem \ref{thm: Lip+W^2,2 Stationary + Lip small implies C^1,alpha}, we have that $u\in C^{1,1/Q}(B^2_1(0);\A_Q(\R))\cap W^{2,2}_{\textnormal{loc}}(B^2_1(0);\A_Q(\R))$. Setting $G:=Du$, it follows from \cite[Lemma 3.10]{LucheseSchauder} that $DG^\ell$ is symmetric $\H^2$-a.e. in $B^2_1(0)$ for all $\ell$. Furthermore, arguing as in the proof of \cite[Lemma 3.10]{LucheseSchauder}, we see that $g(G^\ell):DG^\ell=0$ $\H^2$-a.e. in $B^2_1(0)$ for all $\ell$. The result now follows from Theorem \ref{thm: graph of derivative is calibrated in 2d} and Proposition \ref{prop: branch set of derivative of stationary graph}.
\end{proof}

\section{Applications to the Thin Obstacle Problem}
It is well known that the regularity theory of two-valued functions with variational properties has applications to the thin obstacle problem \cite{AC04ThinObstacle,FS18ThinObstacle,FS20ThinObstacle,FR22TOP,FS24ThinObstacle}. The aim of this section is to elaborate on this connection by discussing the applications of the present methods to the thin obstacle problem. 
\subsection{The Thin Obstacle Problem for the Laplacian}
Fix $n\geq 2$, and for $r>0$ denote $B^\pm_r(0):=B^n_r(0)\cap\{\pm x_n\geq0\}$. Let $H$ denote the set of all $v\in W^{1,2}(B^+_1(0))$ with $v\geq 0$ on $\{x_n=0\}$. We say that $u\in H$ solves the \textit{thin obstacle problem for the Laplacian} if there holds
\begin{equation}\label{eq: thin obstacle problem for Laplacian}
\int|Du|^2\leq\int|Dv|^2,\quad\text{ for all }v\in H\text{ such that }u=v\text{ on }\partial B^n_1(0)\cap\{x_n>0\}.
\end{equation}
The thin obstacle for the Laplacian is often referred to as the Signorini problem. If $u$ solves the thin obstacle problem for the Laplacian, the free boundary, denoted $\Gamma(u)$, is the relative boundary of the set $\{u=0\}$ in $\{x_n=0\}$. 

The following is our main result for the thin obstacle for the Laplacian:
\begin{thm}\label{thm: Derivative of TOP for Lap is locally Dir min}
Suppose that $u\in W^{1,2}(B^+_1(0))$ is a solution to the thin obstacle problem for the Laplacian. Let $w$ denote the even extension of $u$ in $B^n_1(0)$ across the set $\{x_n=0\}$. Then, we have that $G:=\mb{Dw}+\mb{-Dw}\in W^{1,2}_{\textnormal{loc}}(B^n_1(0);\A_2(\R^n))$. Furthermore, $G$ is locally Dirichlet minimising in $B^n_1(0)$ and we have $\Gamma(u)\subseteq\{G=2\mb{0}\}$.
\end{thm}
Theorem \ref{thm: Derivative of TOP for Lap is locally Dir min} implies that the derivative of the two-valued $W^{2,2}_{\textnormal{loc}}$ function $v:=\mb{w}+\mb{-w}$ is locally Dirichlet minimising. Thus, \cite[Theorem 2.13]{Almgren00} and Proposition \ref{prop: W^2,2 + Dw cont implies C^1} recover the known fact that $v$ is a two-valued $C^{1,1/2}(B^n_1(0);\A_2(\R))$  harmonic function. In particular, $u\in C_{\textnormal{loc}}^{1,1/2}(\overline{B^+_1(0)})$, which recovers the result of \cite{AC04ThinObstacle}. Furthermore, since $Dv$ is average-free, provided that $Dv\not\equiv2\mb{0}$\footnote{If $Dv\equiv 2\mb{0}$, then \cite[Proposition 2.12]{DLS11Qvalued} implies that $u\equiv\textnormal{const}$, and thus the free boundary is empty.}, \cite{KW17} and \cite{DLMSVDirRectifiability} imply countable $(n-2)$-rectifiability, together with local finiteness of $\H^{n-2}$ measure and Minkowski-type estimates for the set $\{Dv=2\mb{0}\}$. In particular, since $\Gamma(u)\subseteq \{Dv=2\mb{0}\}$, Theorem \ref{thm: Derivative of TOP for Lap is locally Dir min} recovers the results of \cite[Theorem 1.1 \& 1.2]{FS18ThinObstacle} for the Signorini problem (i.e. when $a=0$). We note that, by \cite[Corollary 1.2]{FS24ThinObstacle}, the countable $(n-2)$-rectifiability conclusion for $\Gamma(u)$ is not optimal, and in fact $\Gamma(u)$ is contained in the countable union of $(n-2)$-dimensional submanifolds of class $C^{1,1}$ and a set of Hausdorff dimension at most $n-3$. Finally, when $n=2$, Corollary \ref{cor: Full Series Expansion in Disc} recovers the known series expansion for $u$ near $\Gamma(u)$ (cf.~\cite{Lewy72CoincidenceSet}, \cite[Theorem 1.1~(b)]{DPSV25}).

We now prove Theorem \ref{thm: Derivative of TOP for Lap is locally Dir min}:
\begin{proof}[Proof of Theorem \ref{thm: Derivative of TOP for Lap is locally Dir min}]
For completeness, we first show the known fact that $u\in W^{2,2}_{\textnormal{loc}}(B^+_1(0))$. To see this, we note that by \eqref{eq: thin obstacle problem for Laplacian}, for every $\phi\in W^{1,2}(B^+_1(0))$ satisfying
\begin{equation}\label{eq: admissable competitors}
\phi=0\text{ on }\partial B^n_1(0)\cap\{x_n>0\}\quad\text{ and }\quad u+t\phi\in H\quad\text{ for all }t\in[0,1]\text{ sufficiently small},
\end{equation}
we have
\begin{equation}\label{eq: sub-harmonic}
\int Du\cdot D\phi\geq 0.
\end{equation}
Fix $i\in\{1,\dots,n-1\},\theta\in(0,1)$ and let $\eta\in C^\infty_c(B^n_1(0))$ be such that $0\leq\eta\leq 1$, $\eta=1$ in a neighbourhood of $B^n_\theta(0)$ and $|D\eta|\leq C(n,\theta)$. For $h\in \R$ such that $0<|h|\leq \dist(\spt(\eta),\partial B^n_1(0))$ and $v\in L^2(B^+_1(0))$, define $v_h,\delta_hv\in L^{2}(B^+_{\theta}(0))$ by setting $v_h(x):=v(x+he_i)$ and $\delta_hv:=v_h-v$. Set $\phi_1:=\eta^2\cdot\delta_hu$. Then $\phi_1=0$ on $\partial B^n_1(0)\cap\{x_n>0\}$, and for every $t\in[0,1]$ we have $u+t\phi_1=(1-t\eta^2)u+t\eta^2u_h\geq 0$ on $\{x_n=0\}$. Hence, $\phi_1$ satisfies the requirements of \eqref{eq: admissable competitors}, and thus
\begin{equation}\label{eq: higher integrability 1}
\int Du\cdot D(\eta^2\delta_hu)\geq 0.
\end{equation}
Similarly, setting $\phi_2(x):=\eta^2(x-he_i)\delta_{-h}u(x)$, for all $|h|$ sufficiently small, we have that $\phi_2=0$ on $\partial B^n_1(0)\cap\{x_n>0\}$, and for every $t\in[0,1]$ we have $u+t\phi_2=(1-t\eta^2_{-h})u+t\eta^2_{-h}u_{-h}\geq 0$ on $\{x_n=0\}$. Hence, $\phi_2$ satisfies the requirements of \eqref{eq: admissable competitors}, and after a change of variables, we have
\begin{equation}\label{eq: higher integrability 2}
\int Du_h\cdot D(\eta^2\delta_hu)\leq 0.
\end{equation}
Combining \eqref{eq: higher integrability 1} with \eqref{eq: higher integrability 2} and applying Young's inequality gives
$$
\int\eta^2|\delta_hDu|^2\leq 4\int|D\eta|^2|\delta_hu|^2.
$$
Dividing by $h^2$ and applying \cite[Section 5.8.2]{Evans10PDE} gives
$$
\int_{B^+_\theta(0)}\biggl|\frac{\delta_hDu}{h}\biggr|^2\leq C(n,\theta)\int_{B^+_1(0)}|Du|^2,
$$
It now follows that $D_iDu\in L^2_{\textnormal{loc}}(B^+_1(0))$, for all $i\in\{1,\dots, n-1\}$. In particular, since $\Delta u=0$ in $B^+_1(0)\cap\{x_n>0\}$, we have that $D_nD_{n}u\in L^2_{\textnormal{loc}}(B^+_1(0))$. Finally, since $D_nD_iu=D_iD_nu$ in $B^+_1(0)\cap\{x_n>0\}$, we have that $D_nDu\in L^2_{\textnormal{loc}}(B^+_1(0))$, as required. 

Let $w$ denote the even extension of $u$ in $B^n_1(0)$ across the set $\{x_n=0\}$, and set $G:=\mb{Dw}+\mb{-Dw}$. Then by the above, $G\in W^{1,2}_{\textnormal{loc}}(B^\pm_1(0);\A_2(\R^n))$. Furthermore, the one-sided traces of $G|_{B^\pm_1(0)}$ agree on $\{x_n=0\}\cap B^n_1(0)$, and it follows from \cite[Theorem 2.4 \& Proposition 2.10]{DLS11Qvalued} that $G\in W^{1,2}_{\textnormal{loc}}(B^n_1(0);\A_2(\R^n))$. Since $\Delta w=0$ in $B^n_1(0)\setminus\{x_n=0\}$, $G$ satisfies \eqref{eq: Hodge structure}, and thus Remark \ref{remark: local Dir minimiser} implies that $G$ is locally Dirichlet minimising. Finally, by \cite[Theorem 2.13]{Almgren00} and Proposition \ref{prop: W^2,2 + Dw cont implies C^1}, we have that $v:=\mb{w}+\mb{-w}\in C^{1,1/2}_{\textnormal{loc}}(B^n_1(0);\A_2(\R))$, and thus $u\in C^1_{\textnormal{loc}}(\overline{B^+_1(0)})$. Suppose that $x=(x'_0,0)\in\Gamma(u)$. Since $u\geq 0$ on $B^n_1(0)\cap\{x_n=0\}$, the $C^1$ function $f(x'):=u(x',0)$ has a local minimum at $x'_0$. In particular, $\nabla f(x'_0)=0$, which implies that $D_iu(x)=0$ for all $i\in\{1,\dots,n-1\}$. Furthermore, by definition of $\Gamma(u)$, there exists a sequence $(x_j)_j\subset (B^n_1(0)\cap\{x_n=0\})\setminus\{u=0\}$  which converges to $x$. Since $u\in C^1_{\textnormal{loc}}(\overline{B^+_1(0)})$, \eqref{eq: sub-harmonic} and the divergence theorem implies that $D_nu(x_j)=0$ for each $j\geq 1$. This gives $Du(x)=0$, which completes the proof.
\end{proof}

\subsection{The Thin Obstacle Problem for Minimal Surfaces}
Fix $n\geq 2$. Let $S$ denote the set of all $v\in W^{1,\infty}(B^+_1(0))$ with $v\geq 0$ on $\{x_n=0\}$. We say that $u\in S$ solves the \textit{thin obstacle problem for minimal surfaces} if there holds
\begin{equation}\label{eq: Thin Obstacle Problem for MSE}
\int\sqrt{1+|Du|^2}\leq\int\sqrt{1+|Dv|^2},\quad\text{ for all }v\in S\text{ such that }u=v\text{ on }\partial B^n_1(0)\cap\{x_n>0\}.
\end{equation}
If $u$ solves the thin obstacle problem for minimal surfaces, the free boundary, denoted $\Gamma(u)$, is the relative boundary of the set $\{u=0\}$ in $\{x_n=0\}$. 

The following is our main result for the thin obstacle problem for minimal surfaces:
\begin{thm}\label{thm: Derivative of Solution to TOP for MSE}
Suppose that $u\in W^{1,\infty}(B^+_1(0))$ is a solution to the thin obstacle problem for minimal surfaces. Let $w$ denote the even extension of $u$ in $B^n_1(0)$ across the set $\{x_n=0\}$. Then, we have that $G:=\mb{Dw}+\mb{-Dw}\in W^{1,2}_{\textnormal{loc}}(B^n_1(0);\A_2(\R^n))$, and for any open ball $B\Subset B^n_1(0)$, $G$ minimises the functional
\begin{equation}\label{eq: Functional minimised by derivative of solution to TOP for MSE}
G|_{\partial B}+W^{1,2}_0(B;\A_2(\R^n))\ni F\mapsto \int_B\sum_{\ell=1}^2(|DF^\ell|^2+(g(F^\ell):DF^\ell)^2-(\textnormal{tr}DF^\ell)^2),
\end{equation}
where $g:\R^{n}\to\R^{n\times n}$ is the map $g(X):=(I_n+X^TX)^{-1}$. Furthermore, we have $\Gamma(u)\subseteq\{G=2\mb{0}\}$.
\end{thm}

Our methods also give a new proof of the previously known fact that a solution to the thin obstacle problem for minimal surfaces is $C^{1,\al}_{\textnormal{loc}}(\overline{B^+_1(0)})$ for some $\al\in(0,1)$ (cf.~\cite[Corollary 4.4]{FS20ThinObstacle}). Indeed, suppose that $u\in W^{1,\infty}(B^+_1(0))$ is a solution to the thin obstacle problem for minimal surfaces. By \cite[Proposition 3.1]{FS20ThinObstacle}, we have that $u\in C^1_{\textnormal{loc}}(\overline{B^+_1(0)})$. Thus, combining the proof of Theorem \ref{thm: Derivative of Solution to TOP for MSE} with that of Theorem \ref{thm: C^1+W^2,2 implies C^1,alpha}, we get that $v:=\mb{w}+\mb{-w}\in C^{1,\al}(B^n_1(0);\A_2(\R))$. In fact, by the argument of \cite[Section 5]{FS20ThinObstacle}, the integral varifold $\bf{v}(v)$ is stationary in $\R\times B^n_1(0)$, and it thus follows from \cite{SW16} that $v\in C^{1,1/2}(B^n_1(0);\A_2(\R))$. Furthermore, Theorem \ref{thm: Derivative of Solution to TOP for MSE} implies that $\Gamma(u)\subseteq\CK_v$, and it thus follows from \cite{SW16,KW21} that $\Gamma(u)$ is discrete when $n=2$, and countably rectifiable when $n\geq 3$, which recovers part of the conclusions of \cite[Theorem 1.1]{FS20ThinObstacle}.

We now prove Theorem \ref{thm: Derivative of Solution to TOP for MSE}:
\begin{proof}[Proof of Theorem \ref{thm: Derivative of Solution to TOP for MSE}]
The proof follows analogously to the proof of Theorem \ref{thm: Derivative of TOP for Lap is locally Dir min}. Indeed, we have that $u\in W^{2,2}_{\textnormal{loc}}(B^+_1(0))$ (cf.~\cite[Proposition 4.2]{FS20ThinObstacle}). One way to see this, is to note that by \eqref{eq: Thin Obstacle Problem for MSE}, for every $\phi\in W^{1,\infty}(B^+_1(0))$ satisfying $\phi=0$ on $\partial B^n_1(0)\cap\{x_n>0\}$ and $u+t\phi\in S$ for all $t\in[0,1]$ sufficiently small, we have
\begin{equation*}
\int \frac{Du\cdot D\phi}{\sqrt{1+|Du|^2}}\geq 0.
\end{equation*}
Furthermore, the function $A:\R^n\to\R^n$ defined by $A(p):=p/\sqrt{1+|p|^2}$ satisfies
$$
|A(p)-A(q)|^2\leq |p-q|^2\leq (1+|p|^2+|q|^2)^{3/2}(A(p)-A(q))\cdot(p-q),\quad\text{ for all }p,q\in\R^n.
$$
It now follows from a standard difference quotient argument analogous to that in the proof of Theorem \ref{thm: Derivative of TOP for Lap is locally Dir min} that $D_iDu\in L^2_{\textnormal{loc}}(B^+_1(0))$ for all $i\in\{1,\dots n-1\}$. Since $u$ solves the minimal surface equation in $B^n_1(0)\cap\{x_n>0\}$ and $Du\in L^\infty(B^+_1(0))$, this then implies that $D_nD_nu\in L^2_{\textnormal{loc}}(B^+_1(0))$. The claimed $W^{2,2}_{\textnormal{loc}}(B^+_1(0))$ regularity of $u$ now follows from the fact that $D_nD_iu=D_iD_nu$ in $B^+_1(0)\cap\{x_n>0\}$. Let $w$ denote the even extension of $u$ in $B^n_1(0)$ across the set $\{x_n=0\}$. Then arguing as in the proof of Theorem \ref{thm: Derivative of TOP for Lap is locally Dir min}, we see that $G\in W^{1,2}_{\textnormal{loc}}(B^n_1(0);\A_2(\R^n))$. Furthermore, since $G$ is the derivative of a two-valued $W^{2,2}_{\textnormal{loc}}$ function, \cite[Lemma 3.10]{LucheseSchauder} implies that each $DG^\ell$ is symmetric. Finally, since $u$ solves the minimal surface equation in $B^n_1(0)\cap\{x_n>0\}$, for each $\ell=1,2$, we have 
$$
g(G^\ell):DG^\ell=0,\quad\H^n\text{-a.e. in }B^n_1(0).
$$
The fact that $G$ locally minimises \eqref{eq: Functional minimised by derivative of solution to TOP for MSE} now follows from the proof of Theorem \ref{thm: Lip+W^2,2 Stationary + Lip small implies C^1,alpha}. Finally, by \cite[Proposition 3.1]{FS20ThinObstacle}, we have that $u\in C^1_{\textnormal{loc}}(\overline{B^+_1(0)})$. Arguing analogously to the last part of proof of Theorem \ref{thm: Derivative of TOP for Lap is locally Dir min} gives that $\Gamma(u)\subseteq\{G=2\mb{0}\}$. This completes the proof.
\end{proof}
\appendix

\section{Branch Sets of Derivatives}\label{app: containments of branch sets}
Fix $n\geq 2$ and $k\geq 1$. In this section, we record some facts about the relationship between the branch set of a $Q$-valued $C^1$ function and that of its derivative. Our first result shows that the former is not necessarily contained in the latter when $Q\geq 3$.
\begin{prop}\label{prop: branch set is not necessarily contained in that of derivative}
For every $Q\geq 3$, there exists a $Q$-valued $C^1$ function $w_Q:B^n_1(0)\to\A_Q(\R^k)$ such that $\B_{w_Q}$ is not contained in $\B_{Dw_Q}$. 
\end{prop}
\begin{remark}
The proof given below shows that we can in fact choose $w_Q$ such that $\B_{w_Q}$ is non-empty and $Dw_Q$ decomposes globally in $B^n_1(0)$ as the sum of smooth functions.
\end{remark}
\begin{proof}
It suffices to show the result when $(n,k,Q)=(2,1,3)$. Let $\{\zeta_i\}_{i=1}^3\subset\partial \D$ be the roots of $\zeta^3=1$. For $i\in\{1,2,3\}$, define $f_i:\D\to\R$ by setting $f_i(z):=\rho(z\cdot\zeta_i)$, where $\rho:\R\to\R$ is defined by $\rho(t):=\exp(-1/t^2)$ for $t>0$ and $\rho(t):=0$ for $t\leq 0$. Then each $f_i$ is smooth and hence $C^1$. Furthermore, every $z\in\D\setminus\{0\}$ is contained in the interior of a set $\{f_i=0\}$ for some $i\in\{1,2,3\}$. Let $v:\D\to\A_2(\R)$ be the function uniquely defined by $v(z)+\mb{0}=\sum_{i=1}^3\mb{f_i(z)}$. Define $w:\D\to\A_3(\R)$ by setting $w:=v+\mb{2}$. Then by the above $\B_w\subseteq\{0\}$ and $Dw=\sum_{i=1}^3\mb{Df_i}$. In particular, the set $\B_{Dw}$ is empty. Let $\{\xi_i\}_{i=1}^3\subset\partial \D$ be the roots of $\xi^3=-1$. Then 
$$
\{\textnormal{card}(v)=1\}=\{t\xi_i:t\in[0,1),i\in\{1,2,3\}\}=:\Sigma,
$$
and $\textnormal{card}(w(z))=3$ for all $z\in\D\setminus\Sigma$. Hence, $\graph(v)$ does not admit an alternating orientation (see e.g.~\cite[Definition 1.7]{MS26}) and thus, arguing as in the proof of \cite[Theorem A.2~(a)]{MS26}, we see that $0\in\B_v$. In turn, this implies that $0\in\B_w$, which completes the proof. 
\end{proof}
In contrast to Proposition \ref{prop: branch set is not necessarily contained in that of derivative}, the branch set of a two-valued $C^1$ function is always contained in that of its derivative:
\begin{prop}\label{prop: Branch set of derivatives two-valued}
If $w:B^n_1(0)\to\A_2(\R^k)$ is a two-valued $C^1$ function, then $\B_w\subseteq\B_{Dw}$. 
\end{prop}
\begin{proof}
Assume without loss of generality that $w$ is average-free. Suppose $B\subseteq B^n_1(0)$ is an open ball and $G\in C^1(B;\R^{n\times k})$ such that $Dw\equiv \mb{G}+\mb{-G}$ in $B$. Write $G=(G_1,\dots,G_n)$ for some $C^1$ functions $G_1,\dots,G_n:B\to\R^k$. Note that, by Proposition \ref{prop: separation principle}, we have that $\B_w\subseteq\CK_w$. In particular, for every $z\in\{G\neq0\}$ there exists $r>0$ and a $C^2$ function $u_z:B^n_r(z)\to\R^k$ such that 
\begin{equation}\label{eq: selection of 1-jet away from zeros of G}
j_1w\equiv\mb{(u_z,G)}+\mb{(-u_z,-G)}\quad\text{ in }B^n_r(z).
\end{equation}
Hence, $D_i G_j\equiv D_jG_i$ on $\{G\neq 0\}$, and it follows from the continuity of $DG$ that
\begin{equation}\label{eq: symmetry of G}
D_i G_j\equiv D_j G_i\quad\text{ on }B.
\end{equation}
Fix a connected $\Omega\subset B\setminus\CK_w$ and let $\gamma:[0,1]\to\Omega$ be a $C^1$ loop. By \cite[Lemma 2.4]{LucheseSchauder}, there exist a $C^1$ function $a:[0,1]\to\R^k$ and a continuous function $A:[0,1]\to\R^{n\times k}$ such that 
\begin{equation}\label{eq: selection along loop}
j_1w(\gamma(t))=\mb{(a(t),A(t))}+\mb{(-a(t),-A(t))}\quad\text{ and }\quad\dot{a}(t)=A(t)\dot{\gamma}(t),\quad\text{ for all }t\in[0,1].
\end{equation}
Furthermore, we have
\begin{equation}\label{eq: relation between A and G}
A(t)\in\{G(\gamma(t)),-G(\gamma(t))\},\quad\text{ for all }t\in[0,1].
\end{equation}
We claim that we have
\begin{equation}\label{eq: claim for loops}
(a(0),A(0))=(a(1),A(1)).
\end{equation}
To see this, suppose for a contradiction that $(a(0),A(0))=(-a(1),-A(1))$. Then by \eqref{eq: selection of 1-jet away from zeros of G}, the set $\gamma([0,1])\cap\{G=0\}$ must be non-empty. After reparametrising, we can thus assume that $G(\gamma(0))=0$ and $a(0)=-a(1)\neq0$. Set $T:=\{t\in[0,1]:G(\gamma(t))=0\}$. Then \eqref{eq: selection along loop} implies $\dot{a}\equiv0$ on $T$, and thus $\H^1(a(T))=0$. In particular, there exists $\xi\in\R^k$ such that $\xi\cdot a(t)\neq0$, for all $t\in T$. Let $f(t):=\xi\cdot a(t)$. Then $f\neq0$ on $T$ and $f(1)=-f(0)\neq 0$. Hence, the set $\Sigma:=\{f=0\}$ is a compact subset of $(0,1)\setminus T$. In particular, there exists a finite collection of pairwise disjoint intervals $\{(\al_j,\beta_j)\}_{j=1}^J$ such that $I_j:=[\al_j,\beta_j]\subset(0,1)\setminus T$ and $\Sigma\subset\cup_{j=1}^J(\al_j,\beta_j)$. By \eqref{eq: relation between A and G}, for each $j\in\{1,\dots,J\}$ there exists $\sigma_j\in\{-1,1\}$ such that 
\begin{equation}\label{eq: sign for A on I_j}
A\equiv\sigma_j G\circ\gamma\quad\text{ on }I_j.
\end{equation}
Set $\eps:=(\inf_{[0,1]\setminus\cup_{j=1}^JI_j}|f|)/2>0$ and let $\chi\in C^\infty_c(-\eps,\eps)$ be even and non-negative such that $\int_\R\chi=1$. Define $F:B\to\R$ by setting
$$
F(x):=\chi(\xi\cdot \nu),\quad\text{ where }\nu\in\R^k\text{ is such that }w(x)=\mb{\nu}+\mb{-\nu}.
$$
Note that $F$ is well-defined. We claim that $F$ is $C^1$ in $B$. To see this, set $\CZ:=\{w=2\mb{0}\}$ and note that, by \eqref{eq: selection of 1-jet away from zeros of G}, for every $x\in B\setminus\CZ$ and $\nu\in\R^k$ such that $w(x)=\mb{\nu}+\mb{-\nu}$, we have
\begin{equation}\label{eq: dF}
dF(x)=\pm\chi'(\xi\cdot \nu)\sum_{i=1}^n(\xi\cdot G_i)dx_i.
\end{equation}
Furthermore, if $x\in\CZ$, then, since $\chi'(0)=0$, we have $|F(x+h)-F(x)|=O(|h|^2)$ for all $h\in\R^n$ sufficiently small. In particular, $dF(x)=0$ and it follows from \eqref{eq: dF} that $F$ is $C^1$. We now define $\eta:=\sum_{i=1}^n(\xi\cdot G_i)dx_i$ and set $\omega:=F\eta\in C^1$. By \eqref{eq: symmetry of G}, $\eta$ is closed and thus \eqref{eq: dF} gives
$$
d\omega(x)=dF(x)\wedge\eta(x)+F(x)d\eta(x)=\pm\chi'(\xi\cdot\nu)\eta(x)\wedge\eta(x)=0,\quad\text{ for all }x\in B.
$$
Hence, $\omega$ is a closed $1$-form on $B$. Since $B$ is simply connected, it follows from the Poincaré lemma that $\omega$ is exact. Set $g(s):=\int_0^s\chi(\tau)\,\d\tau$, and note that $g\equiv 1/2$ on $[\eps,\infty)$ and $g\equiv-1/2$ on $(-\infty,-\eps]$. By the choice of $\eps$, we have that $g(f(\beta_j))-g(f(\al_j))\in\{-1,0,1\}$ for all $j$. Furthermore, after relabelling, $\beta_j<\al_{j+1}$ and $g(f(\beta_j))=g(f(\al_{j+1}))$ for all $j\in\{1,\dots,J-1\}$. In particular, it follows from \eqref{eq: sign for A on I_j} that
\begin{align*}
\int_\gamma\omega&=\sum_{1\leq j\leq J}\sigma_j\int_{I_j}\chi(f)f'=\sum_{1\leq j\leq J}\sigma_j(g(f(\beta_j))-g(f(\al_j)))\\
&\equiv\sum_{1\leq j\leq J}(g(f(\beta_j))-g(f(\al_j)))\equiv g(f(1))-g(f(0))\equiv 1\quad\quad\text{ (mod 2)},
\end{align*}
which, by Stokes' theorem, contradicts the exactness of $\omega$. We have thus established \eqref{eq: claim for loops}, and can define a $C^1$ function $u_\Omega:\Omega\to\R^k$ such that $w\equiv \mb{u_\Omega}+\mb{-u_\Omega}$ in $\Omega$. Finally, define $u:B\to\R^k$ by
$$
u(x):=
\begin{cases}
u_\Omega(x),&\text{ if }x\in\Omega\text{ for some connected component }\Omega\subset B\setminus\CK_w,\\
0,&\text{ if }x\in\CK_w.
\end{cases}
$$
Then by the above, $u\in C^1$ and $w\equiv\mb{u}+\mb{-u}$ in $B$. This completes the proof.
\end{proof}

\section{Properties of Multi-Valued Functions}\label{app: Multi-Valued Functions}

\subsection{Proof of Lemma \ref{lem: well-defined derivative}}

In this section, we give an elementary proof of Lemma \ref{lem: well-defined derivative}. We will use the following lemma:
\begin{lem}\label{lem: well-defined derivative preliminary}
Fix $Q\geq 1$. Let $w^1,\dots,w^Q\in\R^k$ and for each $j\in\{1,2\}$, let $L^j_1,\dots,L^j_Q:\R^n\to\R^k$ be a collection of linear maps. Suppose that for every $\xi\in S^{n-1}$ there exists a permutation $\sigma_\xi\in S_Q$ such that
\begin{equation*}
(w^i,L^1_i\cdot\xi)=(w^{\sigma_\xi(i)},L^2_{\sigma_\xi(i)}\cdot\xi),\quad\text{ for all }i\in\{1,\dots,Q\}.
\end{equation*} 
Then, there exists a permutation $\sigma\in S_Q$ such that $(w^i,L^1_i)=(w^{\sigma(i)},L^2_{\sigma(i)})$, for all $i\in\{1,\dots,Q\}$.
\end{lem}
\begin{proof}
Fix $\nu\in\R^k$, and for $j\in\{1,2\}$, define the polynomial $P^j_\nu:\R\times\R^n\to\R$ by setting
$$
P^j_\nu(t,z):=\prod_{i=1}^Q(t-\nu\cdot (w^i+L^j_i\cdot z)).
$$
By assumption, $P^1_\nu=P^2_\nu$ in $\R[t,z_1,\dots,z_{n}]$, for all $\nu\in\R^k$. Hence, by the unique factorisation of $\R[t,z_1,\dots,z_{n}]$, for every $\nu\in\R^k$ there exists a permutation $\sigma_\nu\in S_Q$ such that $\nu\cdot (w^i+L^1_i)\equiv\nu\cdot (w^{\sigma_\nu(i)}+L^2_{\sigma_\nu(i)})$ for all $i\in\{1,\dots,Q\}$. For each $\sigma\in S_Q$, let $V_\sigma$ denote the subspace
$$
V_\sigma:=\{\nu\in\R^k:\nu\cdot (w^i+L^1_i)\equiv\nu\cdot (w^{\sigma(i)}+L^2_{\sigma(i)}),\text{ for all }i\in\{1,\dots,Q\}\}.
$$
Then by the above, $\R^k$ decomposes as the union $\R^k=\cup_{\sigma\in S_Q} V_\sigma$. In particular, there exists $\sigma\in S_Q$ such that $V_\sigma=\R^k$. This completes the proof.
\end{proof}
We now prove Lemma \ref{lem: well-defined derivative}:
\begin{proof}[Proof of Lemma \ref{lem: well-defined derivative}]
Assume without loss of generality that $Q\geq 2$. Fix $x\in B^n_1(0)$ and suppose $w:B^n_1(0)\to\A_Q(\R^k)$ is a $Q$-valued function and that for each $j\in\{1,2\}$ there exist $Q$ linear maps $L^j_1,\dots,L^j_Q:\R^n\to\R^k$ such that
\begin{equation}\label{eq: proof of well-defined derivative}
\lim_{y\to x}|x-y|^{-1}\G\biggl(w(y),\sum_{i=1}^Q\mb{w^i(x)+L^j_i\cdot(y-x)}\biggr)=0.
\end{equation}
By Lemma \ref{lem: well-defined derivative preliminary}, it suffices to show that for each $\xi\in S^{n-1}$, there exists $\sigma_\xi\in S_Q$ such that 
\begin{equation}\label{eq: proof of well-defined derivative-1}
(w^i(x),L^1_i\cdot\xi)=(w^{\sigma_\xi(i)}(x),L^2_{\sigma_\xi(i)}\cdot\xi),\quad\text{ for all }i\in\{1,\dots,Q\}.
\end{equation} 
So fix $\xi\in S^{n-1}$. For $t\in(0,1-|x|)$, define $y_t:=x+t\xi$, and for each $j\in\{1,2\}$ let $\sigma^j_t\in S_Q$ be a permutation such that
$$
\G\biggl(w(y_t),\sum_{i=1}^Q\mb{w^i(x)+tL^j_i\cdot \xi}\biggr)^2=\sum_{i=1}^Q|w^{\sigma^j_t(i)}(y_t)-(w^i(x)+tL^j_i\cdot\xi)|^2.
$$
For each $\ell\geq 1$ sufficiently large, let $t_\ell:=\ell^{-1}$ and, for $j\in\{1,2\}$, set $\sigma^j_\ell:=\sigma^j_{t_\ell}$. By \eqref{eq: proof of well-defined derivative}, provided $\ell$ is sufficiently large, we must have $w^{(\sigma^1_\ell)^{-1}(i)}(x)=w^{(\sigma^2_\ell)^{-1}(i)}(x)$, for all $i\in\{1,\dots,Q\}$. Taking $\sigma_\ell:=(\sigma^2_\ell)^{-1}\circ\sigma^1_\ell$, it follows from \eqref{eq: proof of well-defined derivative} that we have
\begin{equation}\label{eq: proof of well-defined derivative 1}
\lim_{\ell\to\infty}\sum_{i=1}^Q|L^1_i\cdot\xi-L^2_{\sigma_\ell(i)}\cdot\xi|=0.
\end{equation}
After taking a subsequence and relabelling, there exist permutations $\sigma^1_\xi,\sigma^2_\xi\in S_Q$ such that $\sigma^j_\ell=\sigma^j_\xi$, for all $\ell\geq 1$ and each $j\in\{1,2\}$. Setting $\sigma_\xi:=(\sigma^2_\xi)^{-1}\circ\sigma^1_\xi$, we see that \eqref{eq: proof of well-defined derivative-1} follows from \eqref{eq: proof of well-defined derivative 1}. This establishes \eqref{eq: proof of well-defined derivative-1}, and completes the proof.
\end{proof}

\subsection{Proof of Proposition \ref{prop: W^2,2 + Dw cont implies C^1}}
In this section, we prove Proposition \ref{prop: W^2,2 + Dw cont implies C^1}.
\begin{proof}[Proof of Proposition \ref{prop: W^2,2 + Dw cont implies C^1}]
Suppose $w\in W^{2,2}_{\textnormal{loc}}(B^n_1(0);\A_Q(\R^k))$ and $Dw\in C^0(B^n_1(0);\A_Q(\R^{n\times k}))$. Then $J:=j_1w\in W^{1,2}_{\textnormal{loc}}(B^n_1(0);\A_Q(\R^k\times\R^{n\times k}))$. By \cite[Proposition 2.14]{DLS11Qvalued}, $w$ has a locally Lipschitz representative (which we also denote by $w$). We first show that $J$ has a continuous representative. Let $G$ be a continuous representative of $Dw$. Fix $x\in B^n_1(0)$ and write $w(x)=\sum_{i=1}^Q\mb{a_i}$ and $G(x)=\sum_{i=1}^Q\mb{A_i}$ for some $a_1,\dots,a_Q\in\R^k$ and $A_1,\dots,A_Q\in\R^{n\times k}$. Let $\Sigma_x$ denote the set
$$
\Sigma_x:=\biggl\{\sum_{i=1}^Q\mb{(a_i,A_{\sigma(i)})}:\sigma\in S_Q\biggr\}.
$$
Note that we have
\begin{equation}\label{eq: proof of W^2,2 + Dw cont implies C^1 distance estimate}
\dist(J(y),\Sigma_x)^2\leq \G(w(x),w(y))^2+\G(G(x),G(y))^2,\quad\text{ for }\H^n\text{-a.e.} y\in B^n_1(0).
\end{equation}
Let $M_x:(0,1-|x|)\to\R$ be defined by
$$
M_x(r)^2:=\sup_{y\in B^n_r(x)}(\G(w(x),w(y))^2+\G(G(x),G(y))^2).
$$
Let $N\in\{1,\dots,Q!\}$ and $T_1,\dots,T_N\in\A_Q(\R^k\times\R^{n\times k})$ such that $\Sigma_x=\{T_1,\dots,T_N\}$. If $N=1$, set $\widehat{J}(x):=T_1$. If $N\geq 2$, set $\eps:=\inf_{i\neq j}\G(T_i,T_j)/4>0$. Since $w$ and $G$ are continuous, it follows from \eqref{eq: proof of W^2,2 + Dw cont implies C^1 distance estimate} that there exists $r_x\in(0,1-|x|)$ such that 
\begin{equation}\label{eq: proof of W^2,2 + Dw cont implies C^1 containment}
\min_{i=1,\dots,N}\G(J(y),T_i)\leq\eps,\quad\text{ for }\H^n\text{-a.e. }y\in B^n_{r_x}(x).
\end{equation}
For each $i\in\{1,\dots,N\}$, let $\chi_i\in C^\infty_c(-2\eps,2\eps)$ such that $0\leq\chi_i\leq 1$ and $\chi_i\equiv 1$ in a neighbourhood of $(-\eps,\eps)$. Set $F_i:=\chi_i\circ \G(J,T_i)\in W^{1,2}_{\textnormal{loc}}(B^n_{r_x}(x))$. Note that $F_i(y)\in\{0,1\}$ for $\H^n$-a.e. $y\in B^n_{r_x}(x)$, and thus each $F_i$ is constant. By \eqref{eq: proof of W^2,2 + Dw cont implies C^1 containment}, we have that $\sum_{i=1}^NF_i(y)=1$ for $\H^n$-a.e. $y\in B^n_{r_x}(x)$. Thus, after relabelling, we have that $F_1(y)=1$ and $F_j(y)=0$ for $\H^n$-a.e. $y\in B^n_{r_x}(x)$ and all $j\in\{2,\dots,N\}$. Setting $\widehat{J}(x):=T_1$, it now follows from \eqref{eq: proof of W^2,2 + Dw cont implies C^1 distance estimate} that for any $r\in(0,r_x)$ we have
\begin{equation}\label{eq: proof of W^2,2 + Dw cont implies C^1 rep of J cont estimate}
\G(J(y),\widehat{J}(x))\leq M_x(r),\quad\text{ for }\H^n\text{-a.e. }y\in B^n_r(x).
\end{equation}
Since $J(y)=\widehat{J}(y)$ for $\H^n$-a.e. $y\in B^n_1(0)$, it now follows from \eqref{eq: proof of W^2,2 + Dw cont implies C^1 rep of J cont estimate} that
$$
\G(\widehat{J}(y),\widehat{J}(x))\leq 2M_x(r),\quad\text{ for all }y\in B^n_r(x),\text{ for all }r\in(0,r_x).
$$
In particular, $\widehat{J}$ is continuous. It only remains to show that $w$ is differentiable with derivative equal to $G$. Fix a ball $B^n_{2R}(x_0)\Subset B^n_1(0)$ and let $\omega$ denote the modulus of continuity of $\widehat{J}$ in $B^n_{2R}(x_0)$. Let $h\in \R^n$ be such that $0<|h|<R/2$. By Fubini, \cite[Proposition~1.2]{DLS11Qvalued} and the definition of $J$, for $\H^{n-1}$-a.e $x\in (x_0+h^\perp)\cap B^n_{R/2}(x_0)=:I_h$, there exist absolutely continuous functions $a_1,\dots,a_Q:[0,1]\to\R^k$ and $A_1,\dots,A_Q:[0,1]\to\R^{n\times k}$ such that
$$
\widehat{J}(x+th)=\sum_{\ell=1}^Q\mb{(a_\ell(t),A_\ell(t))}\quad\text{ and }\quad \dot{a}_\ell(t)=A_\ell(t)h,\quad\text{ for a.e. }t\in[0,1].
$$
Furthermore, for any $t\in[0,1]$, we have
$$
A_\ell(s)\in\bigcup_{j=1}^Q\overline{B^{n\times k}_{\omega(|h|)}(A_j(0))},\quad\text{ for all }s\in[0,t], \ell\in\{1,\dots,Q\}.
$$
Since $A_\ell$ is continuous, the set $A_\ell([0,t])$ is contained in the connected component of the above union containing $A_\ell(0)$. Now, a connected component of a union of $Q$ balls of radius $\rho$ has diameter at most $2Q\rho$, and consequently
$$
|A_\ell(t)-A_\ell(0)|\leq 2Q \omega(|h|),\quad\text{ for all }t\in[0,1], \ell\in\{1,\dots,Q\}.
$$
Integrating and summing, we have
\begin{equation}\label{eq: C^1,alpha decay estimate}
\G\bigl(w(x+h),\widehat{J}(x)[h]\bigr)\leq C \omega(|h|)|h|,\quad\text{ for }\H^{n-1}\text{-a.e. }x\in I_h,
\end{equation}
where $C\in(0,\infty)$ is independent of $h$. Since $w$ and $\widehat{J}$ are continuous, \eqref{eq: C^1,alpha decay estimate} holds with $x_0$ in place of $x$. It now follows that $w$ is differentiable, with derivative equal to $G$. This completes the proof.
\end{proof}

\subsection{Proof of Proposition \ref{prop: current of Q-valued W^1,2 function}}
In this section, we prove Proposition \ref{prop: current of Q-valued W^1,2 function}. 
\begin{proof}[Proof of Proposition \ref{prop: current of Q-valued W^1,2 function}]
Suppose without loss of generality that $f\in W^{1,r}(B^n_1(0);\A_Q(\R^k))$. Let $e_1,\dots,e_n$ denote the standard basis of $\R^n$. Define $\zeta:\R^{n\times k}\to\Lambda^n(\R^{n+k})$ by setting
$$
\zeta(P):=(Pe_1,e_1)\wedge\cdots\wedge(Pe_n,e_n).
$$
Then $|\zeta(P)|^2=\det(I_n+P^TP)$, for all $P\in\R^{n\times k}$. Fix $P\in\R^{n\times k}$. Note that $\textnormal{ker}(P^TP)=\textnormal{ker}(P)$, and it follows from the rank-nullity theorem that $\textnormal{rank}(P^TP)=\textnormal{rank}(P)\leq \min\{n,k\}=:r$. Furthermore, $P^TP$ is symmetric and positive semi-definite and thus has $n$ non-negative eigenvalues $\lambda_1\geq\cdots\geq\lambda_n$. By the above, we must have that $\lambda_i=0$ for all $i\geq r+1$. In particular, the AM-GM inequality implies that
$$
|\zeta(P)|=\biggl(\prod_{i=1}^r(1+\lambda_i)\biggr)^{1/2}
\leq r^{-r/2}\biggl(\sum_{i=1}^r(1+\lambda_i)\biggr)^{r/2}
=r^{-r/2}(r+\textnormal{tr}(P^TP))^{r/2}=r^{-r/2}(r+|P|^2)^{r/2}.
$$
It now follows that for any $\omega\in\CD^n$, we have
\begin{align}
\nonumber|\bf{T}_f(\omega)|&\leq\int_{B^n_1(0)}\sum_{\ell=1}^Q|\omega(f^\ell(x),x)||\zeta(Df^\ell(x))|\leq r^{-r/2}Q\sup_{\R^k\times B^n_1(0)}|\omega|\int_{B^n_1(0)}(r+|Df|^2)^{r/2}\\
&\leq 2^{r/2}r^{-r/2}Q\sup_{\R^k\times B^n_1(0)}|\omega|\int_{B^n_1(0)}(r^{r/2}+|Df|^r).\label{eq: mass-bound for current of multi-valued sobolev function}
\end{align}
Thus, $\bf{T}_f$ is well-defined and $\mathbb{M}(\bf{T}_f)<\infty$. By \cite[Proposition 2.5]{DLS11Qvalued}, for every $\lambda>0$ there exists a Lipschitz function $f_\lambda:B^n_1(0)\to\A_Q(\R^k)$ such that $\Lip(f_\lambda)\leq\lambda$ and
$$
\H^n(\{f\neq f_\lambda\})=o(\lambda^{-r}),\quad\text{ as }\lambda\uparrow\infty.
$$
Note that for each $\ell$, we have $Df^\ell=Df^\ell_\lambda$ $\H^n$-a.e. on $E_\lambda:=\{f=f_\lambda\}$. In particular, for any $\omega\in\CD^n$, \eqref{eq: mass-bound for current of multi-valued sobolev function} implies that
$$
|\bf{T}_f(\omega)-\bf{T}_{f_\lambda}(\omega)|\leq Q\sup_{\R^k\times B^n_1(0)}|\omega|\int_{B^n_1(0)\setminus E_\lambda}(r^{r/2}+|Df|^r+|Df_\lambda|^r)\longrightarrow0\quad\text{ as }\lambda\uparrow\infty.
$$
By \cite[Proposition 1.4 \& Theorem 2.1]{DLS15MultiVldFcts}, each $\bf{T}_{f_\lambda}$ is an $n$-dimensional integer rectifiable current with $\partial \bf{T}_\lambda=0$ in $\R^k\times B^n_1(0)$. It thus follows from Federer-Fleming's compactness theorem \cite[4.2.17]{Fed69} that $\bf{T}_f$ is integer rectifiable with $\partial\bf{T}_f=0$ in $\R^k\times B^n_1(0)$. This completes the proof.
\end{proof}

\section{The Obstruction to Calibrations in Higher Dimensions}\label{app: Obstruction to Calibrations in Higher Dimensions}
In this section, we discuss the obstruction to generalising Theorem \ref{thm: graph of derivative is calibrated in 2d} to higher dimensions. We first record the following result (cf.~Remark \ref{remark: Obstruction to Calibrations in Higher Dimensions}):
\begin{prop}\label{prop: Obstruction to a Universal Calibration in Higher Dimensions}
Fix $n\geq 3$. Suppose $h$ is a metric on $\R^n\times B^n_1(0)$ and $\omega$ is a semi-calibration on $(\R^n\times B^n_1(0),h)$. Then there exists a smooth function $u:B^n_1(0)\to\R$ solving the minimal surface equation such that $\omega$ does not semi-calibrate $\bf{T}_{Du}$ in $(\R^n\times B^n_1(0),h)$.
\end{prop}
\begin{proof}
It suffices to prove the result when $n=3$. Suppose that there exists a metric $h$ on $\R^3\times B^3_1(0)$ and a semi-calibration $\omega$ on $(\R^3\times B^3_1(0),h)$ which contradicts the claim. Let $A$ be the $3\times 3$ diagonal matrix defined by $A:=\textnormal{diag}(1,1,-2)$. By Cauchy-Kovalevskaya, there exists $\delta>0$ and a smooth function $v:B^3_\delta(0)\to\R$ solving the minimal surface equation such that $Dv(0)=0$ and $D^2v(0)=A$. For $t\in(0,\delta)$, define $u_t:B^3_{\delta/t}(0)\to\R$ by setting $u_t(x):=v(tx)/t$. Then each $u_t$ solves the minimal surface equation and satisfies $Du_t(0)=0$ and $D^2u_t(0)=tA$. Let $E_1,E_2,E_3,e_1,e_2,e_3$ denote the standard basis for $\R^3\times B^3_1(0)$. For $t\in \R$, set $G_t:=Du_t$. Then by the above, the orientation form of $\bf{T}_{G_t}$ at $(0,0)$ is given by
$$
\zeta_t:= (tE_1+e_1)\wedge(tE_2+e_2)\wedge(-2tE_3+e_3)=Z_0+tZ_1+t^2Z_2+t^3Z_3,
$$
where $Z_3=-2E_1\wedge E_2\wedge E_3$ is non-zero. In particular, the polynomial $q(t):=|\zeta_t|^2_h$ is a degree six polynomial in $t$. Let $p$ be the polynomial uniquely defined by
$$
p(t):=\langle\omega,\zeta_t\rangle=\langle\omega,Z_0\rangle+t\langle\omega,Z_1\rangle+t^2\langle\omega,Z_2\rangle+t^3\langle\omega,Z_3\rangle.
$$
By assumption, $\omega$ semi-calibrates $\bf{T}_{G_t}$, and thus we must have
\begin{equation}\label{eq: equality of polynomials}
p(t)^2=q(t),\quad\text{ for all }t\in\R.
\end{equation}
In particular, $\langle\omega,Z_3\rangle^2=|Z_3|^2_h\neq 0$, and thus $p$ is a non-vanishing cubic. In particular, $p$ has a real root, and thus \eqref{eq: equality of polynomials} implies that $q$ also has a real root $t_0\in\R$. However, this implies that $\zeta_{t_0}=0$, which is a contradiction since the pushforward of $\zeta_{t_0}$ under the orthogonal projection $\R^3\times \R^3\to\{0\}^3\times\R^3$ is $e_1\wedge e_2\wedge e_3$ which is non-zero. This completes the proof.
\end{proof}

As a consequence of the proof of Proposition \ref{prop: Obstruction to a Universal Calibration in Higher Dimensions}, we obtain the following:
\begin{prop}
For every $n\geq 3$ and $Q\geq 4$ there exists $u\in C^\infty(B^n_1(0);\A_Q(\R))$ such that the integral varifold $\bf{v}(u)$ is stationary in $\R\times B^n_1(0)$ and $\bf{T}_{Du}$ is not semi-calibrated in any smooth metric on $\R^n\times B^n_1(0)$.
\end{prop}
\begin{proof}
It suffices to prove the result when $n=3$. For $t\in(0,\delta)$ let $u_t:B^n_{\delta/t}(0)\to\R$ be as in the proof of Proposition \ref{prop: Obstruction to a Universal Calibration in Higher Dimensions}. Fix $0<t_1<t_2<t_3<t_4<\delta$, and set $u:=\sum_{j=1}^4\mb{u_{t_j}}$. Suppose that $h$ is a metric on $\R^3\times B^3_1(0)$ and $\omega$ is a calibration on $(\R^3\times B^3_1(0),h)$ which calibrates $\bf{T}_{Du}$. Let $p$ and $q$ be polynomials defined as in the proof of Proposition \ref{prop: Obstruction to a Universal Calibration in Higher Dimensions}. Arguing as in \eqref{eq: equality of polynomials}, the polynomial $f(t):=q(t)-p(t)^2$ has at least four distinct real roots. Furthermore, by the proof of Proposition \ref{prop: Obstruction to a Universal Calibration in Higher Dimensions}, $f$ is non-vanishing and of degree at most six. Also, since $\omega$ is a semi-calibration, $f$ is non-negative. Thus, each root of $f$ must have multiplicity at least $2$. Since there are at least four distinct such roots, $f$ must vanish identically. This is a contradiction, and completes the proof.
\end{proof}


\begin{thebibliography}{10}

\bibitem{Aiex26C^1Multi-Valued}
Nicolau~S. Aiex.
\newblock On {A}lmgren's {M}ultiple-valued {F}unctions of {C}lass 1.
\newblock {\em \href{https://arxiv.org/abs/2606.13002}{arXiv:2606.13002}},
  2026.

\bibitem{Almgren00}
Frederick~J Almgren.
\newblock {\em Almgren's big regularity paper}, volume \textbf{1} of {\em World
  Scientific Monograph Series in Mathematics}.
\newblock World Scientific, Singapore, 2000.

\bibitem{AC04ThinObstacle}
Ioannis Athanasopoulos and Luis~A. Caffarelli.
\newblock Optimal regularity of lower dimensional obstacle problems.
\newblock {\em Zap. Nauchn. Sem. S.-Peterburg. Otdel. Mat. Inst. Steklov.
  (POMI)}, \textbf{310}:49--66, 2004.

\bibitem{BK17}
Spencer~T. Becker-Kahn.
\newblock Transverse singularities of minimal two-valued graphs in arbitrary
  codimension.
\newblock {\em J. Differ. Geom.}, \textbf{107}(2):241--325, 2017.

\bibitem{Bel25}
Costante Bellettini.
\newblock Extensions of {S}choen-{S}imon-{Y}au and {S}choen-{S}imon theorems
  via iteration \`a{} la {D}e {G}iorgi.
\newblock {\em Invent. Math.}, \textbf{240}(1):1--34, 2025.

\bibitem{DLMSVDirRectifiability}
Camillo De~Lellis, Andrea Marchese, Emanuele~Nunzio Spadaro, and Daniele
  Valtorta.
\newblock Rectifiability and upper {M}inkowski bounds for singularities of
  harmonic {$Q$}-valued maps.
\newblock {\em Comment. Math. Helv.}, \textbf{93}(4):737--779, 2018.

\bibitem{DLS11Qvalued}
Camillo De~Lellis and Emanuele~Nunzio Spadaro.
\newblock Q-valued functions revisited.
\newblock {\em Memoirs of the American Mathematical Society},
  \textbf{211}(991):vi+79, 2011.

\bibitem{DLS14GradientL^pEstimates}
Camillo De~Lellis and Emanuele~Nunzio Spadaro.
\newblock Regularity of area minimizing currents {I}: gradient ${L}^{p}$
  estimates.
\newblock {\em Geometric and Functional Analysis}, \textbf{24}(6), 2014.

\bibitem{DLS15MultiVldFcts}
Camillo De~Lellis and Emanuele~Nunzio Spadaro.
\newblock Multiple valued functions and integral currents.
\newblock {\em Annali della Scuola Normale Superiore di Pisa. Classe di
  Scienze}, \textbf{14}(4):1239--1269, 2015.

\bibitem{DLSS20}
Camillo De~Lellis, Emanuele~Nunzio Spadaro, and Luca Spolaor.
\newblock Regularity theory for 2-dimensional almost minimal currents {III}:
  blowup.
\newblock {\em J. Differ. Geom.}, \textbf{116}(1):125--185, 2020.

\bibitem{DPSV25}
Guido De~Philippis, Luca Spolaor, and Bozhidar Velichkov.
\newblock ({Q}uasi-)conformal methods in two-dimensional free boundary
  problems.
\newblock {\em J. Eur. Math. Soc.}, \textbf{27}:3369--3406, 2025.

\bibitem{EellsSampson64Harmonic}
James Eells and Joseph~H. Sampson.
\newblock Harmonic mappings of {R}iemannian manifolds.
\newblock {\em American Journal of Mathematics}, \textbf{86}(1):109--160, 1964.

\bibitem{Evans10PDE}
Lawrence~C. Evans.
\newblock {\em Partial Differential Equations}, volume \textbf{19} of {\em
  Graduate Studies in Mathematics}.
\newblock American Mathematical Society, Providence, Rhode Island, 2nd edition,
  2010.

\bibitem{Fed69}
Herbert Federer.
\newblock {\em Geometric measure theory}, volume \textbf{153} of {\em Die
  Grundlehren der mathematischen Wissenschaften}.
\newblock Springer-Verlag New York, Inc., New York, 1969.

\bibitem{FR22TOP}
Xavier Fernández-Real.
\newblock The thin obstacle problem: a survey.
\newblock {\em Publ. Mat.}, \textbf{66}(1):3--55, 2022.

\bibitem{FS18ThinObstacle}
Matteo Focardi and Emanuele~Nunzio Spadaro.
\newblock On the measure and the structure of the free boundary of the lower
  dimensional obstacle problem.
\newblock {\em Arch. Ration. Mech. Anal.}, \textbf{230}(1):125--184, 2018.

\bibitem{FS20ThinObstacle}
Matteo Focardi and Emanuele~Nunzio Spadaro.
\newblock How a minimal surface leaves a thin obstacle.
\newblock {\em Annales de l'Institut Henri Poincar{\'e} C, Analyse non
  linéaire}, \textbf{37}(4):1017--1046, 2020.

\bibitem{FS24ThinObstacle}
Federico Franceschini and Joaquim Serra.
\newblock Free boundary partial regularity in the thin obstacle problem.
\newblock {\em Communications on Pure and Applied Mathematics},
  \textbf{77}(1):630--669, 2024.

\bibitem{HartmanWintner1953}
Philip Hartman and Aurel Wintner.
\newblock On the {L}ocal {B}ehaviour of {S}olutions of {N}on-{P}arabolic
  {P}artial {D}ifferential {E}quations.
\newblock {\em American Journal of Mathematics}, \textbf{75}(3):449--476, 1953.

\bibitem{Henneberg1875}
Lebrecht Henneberg.
\newblock {\em Über solche Minimalflächen, welche eine vorgeschriebene ebene
  Curve zur geodätischen Linie haben}.
\newblock Doctoral dissertation, Eidgenössisches Polytechnikum Zürich,
  Z{\"{u}}rich, Switzerland, 1875.

\bibitem{Kru14ConstFreq}
Brian Krummel.
\newblock Constant frequency and the higher regularity of branch sets.
\newblock {\em \href{https://arxiv.org/abs/1410.7339}{arXiv:1410.7339}}, 2014.

\bibitem{KMW}
Brian Krummel, Paul Minter, and Neshan Wickramasekera.
\newblock A branch set stratification for stationary varifolds with
  epsilon-regularity.
\newblock {\em \href{https://arxiv.org/abs/2606.01511}{arXiv:2606.01511}},
  2026.

\bibitem{KW13}
Brian Krummel and Neshan Wickramasekera.
\newblock Fine properties of branch point singularities: Two-valued harmonic
  functions.
\newblock {\em \href{https://arxiv.org/abs/1311.0923}{arXiv:1311.0923}}, 2013.

\bibitem{KW17}
Brian Krummel and Neshan Wickramasekera.
\newblock Fine properties of branch point singularities: Dirichlet energy
  minimizing multi-valued functions.
\newblock {\em \href{https://arxiv.org/abs/1711.06222}{arXiv:1711.06222}},
  2017.

\bibitem{KW21}
Brian Krummel and Neshan Wickramasekera.
\newblock Fine properties of branch point singularities: stationary two-valued
  graphs and stable minimal hypersurfaces near points of density $<3$.
\newblock {\em \href{https://arxiv.org/abs/2111.12246}{arXiv:2111.12246}},
  2021.

\bibitem{Lawson77lectures}
H.~Blaine Lawson, Jr.
\newblock {\em Lectures on {M}inimal {S}ubmanifolds. {V}ol. {I}}, volume
  \textbf{14} of {\em Monografías de Matemática}.
\newblock Instituto de Matemática Pura e Aplicada, Rio de Janeiro, 1977.
\newblock MR 527121.

\bibitem{Lewy72CoincidenceSet}
Hans Lewy.
\newblock On the coincidence set in variational inequalities.
\newblock {\em J. Differ. Geom.}, \textbf{6}(4):497--501, 1972.

\bibitem{Luch26SmallLap}
Mattia Luchese.
\newblock ${C}^{1,\al}$ {R}egularity of ${W}^{2,2}$ ${Q}$-{V}alued {M}aps with
  {S}mall {L}aplacian.
\newblock {\em \textnormal{Preprint}}, 2026.

\bibitem{LucheseSchauder}
Mattia Luchese.
\newblock Perturbative {S}chauder {E}stimates for {Q}-{V}alued {Q}uasilinear
  {E}lliptic {S}ystems and a {S}harp {D}imension {B}ound for {B}ranch {S}ets of
  {S}tationary {G}raphs.
\newblock {\em \href{https://arxiv.org/abs/2609.12950v1}{arXiv:2609.12950}},
  2026.

\bibitem{MinterCampanato}
Paul Minter.
\newblock A {C}ampanato regularity theory for multi-valued functions with
  applications to minimal surface regularity theory.
\newblock {\em Journal of Functional Analysis}, \textbf{289}(1), 2025.

\bibitem{MS26}
Paul Minter and Sidney Stanbury.
\newblock Area minimising hypersurfaces mod $p$ do not admit immersed branch
  points.
\newblock {\em \href{https://arxiv.org/abs/2603.03100}{arXiv:2603.03100}},
  2026.

\bibitem{Mor66MultipleIntegrals}
Charles~Bradfield Morrey.
\newblock {\em Multiple Integrals in the Calculus of Variations}.
\newblock Springer, 1966.

\bibitem{MoyaPerez23MinSurfaces}
David Moya and Joaquín Pérez.
\newblock Generalized {H}enneberg {S}table {M}inimal {S}urfaces.
\newblock {\em Results in Mathematics}, \textbf{78}(2):53, 2023.

\bibitem{Osserman86Survey}
Robert Osserman.
\newblock {\em A {S}urvey of {M}inimal {S}urfaces}.
\newblock Dover Publications, New York, second edition, 1986.

\bibitem{Simon93Cylindrical}
Leon Simon.
\newblock Cylindrical tangent cones and the singular set of minimal
  submanifolds.
\newblock {\em J. Differ. Geom.}, \textbf{38}(3):585--652, 1993.

\bibitem{Simon96}
Leon Simon.
\newblock {\em Theorems on Regularity and Singularity of Energy Minimizing
  Maps}.
\newblock Lectures in Mathematics. ETH Zürich. Birkhäuser Basel, 1996.

\bibitem{SW07}
Leon Simon and Neshan Wickramasekera.
\newblock Stable branched minimal immersions with prescribed boundary.
\newblock {\em J. Differ. Geom.}, \textbf{75}:143--173, 2007.

\bibitem{SW16}
Leon Simon and Neshan Wickramasekera.
\newblock A frequency function and singular set bounds for branched minimal
  immersions.
\newblock {\em Communications on Pure and Applied Mathematics},
  \textbf{69}(7):1213--1258, 2016.

\bibitem{Spadaro2010ComplexVarieties}
Emanuele~Nunzio Spadaro.
\newblock Complex varieties and higher integrability of {D}ir-minimizing
  {Q}-valued functions.
\newblock {\em manuscripta mathematica}, \textbf{132}(3), 2010.

\bibitem{WZ26StableMinHypNonBranchedSingSet}
Gaoming Wang and Xuwen Zhang.
\newblock Regularity of branched stable minimal immersed hypersurfaces.
\newblock {\em \href{https://arxiv.org/abs/2607.29632}{arXiv:2607.29632}},
  2026.

\end{thebibliography}
\end{document}